\documentclass[manuscript,screen,review=false]{acmart}

\usepackage{tabularx,colortbl}
\usepackage[capitalize]{cleveref}
\usepackage{multicol,multirow}
\usepackage{tikz}
\usepackage[linesnumbered,ruled,commentsnumbered]{algorithm2e}
\usepackage{algorithmicx}
\usepackage{algpseudocode}
\usepackage{slashbox}
\usepackage{rotating}
\usepackage{graphicx}
\usepackage{longtable}
\usepackage{tikz}
\usepackage{pgfplots}
\usepgfplotslibrary{groupplots}
\usepackage{amsmath}
\usetikzlibrary{patterns}
\usepackage{slashbox}
\usepackage{pgfplots}
\usepackage{hyperref}
\pgfplotsset{compat = 1.8}
\usetikzlibrary{arrows.meta, shapes.geometric}
\usetikzlibrary{decorations}
\usetikzlibrary{decorations.markings}
\usetikzlibrary{shapes,snakes}
\usetikzlibrary{intersections}
\usetikzlibrary{patterns}
\usepackage[frozencache=true,cachedir=minted-cache]{minted} 
\usepackage{array}
\newcolumntype{L}[1]{>{\raggedright\let\newline\\\arraybackslash\hspace{0pt}}m{#1}}
\newcolumntype{C}[1]{>{\centering\let\newline\\\arraybackslash\hspace{0pt}}m{#1}}
\newcolumntype{R}[1]{>{\raggedleft\let\newline\\\arraybackslash\hspace{0pt}}m{#1}}

\emergencystretch=\maxdimen
\newcommand{\glsolve}{\texttt{TOMLAB/glcSolve}}
\newcommand{\glbolve}{\texttt{TOMLAB/glbSolve}}
\newcommand{\glccluster}{\texttt{TOMLAB/glcCluster}}
\newcommand{\tomlab}{\texttt{TOMLAB}}

\newcommand{\direct}{\texttt{DIRECT}}
\newcommand{\directrev}{\texttt{DIRECT-rev}}
\newcommand{\adc}{\texttt{ADC}}
\newcommand{\directrest}{\texttt{DIRECT-restart}}
\newcommand{\birect}{\texttt{BIRECT}}
\newcommand{\birectgb}{\texttt{Gb-BIRECT}}
\newcommand{\birmin}{\texttt{BIRMIN}}
\newcommand{\dirmin}{\texttt{DIRMIN}}
\newcommand{\directa}{\texttt{Aggressive DIRECT}}

\newcommand{\directav}{\texttt{DIRECT-a}}
\newcommand{\directm}{\texttt{DIRECT-m}}
\newcommand{\directz}{\texttt{DIRECT-l}}
\newcommand{\directgl}{\texttt{DIRECT-GL}}
\newcommand{\directg}{\texttt{DIRECT-G}}
\newcommand{\directl}{\texttt{DIRECT-L}}
\newcommand{\plor}{\texttt{PLOR}}
\newcommand{\glbsolve}{\texttt{glbSolve}}
\newcommand{\gbdirect}{\texttt{Gb-glbSolve}}
\newcommand{\symdirec}{\texttt{glbSolve-sym}}
\newcommand{\symdirect}{\texttt{glbSolve-sym2}}
\newcommand{\disimplc}{\texttt{DISIMPL-C}}
\newcommand{\disimplv}{\texttt{DISIMPL-V}}
\newcommand{\disimplcgb}{\texttt{GB-DISIMPL-C}}
\newcommand{\disimplvgb}{\texttt{GB-DISIMPL-V}}
\newcommand{\directmr}{\texttt{MrDIRECT}}
\newcommand{\directmro}{\texttt{MrDIRECT}}

\newcommand{\directh}{\texttt{DIRECT-GLh}}
\newcommand{\directbarrier}{\texttt{DIRECT-Barrier}}
\newcommand{\directnas}{\texttt{DIRECT-NAS}}
\newcommand{\directsub}{\texttt{subDIRECT-Barrier}}

\newcommand{\directll}{\texttt{DIRECT-L1}}
\newcommand{\directce}{\texttt{DIRECT-GLce}}
\newcommand{\directcemin}{\texttt{DIRECT-GLce-min}}
\newcommand{\directc}{\texttt{DIRECT-GLc}}
\newcommand{\disimpllc}{\texttt{Lc-DISIMPL-C}}
\newcommand{\disimpllv}{\texttt{Lc-DISIMPL-V}}

\newcommand{\matlab}{\texttt{MATLAB}}
\newcommand{\directfinkel}{\texttt{DIRECT v4.0}}
\newcommand{\toolbox}{\texttt{DIRECTGO}}
\newcommand{\toolboxa}{\textbf{DIRECTGO.mltbx}}
\newcommand{\toolboxb}{\textbf{DIRECTGO.mlappinstall}}
\newcommand{\st}{\mathrm{\; s.t.}}

\newcommand{\indexsett}{\mathbb{I}}
\newcommand{\fmincon}{\texttt{fmincon}}
\newcommand{\directlib}{\texttt{DIRECTGOLib v1.0}}

\definecolor{bg}{rgb}{0.95,0.95,0.95}
\definecolor{ao}{rgb}{0.0, 0.5, 0.0}
\definecolor{arsenic}{rgb}{0.23, 0.27, 0.29}
\definecolor{armygreen}{rgb}{0.29, 0.33, 0.13}
\definecolor{antiquebrass}{rgb}{0.8, 0.58, 0.46}
\definecolor{DarkRed}{rgb}{0.55, 0.0, 0.0}
\definecolor{darkblue}{rgb}{0.0, 0.0, 0.55}
\definecolor{blueryb}{rgb}{0.01, 0.28, 1.0}
\definecolor{bluebell}{rgb}{0.64, 0.64, 0.82}
\definecolor{red}{rgb}{1.0, 0.0, 0.0}
\definecolor{redwood}{rgb}{0.67, 0.31, 0.32}
\definecolor{rose}{rgb}{1.0, 0.0, 0.5}
\definecolor{rosybrown}{rgb}{0.74, 0.56, 0.56}
\definecolor{rosewood}{rgb}{0.4, 0.0, 0.04}
\definecolor{saffron}{rgb}{0.96, 0.77, 0.19}
\definecolor{schoolbusyellow}{rgb}{1.0, 0.85, 0.0}
\definecolor{skyblue}{rgb}{0.53, 0.81, 0.92}
\definecolor{unmellowyellow}{rgb}{1.0, 1.0, 0.4}
\definecolor{wheat}{rgb}{0.96, 0.87, 0.7}
\definecolor{aureolin}{rgb}{0.99, 0.93, 0.0}
\definecolor{persianblue}{rgb}{0.11, 0.22, 0.73}
\definecolor{browna}{rgb}{0.59, 0.29, 0.0}
\definecolor{ufogreen}{rgb}{0.24, 0.82, 0.44}
\definecolor{forestgreen}{rgb}{0.13, 0.55, 0.13}
\definecolor{radicalred}{rgb}{1.0, 0.21, 0.37}
\definecolor{LightGreen}{rgb}{0.56, 0.93, 0.56}
\definecolor{LightCoral}{rgb}{0.94, 0.5, 0.5}
\definecolor{LightBlue}{rgb}{0.68, 0.85, 0.9}
\definecolor{DarkGreen}{rgb}{0.0, 0.2, 0.13}
\definecolor{LimeGreen}{rgb}{0.2, 0.8, 0.2}
\definecolor{DarkRed}{rgb}{0.55, 0.0, 0.0}
\definecolor{Tomato}{rgb}{1.0, 0.39, 0.28}
\definecolor{DarkBlue}{rgb}{0.0, 0.0, 0.55}
\definecolor{forestgreen}{rgb}{0.13, 0.55, 0.13}
\definecolor{zaffre}{rgb}{0.0, 0.08, 0.66}
\definecolor{wildstrawberry}{rgb}{1.0, 0.26, 0.64}
\definecolor{venetianred}{rgb}{0.78, 0.03, 0.08}
\definecolor{selectiveyellow}{rgb}{1.0, 0.73, 0.0}
\definecolor{yaleblue}{rgb}{0.06, 0.3, 0.57}

\begin{document}

	\title{\texttt{DIRECTGO}: A new \direct-type \matlab{} toolbox for derivative-free global optimization}

	\author{Linas Stripinis}
	\authornote{Both authors contributed equally to this research.}
	\email{linas.stripinis@mif.vu.lt}
	\orcid{0000-0001-9680-5847}
	\author{Remigijus Paulavi\v{c}ius}
	\authornotemark[1]
	\email{remigijus.paulavicius@mif.vu.lt}
	\orcid{0000-0003-2057-2922}
	\affiliation{%
		\institution{Vilnius University Institute of Data Science and Digital Technologies}
		\streetaddress{Akademijos 4}
		\city{Vilnius}
		\country{Lithuania}
		\postcode{LT-08663}
	}

	\renewcommand{\shortauthors}{Stripinis and Paulavi\v{c}ius}

	\begin{abstract}
		In this work, we introduce \toolbox, a new \matlab{} toolbox for derivative-free global optimization.
		\toolbox{} collects various deterministic derivative-free \direct{}-type algorithms for box-constrained, generally-constrained, and problems with hidden constraints.
		Each sequential algorithm is implemented in two ways: using static and dynamic data structures for more efficient information storage and organization.
		Furthermore, parallel schemes are applied to some promising algorithms within \toolbox. 
		The toolbox is equipped with a graphical user interface (GUI), ensuring the user-friendly use of all functionalities available in \toolbox.
		Available features are demonstrated in detailed computational studies using a comprehensive \directlib{} library of global optimization test problems.
		Additionally, eleven classical engineering design problems illustrate the potential of \toolbox{} to solve challenging real-world problems.
		Finally, the appendix gives examples of accompanying \matlab{} programs and provides a synopsis of its use on the test problems with box and general constraints.

	\end{abstract}


	\keywords{Global optimization, Derivative-free optimization, \direct-type algorithms, Benchmarking, Optimization software, \texttt{MATLAB}, \texttt{TOMLAB}}

	\maketitle

	\section{Introduction}
	\label{sec:intro}

	The \direct{} (\texttt{DI}viding \texttt{RECT}angles) algorithm~\cite{Jones1993} is a well-known and widely used solution technique for derivative-free global optimization problems.
	The \direct{} algorithm extends classical Lipschitz optimization~\cite{Paulavicius2006,Paulavicius2007,Paulavicius2008,Paulavicius2010:ol,Pinter1996book,Piyavskii1967,Sergeyev2011,Shubert1972}, where the Lipschitz constant is not assumed to be known.
	This property makes \direct-type methods especially attractive for solving various real-world optimization problems (see, e.g.,~\cite{Baker2000,Bartholomew2002,Carter2001,Cox2001,Serafino2011,Gablonsky2001,Liuzzi2010,Paulavicius2019:eswa,Paulavicius2014:book,Stripinis2018b} and the references given therein).
	Moreover, a recent review and comparison in~\cite{Rios2013} revealed that, on average, \direct-type algorithms performance is one of the best among all tested state-of-the-art derivative-free global optimization approaches.
	The \direct-type algorithms often outperform algorithms belonging to other well-known classes, such as Genetic~\cite{John1975}, Simulated annealing~\cite{Kirkpatrick1983}, and Particle swarm optimization~\cite{Kennedy1995}.

	While the original \direct{} addresses only box-constrained optimization problems, various \direct-type modifications and extensions have been proposed.
	Based on the type of constraints, \direct-type algorithms can be classified into four main categories:
	\begin{itemize}
		\item Box-constrained (see, e.g.,~\cite{Finkel2004aa,Finkel2006,Jones2001,Gablonsky2001,Jones1993,Liu2016:GO,Liu2015,Liuzzi2010,Liuzzi2010:jogo,Paulavicius2014:book,Sergeyev2006} and the references given therein);
		\item Linearly-constrained/symmetric (see, e.g.,~\cite{Grbic2013,Paulavicius2014:jogo,Paulavicius2013:jogo,Paulavicius2014:book,Paulavicius2016:ol} and the references given therein);
		\item Generally-constrained (see, e.g.,~\cite{Costa2017,Finkel2005,Jones2001,Liu2017,Stripinis2018b} and the references given therein);
		\item Containing hidden constraints (see, e.g.,~\cite{Gablonsky2001:phd,Na2017,Stripinis2021} and the references given therein).
	\end{itemize}
	\matlab~\cite{Matlab2020} is one of the most broadly used mathematical computing environments in scientific and technical computing (see, e.g., \cite{Burgel2019,Didier2003,Ye2017}).
	Many widely used implementations of the original \direct{} algorithm (see, e.g.,~\cite{Bjorkman1999,Finkel2004,Gablonsky2001:phd}) as well as various later introduced \direct-type extensions (see, e.g.,~\cite{Liu2017,Liu2014,Liu2016:GO,Paulavicius2014:book}), were developed using \matlab.
	Motivated by this, we developed a \direct{}-type global optimization toolbox (\toolbox) within the \matlab{} environment.
	The \toolbox{} toolbox is equipped with a graphical user interface (GUI), which links to a \directlib~\cite{DIRECTGOLibZenodov1.0,DIRECTGOLibv1.0} library and ensures the user-friendly use of all functionalities available in \toolbox.
   The \directlib{} library is a continuation of our previous \texttt{DIRECTLib}~\cite{DIRECTLib2018}, which was widely used in our different previous studies (see, e.g., \cite{Stripinis2018a,Stripinis2018b,Stripinis2020}).
   However, \texttt{DIRECTLib} was designed as a static library and did not offer the global optimization community comfortable opportunities to contribute to it.
   Therefore, a new \directlib{} is designed as an open-source GitHub repository to which other researchers can easily contribute.


	The first publicly available \direct{} implementations and many others introduced later typically use static data structures for storage and organization.
  Our recent work~\cite{Stripinis2020} showed that the \matlab{} implementation of the same \direct-type algorithm based on dynamic data structure often has a significant advantage over implementation based on the static data structure.
	Therefore, each algorithm in \toolbox{} is implemented using both static and dynamic data structures.
	As various applications can benefit from parallel computing, the SPMD (Single Program Multiple Data) parallel scheme (see~\cite{Stripinis2020} for more information) is used to implement some approaches.


	\subsection{Contributions and structure}
	We summarize our main contributions below:
	\begin{itemize}
		\itemsep=0cm
		\item We develop a new \matlab{} toolbox (\toolbox) for derivative-free global optimization, consisting of $36$ different \direct-type algorithms (see \cref{tab:direct_classification} for the details).
		\item We implement each \direct-type algorithm using two types of data structures, static and dynamic~\cite{He2002,Stripinis2020}.
		\item We adapt the SPMD parallel scheme~\cite{Stripinis2020} for selected \direct-type algorithms.
		\item We create a new library (\directlib~\cite{DIRECTGOLibZenodov1.0}) of test and engineering global optimization problems for usage with \toolbox{} and convenient contribution to it through GitHub~\cite{DIRECTGOLibv1.0}.
		\item We perform a comprehensive experimental study on the effectiveness of various \direct-type approaches.
		\item We design a user-friendly application with a graphical user interface (GUI).
		\item We make \toolbox{} open-source, i.e., freely available to anyone~\cite{DGO2021}.
	\end{itemize}

	The rest of the paper is organized as follows.
	\cref{sec:theory} provides the classification of existing \direct-type algorithms and describes the algorithms implemented within our toolbox in more detail.
	The parallel scheme used in implementing some algorithms is also discussed here.
	\toolbox{} toolbox is introduced and described in \cref{sec:tooblox}.
	The detailed computational study of the \toolbox{} toolbox using classical global optimization test and engineering design problems \directlib~\cite{DIRECTGOLibZenodov1.0} are presented in \cref{sec:experimentstests,sec:aplications}, respectively.
	Finally, in \cref{sec:conclusions}, we conclude the work and discuss the possible future directions.

	\begin{table}
		\scriptsize
		\caption{Classification of \direct-type implementations (within the \toolbox{} toolbox) based on the type of constraints.}
		\begin{tabularx}{\textwidth}{p{1.2cm}p{2.6cm}p{0.2cm}p{0.2cm}p{0.2cm}X}
			\toprule
			\multirow{2}{*}{Problem type} & \multirow{2}{*}{Algorithm name} & \multicolumn{3}{c}{Implementation} & \multirow{2}{*}{Description and References} \\
			 	 &  		& st & dy & pa 					     &  \\
			\midrule
			& \directfinkel 						& $+$ & $+$ & $+$ & Finkel’s implementation \cite{Finkel2004} of the original \direct{} \cite{Jones1993} algorithm \\
			& \directrest 					& $+$ & $+$ & $+$ & Our implementation of the algorithm from~\cite{Finkel2004aa} (based on Finkel’s \direct{} \cite{Finkel2004} implementation)\\
			& \directm 		 				& $+$ & $+$ & $+$ & Our implementation of the algorithm from~\cite{Finkel2006} (based on Finkel’s \direct{} \cite{Finkel2004} implementation)\\
			& \directz 						& $+$ & $+$ & $+$ & Our implementation of the algorithm from~\cite{Gablonsky2001} (based on Finkel’s \direct{} \cite{Finkel2004} implementation)\\
			& \directrev 					& $+$ & $+$ & $+$ & Our implementation of the algorithm from~\cite{Jones2001} (based on Finkel’s \direct{} \cite{Finkel2004} implementation)\\
			& \directav 					& $+$ & $+$ & $+$ & Our implementation of the algorithm from~\cite{Liu2013} (based on Finkel’s \direct{} \cite{Finkel2004} implementation)\\
			& \dirmin 						& $+$ & $+$ & $+$ & Our implementation of the algorithm from~\cite{Liuzzi2010} (based on Finkel’s \direct{} \cite{Finkel2004} implementation)\\
			& \plor 						& $+$ & $+$ & $+$ & Our implementation of the algorithm from~\cite{Mockus2017} (based on Finkel’s \direct{} \cite{Finkel2004} implementation)\\

			Box 		& \glbsolve 					& $+$ & $+$ & $-$ & Bj\"{o}rkman’s implementation \cite{Bjorkman1999} of the original \direct{} \cite{Jones1993} algorithm \\
			constrained & \symdirec,~\symdirect 		& $+$ & $+$ & $-$ & Our implementation of algorithms from~\cite{Grbic2013} (based on Bj\"{o}rkman’s \glbsolve{} \cite{Bjorkman1999} implementation)\\
			& \directmr,~\directmro$_{075}$ & $+$ & $+$ & $-$ & Our implementation of algorithms from~\cite{Liu2016:GO,Liu2015} (based on Bj\"{o}rkman’s \glbsolve{} \cite{Bjorkman1999} implementation)\\
			& \birect  						& $+$ & $+$ & $-$ & Our implementation of the algorithm from~\cite{Paulavicius2016:jogo} (based on Bj\"{o}rkman’s \glbsolve{} \cite{Bjorkman1999} implementation)\\
			& \disimplcgb,~\disimplvgb  	& $+$ & $+$ & $-$ & Our implementation of algorithms from~\cite{Paulavicius2014:jogo} (based on Bj\"{o}rkman’s \glbsolve{} \cite{Bjorkman1999} implementation)\\
			& \birectgb,~\birmin,~\gbdirect & $+$ & $+$ & $-$ & Our implementation of algorithms from~\cite{Paulavicius2019:eswa} (based on Bj\"{o}rkman’s \glbsolve{} \cite{Bjorkman1999} implementation)\\
			& \disimplc,~\disimplv 			& $+$ & $+$ & $-$ & Our implementation of algorithms from~\cite{Paulavicius2013:jogo} (based on Bj\"{o}rkman’s \glbsolve{} \cite{Bjorkman1999} implementation)\\
			& \adc  						& $+$ & $+$ & $-$ & Our implementation of the algorithm from~\cite{Sergeyev2006} (based on Bj\"{o}rkman’s \glbsolve{} \cite{Bjorkman1999} implementation)\\
			& \directa 						& $+$ & $+$ & $+$ & Our implementation of the algorithm from~\cite{Baker2000} \\
			& \directg,~\directl,~\directgl & $+$ & $+$ & $+$ & Our implementation of algorithms from~\cite{Stripinis2018a} \\

			\midrule
			Linearly  		& \multirow{2}{*}{\disimpllc,~\disimpllv} & $+$ & $+$ & $-$ & \multirow{2}{*}{ Our implementation of algorithms from~\cite{Paulavicius2014:book,Paulavicius2016:ol} (based on Bj\"{o}rkman’s \glbsolve{} \cite{Bjorkman1999} implementation)}  \\
			constrained  \\
			\midrule
			Generally	 & \directll  		& $+$ & $+$ & $+$ & Finkel’s implementation of the algorithm from~\cite{Finkel2004} \\
			constrained  & \directc,~\directce, 	& $+$ & $+$ & $+$ & \multirow{2}{*}{Our implementation of algorithms from~\cite{Stripinis2018b} (based on our \directgl{} \cite{Stripinis2018a} implementation)}\\
			& \directcemin & \\
			\midrule
			& \directnas		& $+$ & $+$ & $-$ & Finkel’s implementation of the algorithm from~\cite{Gablonsky2001:phd} \\
			Hidden	 		& \directbarrier 	& $+$ & $+$ & $-$ & Our implementation of the algorithm from~\cite{Gablonsky2001:phd} (based on Finkel’s \direct{} \cite{Finkel2004} implementation)\\
			constraints 	& \directsub 		& $+$ & $+$ & $-$ & Our implementation of the algorithm from~\cite{Na2017} (based on Finkel’s \direct{} \cite{Finkel2004} implementation)\\
			& \directh 			& $+$ & $+$ & $-$ & Our implementation of the algorithm from~\cite{Stripinis2021} (based on our \directgl{} \cite{Stripinis2018a} implementation) \\
			\bottomrule
			\multicolumn{6}{l}{\textit{st} - implementation using static data structures.} \\
			\multicolumn{6}{l}{\textit{dy} - implementation using dynamic data structures.} \\
			\multicolumn{6}{l}{\textit{pa} - parallel implementation using dynamic data structures.} \\
		\end{tabularx}
		\label{tab:direct_classification}
	\end{table}

	\section{Theoretical and algorithmic backgrounds}
	\label{sec:theory}

	This section provides the classification of existing \direct-type algorithms and describes the algorithms implemented within our \toolbox{} toolbox in more detail.
	For a thorough review, we refer to a recent survey~\cite{Jones2021}.

	The derivative-free \direct{} algorithm~\cite{Jones1993} is an efficient deterministic technique to solve global optimization~\cite{Horst1995:book,Sergeyev2017:book,Strongin2000:book} problems subject to simple box constraints
	\begin{equation}
		\label{eq:opt-problem1}
		\begin{aligned}
			& \min_{\mathbf{x}\in D} && f(\mathbf{x}), 
		\end{aligned}
	\end{equation}
	where $f:\mathbb{R}^n \rightarrow \mathbb{R}$ denotes the objective function and the feasible region is an $n$-dimensional hyper-rectangle $ D = [ \mathbf{a},  \mathbf{b}] = \{ \mathbf{x} \in \mathbb{R}^n: a^j \leq x^j \leq b^j, j = 1, \dots, n\} $.
	The objective function $f(\mathbf{x})$ is supposed to be Lipschitz-continuous (with an unknown Lipschitz constant) but can be non-linear, non-differentiable, non-convex, and multi-modal.

\begin{figure}[ht]
		\resizebox{.9\textwidth}{!}{
			\begin{tikzpicture}
				\begin{axis}[
					width=0.5\textwidth,height=0.5\textwidth,
					every axis/.append style={font=\large},
					ylabel style={yshift=-0.2cm},
					xlabel = {$x^1$},
					ylabel = {$x^2$},
					label style={font=\large},
					tick label style={font=\large},
					enlargelimits=0.05,
					ytick distance=0.2,
					xtick distance=0.2,
					title={Initialization},
					legend style={draw=none},
					legend columns=1,
					legend style={at={(0.8,-0.15)},font=\normalsize},
					]
					\addlegendimage{only marks,mark=*,color=black}
					\addlegendentry{Sampling point}
					\addplot[thick,patch,mesh,draw,black,patch type=rectangle,line width=0.5mm] coordinates {(0,0) (1,0) (1,1) (0,1)} ;
					\draw [black, thick, mark size=0.1pt, fill=blue!50,opacity=0.4,line width=0.5mm] (axis cs:0,0) rectangle (axis cs:1,1);
					\addplot[only marks,mark=*,black] coordinates {(0.5,0.5)} node[yshift=-8pt] {$c_1$} node[yshift=+8pt] {\small $1408.5$};
				\end{axis}
			\end{tikzpicture}
			\begin{tikzpicture}
				\begin{axis}[
					width=0.5\textwidth,height=0.5\textwidth,
					every axis/.append style={font=\large},
					ylabel style={yshift=-0.2cm},
					xlabel = {$x^1$},
					ylabel = {$x^2$},
					label style={font=\large},
					tick label style={font=\large},
					enlargelimits=0.05,
					ytick distance=0.2,
					xtick distance=0.2,
					title={Iteration $1$},
					legend style={draw=none},
					legend columns=1,
					legend style={at={(0.8,-0.15)},font=\normalsize},
					]
					\addlegendimage{area legend, fill=blue!50,opacity=0.4}
					\addlegendentry{Selected POH}
					\addplot[thick,patch,mesh,draw,black,patch type=rectangle,line width=0.5mm] coordinates {(0,0) (1,0) (1,1) (0,1)} ;
					\addplot[thick,patch,mesh,draw,black,patch type=rectangle,line width=0.5mm] coordinates {(0,0) (1,0) (1,0.3333) (0,0.3333)};
					\addplot[thick,patch,mesh,draw,black,patch type=rectangle,line width=0.5mm] coordinates {(0,0) (1,0) (1,0.6666) (0,0.6666)};
					\addplot[thick,patch,mesh,draw,black,patch type=rectangle,line width=0.5mm] coordinates {(0,0.3333) (0,0.6666) (0.3333,0.6666) (0.3333,0.3333)};
					\addplot[thick,patch,mesh,draw,black,patch type=rectangle,line width=0.5mm] coordinates {(1,0.3333) (1,0.6666) (0.6666,0.6666) (0.6666,0.3333)};
					\draw [black, thick, mark size=0.1pt, fill=blue!50,opacity=0.4,line width=0.5mm] (axis cs:0,0.6666) rectangle (axis cs:1,1);
					\addplot[only marks,mark=o,black] coordinates {(0.5,0.5)} node[yshift=-8pt] {$c_1$} node[yshift=+8pt] {\small $1408.5$};
					\addplot[only marks,mark=o,black] coordinates {(0.5,1/6)} node[yshift=-8pt] {$c_2$} node[yshift=+8pt] {\small $7658.5$};
					\addplot[only marks,mark=*,black] coordinates {(0.5,5/6)} node[yshift=-8pt] {$c_3$} node[yshift=+8pt] {\small $158.5$};
					\addplot[only marks,mark=o,black] coordinates {(1/6,0.5)} node[yshift=-8pt] {$c_4$} node[yshift=+8pt] {\small $1418.5$};
					\addplot[only marks,mark=o,black] coordinates {(5/6,0.5)} node[yshift=-8pt] {$c_5$} node[yshift=+8pt] {\small $288948.5$};
				\end{axis}
			\end{tikzpicture}
			\begin{tikzpicture}
				\begin{axis}[
					width=0.5\textwidth,height=0.5\textwidth,
					every axis/.append style={font=\large},
					ylabel style={yshift=-0.2cm},
					xlabel = {$x^1$},
					ylabel = {$x^2$},
					label style={font=\large},
					tick label style={font=\large},
					enlargelimits=0.05,
					ytick distance=0.2,
					xtick distance=0.2,
					title={Iteration $2$},
					legend style={draw=none},
					legend columns=1,
					legend style={at={(0.8,-0.15)},font=\normalsize},
					]
					\addlegendimage{area legend,black,fill=white,opacity=0.5}
					\addlegendentry{Unselected region}
					\addplot[thick,patch,mesh,draw,black,patch type=rectangle,line width=0.5mm] coordinates {(0,0) (1,0) (1,1) (0,1)} ;
					\draw [black, thick, mark size=0.1pt,line width=0.5mm] (axis cs:0.3333,0.3333) rectangle (axis cs:0.6666,1);
					\draw [black, thick, mark size=0.1pt, fill=blue!50,opacity=0.4,line width=0.5mm] (axis cs:0,0) rectangle (axis cs:1,0.3333);
					\draw [black, thick, mark size=0.1pt, fill=blue!50,opacity=0.4,line width=0.5mm] (axis cs:1/3,2/3) rectangle (axis cs:2/3,1);
					\draw [black, thick, mark size=0.1pt,line width=0.5mm] (axis cs:0.6666,0.3333) rectangle (axis cs:1,0.6666);
					\draw [black, thick, mark size=0.1pt,line width=0.5mm] (axis cs:0,0.6666) rectangle (axis cs:1,1);
					\draw [black, thick, mark size=0.1pt,line width=0.5mm] (axis cs:0,0.3333) rectangle (axis cs:0.3333,0.6666);
					\addplot[only marks,mark=o,black] coordinates {(0.5,0.5)} node[yshift=-8pt] {$c_1$} node[yshift=+8pt] {\small $1408.5$};
					\addplot[only marks,mark=*,black] coordinates {(0.5,1/6)} node[yshift=-8pt] {$c_2$} node[yshift=+8pt] {\small $7658.5$};
					\addplot[only marks,mark=*,black] coordinates {(0.5,5/6)} node[yshift=-8pt] {$c_3$} node[yshift=+8pt] {\small $158.5$};
					\addplot[only marks,mark=o,black] coordinates {(1/6,0.5)} node[yshift=-8pt] {$c_4$} node[yshift=+8pt] {\small $1418.5$};
					\addplot[only marks,mark=o,black] coordinates {(5/6,0.5)} node[yshift=-8pt] {$c_5$} node[yshift=+8pt] {\small $288948.5$};
					\addplot[only marks,mark=o,black] coordinates {(1/6,5/6)} node[yshift=-8pt] {$c_6$} node[yshift=+8pt] {\small $168.5$};
					\addplot[only marks,mark=o,black] coordinates {(5/6,5/6)} node[yshift=-8pt] {$c_7$} node[yshift=+8pt] {\small $237698.5$};
				\end{axis}
		\end{tikzpicture}}
    \caption{Illustration of selection, sampling, and subdivision (trisection) used in the original \direct{} algorithm~\cite{Jones1993} on a two-dimensional \textit{Rosenbrock} test function in the first two iterations.}
		\label{fig:divide}
	\end{figure}
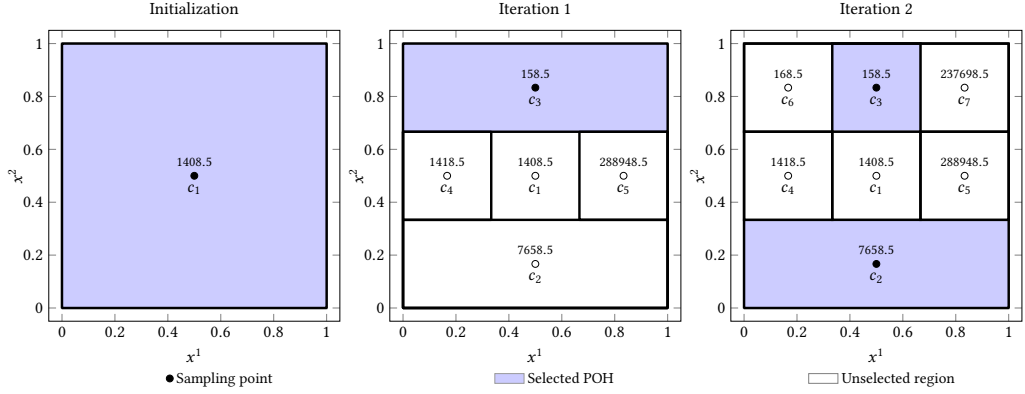

  The \direct{} algorithm includes three main steps: selection, sampling, and partitioning (subdivision).
	At the initial iteration, the \direct{} algorithm normalizes the feasible region $D$ to be the unit hyper-cube $\bar{D}$ and refers to the original space $D$ only when evaluating the objective function.
	Regardless of the dimension, the first evaluation of the objective function is done at the midpoint $\mathbf{c}_1 \in \bar{D}$ (see the left panel of \cref{fig:divide}).
	Then \direct{} selects $\bar{D}$ and samples at $\mathbf{c}_1 \pm \delta e^j, j = 1, \dots, n$, where $e^j$ is the $j$th unit vector and $\delta$ is equal to one-third of the maximum side length of $\bar{D}$.
  The subdivision procedure in \direct{} is based on $n$-dimensional trisection along all longest dimensions (sides).
  When several dimensions have the maximum side length, \direct{} starts trisection from the dimension with the lowest $w^j$ and continues to the highest~\cite{Jones2021,Jones1993}.
  Here $w^j$ is defined as the best function values sampled along dimension $j$
  \begin{equation}
    w^j = \min \{ f(\mathbf{c}_i + \delta_i\mathbf{e}^j),  f(\mathbf{c}_i - \delta_i\mathbf{e}^j) \},
  \end{equation}
  where $\mathbf{c}_i$ is the center of the hyper-rectangle $\bar{D}_i$, and $j \in M$ --- set of dimensions with the maximum side length.
  \Cref{fig:divide} illustrates the \direct{} algorithm's selection, sampling, and subdivision (trisection) for a two-dimensional \textit{Rosenbrock} test function.

	The selection procedure at the initial step is trivial as we have only one candidate $ \bar{D} $.
	To formalize the selection of potentially optimal hyper-rectangles (POHs) in the future iterations, we define the current partition at the iteration $ k $
  	\[
  	\mathcal{P}^k = \{ \bar{D}_i^k : i \in \indexsett^k \},
  	\]
  	where $ \bar{D}_i^k = [\mathbf{a}_i, \mathbf{b}_i ] = \{ \mathbf{x} \in \mathbb{R}^n: 0 \leq a_{i}^{j} \leq x^j \leq b_{i}^{j} \leq 1, j = 1,\dots, n, \forall i \in \indexsett^k \} $ and $ \indexsett^k $ is the index set identifying the current partition $ \mathcal{P}^k $.
  	The next partition $ \mathcal{P}^{k+1} $ is obtained after the subdivision of the selected POHs from the current partition $ \mathcal{P}^k $.
  \direct{} assesses the potentiality based on the lower bound estimates for the objective function $ f $ over each hyper-rectangle $ \bar{D}_i^k $ as stated in~\Cref{def:potOptRect}.

	\begin{definition}
		\label{def:potOptRect}
		(Potentially optimal hyper-rectangle)
		Let $ \mathbf{c}^k_i $ denote the center sampling point and $ \delta^k_i $ be a measure of the hyper-rectangle $ \bar{D}^k_i$.
		Let $ \varepsilon > 0 $ be a positive constant and $f_{\rm min}$ be the best currently found value of the objective function.
		A hyper-rectangle $ \bar{D}^k_j, j \in \indexsett^k $ is said to be potentially optimal if there exists some rate-of-change (Lipschitz) constant $ \tilde{L} > 0$ such that
		\begin{eqnarray}
			f(\mathbf{c}^k_j) - \tilde{L}\delta^k_j & \leq & f(\mathbf{c}^k_i) - \tilde{L}\delta^k_i, \quad \forall i \in \indexsett^k, \label{eqn:potOptRect1} \\
			f(\mathbf{c}^k_j) - \tilde{L}\delta^k_j & \leq & f_{\rm min} - \varepsilon|f_{\rm min}|, \label{eqn:potOptRect2}
		\end{eqnarray}
		where the measure of the hyper-rectangle $ \bar{D}^k_i$ is
		\begin{equation}
			\label{eq:distance}
			\delta^k_i = \frac{1}{2} \| \mathbf{b}^k_i - \mathbf{a}^k_i \|_2.
		\end{equation}
	\end{definition}

	The hyper-rectangle $ D_j^k $ is potentially optimal if the lower Lipschitz bound for the objective function computed by the left-hand side of \eqref{eqn:potOptRect1} is the smallest one with some positive constant $\tilde{L}$ among the hyper-rectangles of the current partition $ \mathcal{P}^k $.
	In~\eqref{eqn:potOptRect2}, the parameter $\varepsilon$ is used to protect from an excessive refinement of the local minima~\cite{Jones1993,Paulavicius2014:jogo}.
	Authors obtained good results for $\varepsilon$ values ranging from $10^{-3}$ to $10^{-7}$ in~\cite{Jones1993}.

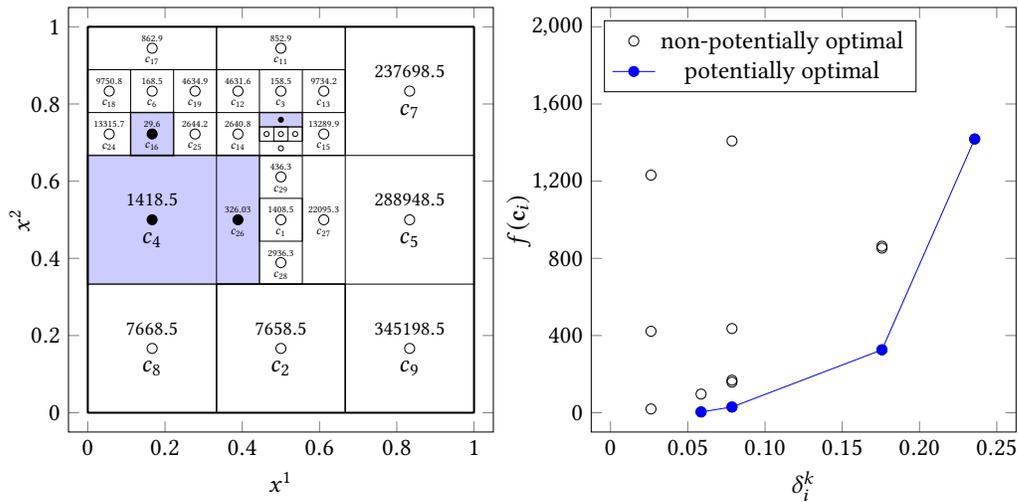
\begin{figure}
		\centering
		\resizebox{.9\textwidth}{!}{%
			\begin{tikzpicture}[baseline]
				\begin{axis}[
					width=0.5\textwidth,height=0.5\textwidth,
					every axis/.append style={font=\large},
					ylabel style={yshift=-0.2cm},
					xlabel = {$x^1$},
					ylabel = {$x^2$},
					label style={font=\large},
					tick label style={font=\large},
					enlargelimits=0.05,
					ytick distance=0.2,
					xtick distance=0.2,
					]
					\draw [black, thick, mark size=0.1pt, fill=blue!50,opacity=0.4] (axis cs:0,1/3) rectangle (axis cs:4/9,2/3);
					\draw [black, thick, mark size=0.1pt, fill=blue!50,opacity=0.4] (axis cs:1/9,2/3) rectangle (axis cs:2/9,7/9);
					\draw [black, thick, mark size=0.1pt, fill=blue!50,opacity=0.4] (axis cs:4/9,20/27) rectangle (axis cs:5/9,21/27);

					\addplot[thick,patch,mesh,draw,black,patch type=rectangle] coordinates {(0,0) (1,0) (1,1) (0,1)} ;
					\addplot[only marks,mark=o,black] coordinates {(0.5,0.5)} node[scale=0.5,yshift=-10pt] {$c_1$} node[scale=0.5,yshift=+8pt] {\small $1408.5$};
					\addplot[only marks,mark=o,black] coordinates {(0.5,1/6)} node[yshift=-8pt] {$c_2$} node[yshift=+8pt] {\small $7658.5$};
					\addplot[only marks,mark=o,black] coordinates {(0.5,5/6)} node[scale=0.5, yshift=-10pt] {$c_3$} node[scale=0.5, yshift=+8pt] {\small$158.5$};
					\addplot[only marks,mark=*,black] coordinates {(1/6,0.5)} node[yshift=-8pt] {$c_4$} node[yshift=+8pt] {\small$1418.5$};
					\addplot[only marks,mark=o,black] coordinates {(5/6,0.5)} node[yshift=-8pt] {$c_5$} node[yshift=+8pt] {\small$288948.5$};
					\addplot[only marks,mark=o,black] coordinates {(1/6,5/6)} node[scale=0.5, yshift=-10pt] {$c_6$} node[scale=0.5, yshift=+8pt] {\small $168.5$};
					\addplot[only marks,mark=o,black] coordinates {(5/6,5/6)} node[yshift=-8pt] {$c_7$} node[yshift=+8pt] {\small$237698.5$};
					\addplot[only marks,mark=o,black] coordinates {(1/6,1/6)} node[yshift=-8pt] {$c_8$} node[yshift=+8pt] {\small$7668.5$};
					\addplot[only marks,mark=o,black] coordinates {(5/6,1/6)} node[yshift=-8pt] {$c_9$} node[yshift=+8pt] {\small$345198.5$};
					\addplot[only marks,mark=o,black, mark size=1pt] coordinates {(0.5,13/18)}; 
					\addplot[only marks,mark=o,black] coordinates {(0.5,17/18)} node[scale=0.5, yshift=-10pt] {$c_{11}$} node[scale=0.5, yshift=+8pt] {\small$852.9$};
					\addplot[only marks,mark=o,black] coordinates {(7/18,15/18)} node[scale=0.5, yshift=-10pt] {$c_{12}$} node[scale=0.5, yshift=+8pt] {\small$4631.6$};
					\addplot[only marks,mark=o,black] coordinates {(11/18,15/18)} node[scale=0.5, yshift=-10pt] {$c_{13}$} node[scale=0.5, yshift=+8pt] {\small$9734.2$};
					\addplot[only marks,mark=o,black] coordinates {(7/18,13/18)} node[scale=0.5, yshift=-10pt] {$c_{14}$} node[scale=0.5, yshift=+8pt] {\small $2640.8$};
					\addplot[only marks,mark=o,black] coordinates {(11/18,13/18)} node[scale=0.5, yshift=-10pt] {$c_{15}$} node[scale=0.5, yshift=+8pt] {\small$13289.9$};
					\addplot[only marks,mark=*,black] coordinates {(1/6,13/18)} node[scale=0.5, yshift=-10pt] {$c_{16}$} node[scale=0.5, yshift=+8pt] {\small $29.6$};
					\addplot[only marks,mark=o,black] coordinates {(1/6,17/18)} node[scale=0.5, yshift=-10pt] {$c_{17}$} node[scale=0.5, yshift=+8pt] {\small $862.9$};
					\addplot[only marks,mark=o,black] coordinates {(1/18,5/6)} node[scale=0.5, yshift=-10pt] {$c_{18}$} node[scale=0.5, yshift=+8pt] {\small $9750.8$};
					\addplot[only marks,mark=o,black, mark size=2pt] coordinates {(5/18,5/6)} node[scale=0.5, yshift=-10pt] {$c_{19}$} node[scale=0.5, yshift=+8pt] {\small $4634.9$};

					\addplot[only marks,mark=o,black, mark size=1pt] coordinates {(1/2,0.685185)};
					\addplot[only marks,mark=*,black, mark size=1pt] coordinates {(1/2,0.759259)};
					\addplot[only marks,mark=o,black, mark size=1pt] coordinates {(0.53703703,13/18)};
					\addplot[only marks,mark=o,black, mark size=1pt] coordinates {(0.46296296,13/18)};
					\addplot[only marks,mark=o,black] coordinates {(1/18,13/18)} node[scale=0.5, yshift=-10pt] {$c_{24}$} node[scale=0.5, yshift=+8pt] {\small $13315.7$};
					\addplot[only marks,mark=o,black] coordinates {(5/18,13/18)} node[scale=0.5, yshift=-10pt] {$c_{25}$} node[scale=0.5, yshift=+8pt] {\small $2644.2$};
					\addplot[only marks,mark=*,black] coordinates {(7/18,1/2)} node[scale=0.5, yshift=-10pt] {$c_{26}$} node[scale=0.5, yshift=+8pt] {\small $326.03$};
					\addplot[only marks,mark=o,black] coordinates {(11/18,1/2)} node[scale=0.5, yshift=-10pt] {$c_{27}$} node[scale=0.5, yshift=+8pt] {\small $22095.3$};
					\addplot[only marks,mark=o,black] coordinates {(1/2,7/18)} node[scale=0.5, yshift=-10pt] {$c_{28}$} node[scale=0.5, yshift=+8pt] {\small $2936.3$};
					\addplot[only marks,mark=o,black] coordinates {(1/2,11/18)} node[scale=0.5, yshift=-10pt] {$c_{29}$} node[scale=0.5, yshift=+8pt] {\small $436.3$};

					\addplot[patch,mesh,draw,black,patch type=rectangle] coordinates {
						(0,0) (1/3,0) (1/3,1) (0,1)
						(2/3,0) (1,0) (1,1) (2/3,1)
						(1/3,0) (2/3,0) (2/3,1/3) (1/3,1/3)
						(1/3,0) (2/3,0) (2/3,1/3) (1/3,1/3)
						(0,1/3) (1,1/3) (1,2/3) (0,2/3)
						(1/3,2/3) (2/3,2/3) (2/3,1) (1/3,1)
						(0,8/9) (2/3,8/9) (2/3,7/9) (0,7/9)
						(4/9,6/9) (4/9,8/9) (5/9,8/9) (5/9,6/9)
						(4/9,3/9) (4/9,6/9) (5/9,6/9) (5/9,3/9)
						(1/9,6/9) (1/9,8/9) (2/9,8/9) (2/9,6/9)
						(4/9,4/9) (5/9,4/9) (5/9,5/9) (4/9,5/9)
						(4/9,19/27) (5/9,19/27) (5/9,20/27) (4/9,20/27)
						(13/27,19/27) (13/27,20/27) (14/27,20/27) (14/27,19/27)};
				\end{axis}
			\end{tikzpicture}
			\begin{tikzpicture}[baseline]
				\begin{axis}[
					width=0.5\textwidth,height=0.5\textwidth,
					every axis/.append style={font=\large},
					ylabel style={yshift=-0.2cm},
					xlabel = {$\delta_i^k$},
					ymin=0,ymax=2000,
					xmin=0,xmax=0.25,
					xtick distance=0.05,
					xtick = {0,0.05,0.10,0.15,0.20,0.25},
					xticklabels = {$0$,$0.05$,$0.10$,$0.15$,$0.20$,$0.25$},
					ytick distance=400,
					restrict y to domain=0:2000,
					ylabel = {$f(\mathbf{c}_i)$},
					legend style={font=\large},
					legend pos=north west,
					enlargelimits=0.05,
					]

					\addplot[black,only marks,mark=o,mark size=2pt] table[x=Xx,y=Yy] {figure_data/poh.txt};

					\addplot[blue,mark=*,mark size=2pt] coordinates {(0.0586, 4.179012346) (0.0786,29.61111111) (0.1757, 326.0308642) (0.2357,1418.5)};
					\legend{non-potentially optimal, potentially optimal};
				\end{axis}
		\end{tikzpicture}}
		\caption{Visualization of selected potentially optimal rectangles in the fifth iteration of the \direct{} algorithm on a two-dimensional \textit{Rosenbrock} test problem.}
		\label{fig:poh}
	\end{figure}

	A geometric interpretation of the selection procedure is given on the right panel of~\cref{fig:poh}.
	Here, each hyper-rectangle is represented as a point.
  The $x$-axis shows the measure $(\delta^k_i)$ while the $y$-axis -- the objective function value attained at the midpoint $(\mathbf{c}^k_i)$ of a certain hyper-rectangle.
  The hyper-rectangles meeting conditions \eqref{eqn:potOptRect1} and \eqref{eqn:potOptRect2} are points on the lower-right convex hull (highlighted in blue color).
	Condition \eqref{eqn:potOptRect2} prevents wasting function evaluations on tiny hyper-rectangles where only a negligible improvement can be expected.


	Then at each subsequent iteration, \direct{} performs a selection of POHs, which are sampled, evaluated, and trisected.
	Almost all \direct-type extensions and modifications follow the same algorithmic framework, summarized in Algorithm~\ref{alg:direct}.

	\begin{algorithm}
		\textbf{Initialization}. Normalize the search space $D$ to be the unit hyper-rectangle $\bar{D}$, but refer to the original space $D$ when making function calls.
		Evaluate the objective $f$ at the center point $\mathbf{c}_1$.
		Set $f_{\rm min} = f(\mathbf{c}_1)$, $c_{\rm min} = \mathbf{c}_1$. Initialize algorithmic performance measures, and \textit{stopping criteria}.\label{alg:initialization}

		\While{stopping criteria are not satisfied}{
			\textbf{Selection}. \textit{Identify} the sets $S$ of POHs (subregions of $\bar{D}$).\label{alg:selection}

			\textbf{Sampling}. For each POH $(\bar{D}_j \in S)$ \textit{sample} and \textit{evaluate} the objective function at new domain points.
			Update $f_{\rm min}, \mathbf{c}_{\rm min}$, and algorithmic performance measures.\label{alg:sampling}

			\textbf{Subdivision}. Each POH $(\bar{D}_j \in S)$ subdivide (trisect) and update the partition $(\mathcal{P})$. \label{alg:subdivision}
		}
		\textbf{Return} $f_{\rm min}, \mathbf{c}_{\rm min}$, and performance measures.
		\caption{Main steps of \direct-type algorithms}
		\label{alg:direct}
	\end{algorithm}

	\subsection{\direct-type algorithms for box-constrained global optimization}
	\label{sec:box}

	Many different \direct{} extensions have been suggested.
	Most of them focused on improving the selection of POHs, while others introduced new partitioning and sampling strategies.
	The summary of all box-constrained proposals considered in the \toolbox{} toolbox is given in~\cref{tab:algorithms}.
	Most algorithms are based on the trisection of $n$-dimensional POHs, and just \adc, \birect, and both \texttt{DISIMPL} versions use different partitioning strategies.
	Below we briefly review the \direct-type approaches for box-constrained global optimization implemented in the current release of the \toolbox{} toolbox.

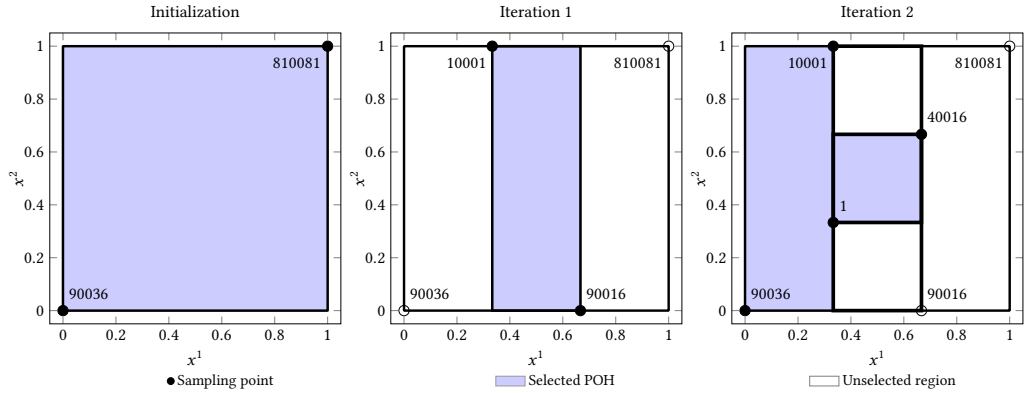
\begin{figure}[ht]
		\resizebox{.9\textwidth}{!}{
			\begin{tikzpicture}
				\begin{axis}[
					width=0.5\textwidth,height=0.5\textwidth,
					every axis/.append style={font=\large},
					ylabel style={yshift=-0.2cm},
					xlabel = {$x^1$},
					ylabel = {$x^2$},
					label style={font=\large},
					tick label style={font=\large},
					enlargelimits=0.05,
					ytick distance=0.2,
					xtick distance=0.2,
					title={Initialization},
					legend style={draw=none},
					legend columns=1,
					legend style={at={(0.8,-0.15)},font=\normalsize},
					]
					\addlegendimage{only marks,mark=*,color=black}
					\addlegendentry{Sampling point}
					\addplot[thick,patch,mesh,draw,black,patch type=rectangle,line width=0.5mm] coordinates {(0,0) (1,0) (1,1) (0,1)} ;
					\draw [black, thick, mark size=0.1pt, fill=blue!50,opacity=0.4,line width=0.5mm] (axis cs:0,0) rectangle (axis cs:1,1);
					\addplot[only marks,mark size=3pt,mark=*,black] coordinates {(0,0)} node[black,yshift=10pt,xshift=15pt] {$90036$};
					\addplot[only marks,mark size=3pt,mark=*,black] coordinates {(1,1)} node[black,yshift=-10pt,xshift=-18pt] {$810081$};
				\end{axis}
			\end{tikzpicture}
			\begin{tikzpicture}
				\begin{axis}[
					width=0.5\textwidth,height=0.5\textwidth,
					every axis/.append style={font=\large},
					ylabel style={yshift=-0.2cm},
					xlabel = {$x^1$},
					ylabel = {$x^2$},
					label style={font=\large},
					tick label style={font=\large},
					enlargelimits=0.05,
					ytick distance=0.2,
					xtick distance=0.2,
					title={Iteration $1$},
					legend style={draw=none},
					legend columns=1,
					legend style={at={(0.8,-0.15)},font=\normalsize},
					]
					\addlegendimage{area legend, fill=blue!50,opacity=0.4}
					\addlegendentry{Selected POH}
					\addplot[thick,patch,mesh,draw,black,patch type=rectangle,line width=0.5mm] coordinates {(0,0) (1,0) (1,1) (0,1)} ;
					\draw [black, thick, mark size=0.1pt, fill=blue!50,opacity=0.4,line width=0.5mm] (axis cs:0.3333,0) rectangle (axis cs:0.6666,1);
					\draw [black, thick, mark size=0.1pt,line width=0.5mm] (axis cs:0.3333,0) rectangle (axis cs:0.6666,1);
					\addplot[only marks,mark size=3pt,mark=o,black] coordinates {(0,0)} node[black,yshift=10pt,xshift=15pt] {$90036$};
					\addplot[only marks,mark size=3pt,mark=o,black] coordinates {(1,1)} node[black,yshift=-10pt,xshift=-18pt] {$810081$};
					\addplot[only marks,mark size=3pt,mark=*,black] coordinates {(1/3,1)} node[black,yshift=-10pt,xshift=-15pt] {$10001$};
					\addplot[only marks,mark size=3pt,mark=*,black] coordinates {(2/3,0)} node[black,yshift=10pt,xshift=15pt] {$90016$};
				\end{axis}
			\end{tikzpicture}
			\begin{tikzpicture}
				\begin{axis}[
					width=0.5\textwidth,height=0.5\textwidth,
					every axis/.append style={font=\large},
					ylabel style={yshift=-0.2cm},
					xlabel = {$x^1$},
					ylabel = {$x^2$},
					label style={font=\large},
					tick label style={font=\large},
					enlargelimits=0.05,
					ytick distance=0.2,
					xtick distance=0.2,
					title={Iteration $2$},
					legend style={draw=none},
					legend columns=1,
					legend style={at={(0.8,-0.15)},font=\normalsize},
					]
					\addlegendimage{area legend,black,fill=white,opacity=0.5}
					\addlegendentry{Unselected region}
					\addplot[thick,patch,mesh,draw,black,patch type=rectangle,line width=0.5mm] coordinates {(0,0) (1,0) (1,1) (0,1)} ;
					\draw [black, thick, mark size=0.1pt, fill=blue!50,opacity=0.4,line width=0.5mm] (axis cs:0,0) rectangle (axis cs:0.3333,1);
					\draw [black, thick, mark size=0.1pt, fill=blue!50,opacity=0.4,line width=0.5mm] (axis cs:0.3333,0.3333) rectangle (axis cs:0.6666,0.6666);
					\draw [black, thick, mark size=0.1pt,line width=0.75mm] (axis cs:0.3333,0) rectangle (axis cs:0.6666,1);
					\draw [black, thick, mark size=0.1pt,line width=0.75mm] (axis cs:0.3333,0.3333) rectangle (axis cs:0.6666,0.6666);
					\addplot[only marks,mark size=3pt,mark=*,black] coordinates {(0,0)} node[black,yshift=10pt,xshift=15pt] {$90036$};
					\addplot[only marks,mark size=3pt,mark=o,black] coordinates {(1,1)} node[black,yshift=-10pt,xshift=-18pt] {$810081$};
					\addplot[only marks,mark size=3pt,mark=*,black] coordinates {(1/3,1)} node[black,yshift=-10pt,xshift=-15pt] {$10001$};
					\addplot[only marks,mark size=3pt,mark=o,black] coordinates {(2/3,0)} node[black,yshift=10pt,xshift=15pt] {$90016$};
					\addplot[only marks,mark size=3pt,mark=*,black] coordinates {(2/3,2/3)} node[black,yshift=10pt,xshift=15pt] {$40016$};
					\addplot[only marks,mark size=3pt,mark=*,black] coordinates {(1/3,1/3)} node[black,yshift=10pt,xshift=6pt] {$1$};
				\end{axis}
		\end{tikzpicture}}
		\caption{Illustration of the diagonal trisection strategy introduced in the \adc{} algorithm on a two-dimensional \textit{Rosenbrock} test function in the first two iterations.}
		\label{fig:divide_adc}
	\end{figure}

  Adaptive diagonal curves (\adc) based algorithm with a new two-phase technique balancing local and global information was introduced in \cite{Sergeyev2006}.
  Independently on  dimensionality, the \adc{} algorithm evaluates the objective function at two vertices $\mathbf{a}^k_i$ and $\mathbf{b}^k_i$ of the main diagonal, as shown in \cref{fig:divide_adc}.
  Notice that up to $2^n$ hyper-rectangles can share the same vertex, leading (in a long sequence) to a smaller number of sampled points than the total number of hyper-rectangles in the current partition.
  Furthermore, as in the revised version of \direct{} \cite{Jones2001}, \adc{} trisects each selected POH along just one of the longest dimensions.
  Such a diagonal scheme potentially obtains more comprehensive information about the objective function than center-based sampling, which sometimes may take many iterations to find the solution.
  For example, a hyper-rectangle containing the optimum with a bad function value at the midpoint makes him undesirable for further selection.
  The \adc{} algorithm intuitively reduces this chance for both sampling points in the hyper-rectangle containing the optimum solution.

\begin{figure}[ht]
		\resizebox{.9\textwidth}{!}{
			\begin{tikzpicture}
				\begin{axis}[
					width=0.5\textwidth,height=0.5\textwidth,
					every axis/.append style={font=\large},
					ylabel style={yshift=-0.2cm},
					xlabel = {$x^1$},
					ylabel = {$x^2$},
					label style={font=\large},
					tick label style={font=\large},
					enlargelimits=0.05,
					ytick distance=0.2,
					xtick distance=0.2,
					title={Initialization},
					legend style={draw=none},
					legend columns=1,
					legend style={at={(0.8,-0.15)},font=\normalsize},
					]
					\addlegendimage{only marks,mark=*,color=black}
					\addlegendentry{Sampling point}
					\addplot[thick,patch,mesh,draw,black,patch type=rectangle,line width=0.5mm] coordinates {(0,0) (1,0) (1,1) (0,1)} ;
					\draw [black, thick, mark size=0.1pt, fill=blue!50,opacity=0.4,line width=0.5mm] (axis cs:0,0) rectangle (axis cs:1,1);
					\addplot[only marks,mark size=3pt,mark=*,black] coordinates {(1/3,1/3)} node[black,yshift=-10pt,xshift=-6pt] {$1$};
					\addplot[only marks,mark size=3pt,mark=*,black] coordinates {(2/3,2/3)} node[black,yshift=10pt,xshift=15pt] {$40016$};
				\end{axis}
			\end{tikzpicture}
			\begin{tikzpicture}
				\begin{axis}[
					width=0.5\textwidth,height=0.5\textwidth,
					every axis/.append style={font=\large},
					ylabel style={yshift=-0.2cm},
					xlabel = {$x^1$},
					ylabel = {$x^2$},
					label style={font=\large},
					tick label style={font=\large},
					enlargelimits=0.05,
					ytick distance=0.2,
					xtick distance=0.2,
					title={Iteration $1$},
					legend style={draw=none},
					legend columns=1,
					legend style={at={(0.8,-0.15)},font=\normalsize},
					]
					\addlegendimage{area legend,black, fill=blue!50,opacity=0.4}
					\addlegendentry{Selected POH}
					\addplot[thick,patch,mesh,draw,black,patch type=rectangle,line width=0.5mm] coordinates {(0,0) (1,0) (1,1) (0,1)} ;
					\draw [black, thick, mark size=0.1pt, fill=blue!50,opacity=0.4,line width=0.5mm] (axis cs:0,0) rectangle (axis cs:0.5,1);
					\draw [black, thick, mark size=0.1pt,line width=0.5mm] (axis cs:0,0) rectangle (axis cs:0.5,1);
					\addplot[only marks,mark size=3pt,mark=*,black] coordinates {(1/3,1/3)} node[black,yshift=-10pt,xshift=-6pt] {$1$};
					\addplot[only marks,mark size=3pt,mark=o,black] coordinates {(2/3,2/3)} node[black,yshift=10pt,xshift=15pt] {$40016$};
					\addplot[only marks,mark size=3pt,mark=*,black] coordinates {(1/6,2/3)} node[black,yshift=-10pt,xshift=15pt] {$168.5$};
					\addplot[only marks,mark size=3pt,mark=o,black] coordinates {(5/6,1/3)} node[black,yshift=-10pt,xshift=-15pt] {$31654$};
				\end{axis}
			\end{tikzpicture}
			\begin{tikzpicture}
				\begin{axis}[
					width=0.5\textwidth,height=0.5\textwidth,
					every axis/.append style={font=\large},
					ylabel style={yshift=-0.2cm},
					xlabel = {$x^1$},
					ylabel = {$x^2$},
					label style={font=\large},
					tick label style={font=\large},
					enlargelimits=0.05,
					ytick distance=0.2,
					xtick distance=0.2,
					title={Iteration $2$},
					legend style={draw=none},
					legend columns=1,
					legend style={at={(0.8,-0.15)},font=\normalsize},
					]
					\addlegendimage{area legend,black,fill=white,opacity=0.5}
					\addlegendentry{Unselected region}
					\addplot[thick,patch,mesh,draw,black,patch type=rectangle,line width=0.5mm] coordinates {(0,0) (1,0) (1,1) (0,1)} ;
					\draw [black, thick, mark size=0.1pt, fill=blue!50,opacity=0.4,line width=0.5mm] (axis cs:0.5,0) rectangle (axis cs:1,1);
					\draw [black, thick, mark size=0.1pt, fill=blue!50,opacity=0.4,line width=0.5mm] (axis cs:0,0) rectangle (axis cs:0.5,0.5);
					\draw [black, thick, mark size=0.1pt ,line width=0.5mm] (axis cs:0,0) rectangle (axis cs:0.5,1);
					\draw [black, thick, mark size=0.1pt ,line width=0.5mm] (axis cs:0,0) rectangle (axis cs:0.5,0.5);
					\addplot[only marks,mark size=3pt,mark=*,black] coordinates {(1/3,1/3)} node[black,yshift=-10pt,xshift=-6pt] {$1$};
					\addplot[only marks,mark size=3pt,mark=*,black] coordinates {(2/3,2/3)} node[black,yshift=10pt,xshift=15pt] {$40016$};
					\addplot[only marks,mark size=3pt,mark=o,black] coordinates {(1/6,2/3)} node[black,yshift=-10pt,xshift=15pt] {$168.5$};
					\addplot[only marks,mark size=3pt,mark=*,black] coordinates {(5/6,1/3)} node[black,yshift=-10pt,xshift=-15pt] {$31654$};
					\addplot[only marks,mark size=3pt,mark=*,black] coordinates {(1/6,1/6)} node[black,yshift=-10pt,xshift=15pt] {$76685$};
					\addplot[only marks,mark size=3pt,mark=o,black] coordinates {(1/3,5/6)} node[black,yshift=10pt,xshift=-15pt] {$56260$};
				\end{axis}
		\end{tikzpicture}}
		\caption{Illustration of the diagonal bisection strategy used in the \birect{} algorithm on a two-dimensional \textit{Rosenbrock} test function in the first two iterations.}
		\label{fig:divide_birect}
	\end{figure}
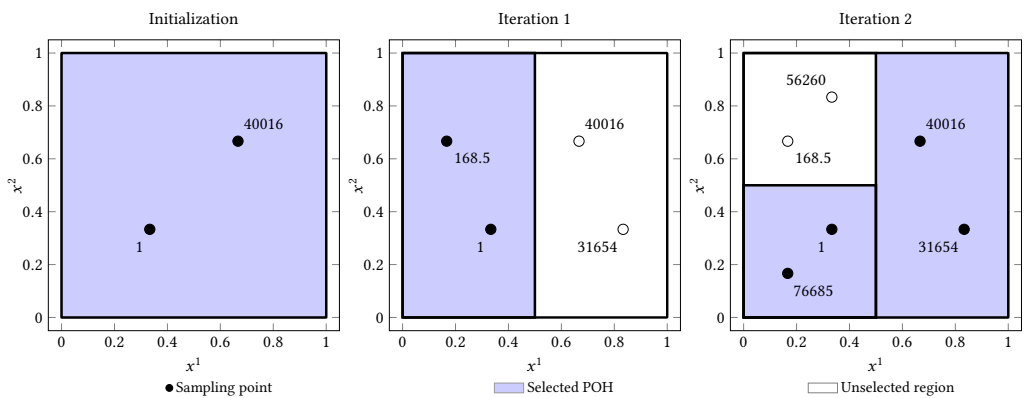

	\birect{} (\texttt{BI}secting \texttt{RECT}angles)~\cite{Paulavicius2016:jogo} is motivated by the diagonal partitioning strategy~\cite{Sergeyev2006,Sergeyev2008:book,Sergeyev2017:book}.
	As the name suggests, the bisection is used instead of a trisection typical for \direct-type algorithms.
	In \birect, the objective function is evaluated at two points on the diagonal equidistant between themselves and a diagonal's vertices, as shown in \cref{fig:divide_birect}.
	Such a sampling strategy enables the reuse of sampling points in descendant hyper-rectangles.
  Moreover, as in the \adc{} case, using two-points-based diagonal sampling, potentially more comprehensive information about the objective function is considered than in the center-based sampling.

	In \texttt{DISIMPL}~\cite{Paulavicius2013:jogo}, simplicial partitions are considered instead of hyper-rectangles.
  The hyper-cube $\bar{D}$ is partitioned into $n!$ simplices by the standard face-to-face simplicial division based on the combinatorial vertex triangulation at the first iteration.
	After this, all simplices share the diagonal of the feasible region and have equal hyper-volume.
	In~\cite{Paulavicius2013:jogo}, we proposed two different sampling strategies.
	Both are included in the \toolbox{} toolbox: i) \disimplc{} evaluating the objective function at the geometric center of the simplex; ii) \disimplv{} evaluating the objective function on all unique vertices of the simplex.
  For box-constrained problems, the total number of initial simplices grows speedily with the dimension increase.
	Therefore, \texttt{DISIMPL} effectively can be used only for small box-constrained problems.
	However, the \texttt{DISIMPL} approach is auspicious (among all \direct-type methods) for symmetric optimization problems~\cite{Paulavicius2013:jogo,Paulavicius2014:book} and problems with linear constraints~\cite{Paulavicius2016:ol}.

	In the \directrest{} algorithm \cite{Finkel2004aa}, the authors introduced an adaptive scheme for the $\varepsilon$ parameter.
	Condition \eqref{eqn:potOptRect2} is needed to stop the \direct{} from wasting function evaluations on minor hyper-rectangles where only a negligible improvement can be expected.
	The \directrest{} algorithm starts with $\varepsilon = 0$, and the same value for $\varepsilon$ is maintained while improvement is achieved.
	However, if five consecutive iterations have no improvement in the best function value, the search may be stagnated around a local optimum.
	Therefore, the algorithm switches to $\varepsilon = 0.01$ value to prevent an excessive local search.
	If the algorithm finds an improvement or fails to see the progress within 50 iterations in this phase, \directrest{} switches to $\varepsilon = 0$.
	Then, if another $50$ iterations pass without improvement, this may indicate that the global minimum has been found, and one should work on refining it to higher accuracy.

	The authors of \directmr~\cite{Liu2015} and \directmro$_{075}$~\cite{Liu2015b} algorithms introduced three different levels to perform the selection procedure:
	\begin{itemize}
		\item At level 2, \direct{} is run as usual, with $\varepsilon = 10^{-5}$.
		\item At level 1, the selection is limited to only $90 \%$ of $\bar{D}^k_i \in \mathcal{P}^k$; $10 \%$ of the largest hyper-rectangles are ignored.
    Here, $\varepsilon = 10^{-7}$ is used.
		\item At level 0, the selection is limited to $10 \%$ of the largest hyper-rectangles (ignored at level 1) using $\varepsilon = 0$.
	\end{itemize}
	Both algorithms cycle through these levels using “\textit{W-cycle}”: $21011012$.
	The main difference between the proposed algorithms is that \directmr{} uses fixed $\varepsilon = 10^{-4}$ value at all levels, while \directmro$_{075}$ follows above-mentioned rules.

	In~\cite{Baker2000}, the authors relaxed the selection criteria of POHs and proposed an aggressive version of the \direct{} algorithm.
	\directa's main idea is to select and divide at least one hyper-rectangle from each group of different diameters $(\delta^k_i)$ having the lowest function value.
	Therefore, using \directa{} in the situation presented in~\cref{fig:poh}, a hyper-rectangle with the slightest measure $\delta_i^k$ would also be selected and divided.
	The aggressive version performs more function evaluations per iteration than other \direct-type methods.
	From the optimization point of view, such an approach seems less favorable since it ``wastes'' function evaluations by exploring unnecessary (non-potentially optimal) hyper-rectangles.
	However, such a strategy is much more appealing in a parallel environment, as was shown in~\cite{He2009part2,He2009part1,He2010,Watson2001}.
  Note that the authors did not specify which hyper-rectangle should be selected from the same group $(\delta^k_i)$ if more than one with identical objective function values exist.
  Thus, the hyper-rectangle with the larger index value was selected in our implementations.

	In~\cite{Gablonsky2001}, the algorithm named \directz{} was proposed.
	In most \direct-type{} algorithms, the measure of the hyper-rectangle is calculated by a half-length of a diagonal (see~\eqref{eq:distance}).
	In \directz{}, this measure is evaluated by the length of its longest side.
	This measure corresponds to the infinity norm and allows the \directz{} algorithm to group more hyper-rectangles with the same measure.
	Thus, there are fewer different measures, so fewer POHs are selected
	Moreover, with \directz{} at most one hyper-rectangle selected from each group, even if there is more than one POH in the same group.
	Such a strategy allows a reduction in the number of divisions within a group.
	Once again, the same rule is adapted \cite{Stripinis2018a} to determine which hyper-rectangle to select from several possible ones.

	In~\cite{Finkel2006}, the authors concluded that the original \direct{} algorithm is sensitive to the objective function's additive scaling.
  Additionally, the algorithm does not operate well when the objective function values are large enough.
  The authors proposed a scaling of function values by subtracting the median ($f_{\rm median}$) of the collected function values to overcome this.
	\directm{} replaces the equation \eqref{eqn:potOptRect2} in \cref{eqn:potOptRect2} to:
	\begin{equation}
		f(\mathbf{c}_j) - \tilde{L}\delta_j \leq f_{\rm min} - \varepsilon|f_{\rm min} - f_{\rm median}|. \label{eqn:potOptRect3}
	\end{equation}
	Similarly, in~\cite{Liu2013}, the authors extended the same idea in \directav{} to reduce the objective function's additive scaling.
  Instead of the median value, the authors proposed to use the average value $(f_{\rm average})$ at each iteration
	\begin{equation}
		f(\mathbf{c}_j) - \tilde{L}\delta_j \leq f_{\rm min} - \varepsilon|f_{\rm min} - f_{\rm average}|. \label{eqn:potOptRect4}
	\end{equation}
	Another extension of the \direct{} algorithm was proposed in~\cite{Grbic2013}.
	The authors introduced \symdirec{} (\symdirect) as \direct{} extensions for symmetric Lipschitz continuous functions.
	When solving symmetric optimization problems, there exist equivalent subregions in the hyper-rectangle.
	The algorithm determines which hyper-rectangles can be safely discarded, considering the problem’s symmetrical nature, and avoids exploration over equivalent subregions.

	In the \plor{} algorithm~\cite{Mockus2017}, the set of POHs is reduced to just two, corresponding to the first and last point on the Pareto front (see the right panel in \cref{fig:poh}).
  Therefore, only hyper-rectangles with the lowest function value and the most extensive measure, breaking ties in favor of a better center-point function value, are selected.

	Our recent extension, \directgl{}~\cite{Stripinis2018a}, introduced a new approach to identifying the extended set of POHs.
	Here, using a novel two-step-based strategy, the set of the best hyper-rectangles is enlarged by adding more medium-measured hyper-rectangles with the smallest function value at their centers and, additionally, closest to the current minimum point.
	The first step of the selection procedure forces the \directgl{} algorithm to work more globally (compared to the selection used in \direct~\cite{Jones1993}).
	In contrast, the second step assures a faster and broader examination around the current minimum point.
	The original \directgl{} version performs a selection of POHs in each iteration twice~\cite{Stripinis2018a}, and the algorithm separately handles the found independent sets $G$ (using Definition 2 from~\cite{Stripinis2018a} - \directg) and $L$ (using Definition 3 from~\cite{Stripinis2018a} - \directl).
	Following the same trend from~\cite{Stripinis2020}, the version used in this paper slightly differs compared to~\cite{Stripinis2018a}.
	In the current version of \directgl{}, identifying these two sets is performed in succession, and the unique union of these two sets ($S = G \cup L$) is used in Algorithm~\ref{alg:direct}, Line~\ref{alg:selection}.
	This modification was introduced to reduce the data communication between the computational units in the parallel algorithm version~\cite{Stripinis2020pdirectglce}.
  At the same time, we found that this way modified \directgl{} was, on average, more effective than the original one.

	Several globally biased (\texttt{Gb-}) versions of \direct-type algorithms were introduced and investigated~\cite{Paulavicius2014:jogo,Paulavicius2019:eswa}.
	Proposed approaches are primarily oriented for solving extremely difficult global optimization problems and contain a phase that constrains itself to large subregions.
	The introduced step performs until a sufficient number of divisions of hyper-rectangles near the current best point is done.
	Once those subdivisions around the current best minima point are performed, the neighborhood contains only small measure hyper-rectangles and all larger ones located far away from it.
	Therefore, the two-phase strategy makes the \direct-type algorithms examine larger hyper-rectangles and return to the general phase only when an improved minimum is obtained.
	The proposed globally biased strategy is combined with \glbsolve, \birect, \disimplc, and \disimplv{} algorithmic frameworks within our \toolbox{} toolbox.

	Finally, three different hybridized \direct-type algorithms are proposed (\directrev~\cite{Jones2001}, \dirmin~\cite{Liuzzi2010}, \birmin~\cite{Paulavicius2019:eswa}).
	In our implementation, all algorithms are combined with the same local search routine -- \textit{fmincon}.
	The \dirmin{} algorithm suggests running a local search starting from the midpoint of every POH.
	However, such an approach likely generates more local searches than necessary, as many start points will converge to the same local optimum.
	The other authors of the \directrev{} and \birmin{} algorithms tried to minimize the usage of local searches.
	They suggested using \textit{fmincon} only when some improvement in the best current solution is obtained.
	The authors in \cite{Jones2001} additionally incorporated the following two enhancements.
	First, in the \directrev{} algorithm, selected hyper-rectangles are trisected only on one longest side.
	Second, only one POH is selected if several equally good exist (the same measure and objective values) in \cref{def:potOptRect}.
	We have applied the same rule from \cite{Stripinis2018a} to determine which hyper-rectangle to select from several ones when needed.

	\begin{table}
		\caption{Summary of main characteristics of \direct-type algorithms for box-constrained global optimization}
		\label{tab:algorithms}
		\small
		\resizebox{\textwidth}{0.6\textwidth}{
			\setlength{\extrarowheight}{4pt}
			\begin{tabular}{|C{1.8cm}|C{2.45cm}|C{2.2cm}|C{0.25cm}|C{1cm}|p{6cm}|C{0.75cm}|}
				\toprule
				\backslashbox{\textbf{Alg.}}{\textbf{Step}}
				& \textit{Partitioning scheme} & \textit{Sampling scheme} & \multicolumn{3}{c|}{\textit{Selection of POH \& comments on additional steps, if any}}  & \textit{Input param.}  \\
				\midrule
				\directa &  & \multirow{18}{*}{\rotatebox{90}{\parbox{9.5cm}{Samples midpoints of the hyper-rectangles}}} & \multicolumn{3}{l|}{\parbox{8cm}{Relaxed criteria of POH selection. At each iteration, the algorithm selects and divides the best hyper-rectangles of each size $(\delta^k_i)$.}} & \multirow{4}{*}{\rotatebox{90}{\parbox{2.7cm}{No input parameters}}} \\[6pt]
				\cline{0-0}\cline{4-6}
				\directg & \multirow{8}{*}{\rotatebox{90}{ \parbox{4.5cm}{Hyper-rectangular partitions-based on $n$-dimensional trisection}}} & &\multicolumn{3}{l|}{\parbox{8cm}{Uses enhanced global selection (Definition 2 in~\cite{Stripinis2018a}).}} &   \\[2pt]
				\cline{0-0}\cline{4-6}
				\directl & & & \multicolumn{3}{l|}{\parbox{8cm}{Uses enhanced local selection (Definition 3 in~\cite{Stripinis2018a}).}} &   \\[2pt]
				\cline{0-0}\cline{4-6}
				\vspace{10pt}
				\directgl &  & & \multicolumn{3}{l|}{\parbox{8cm}{Uses the unique union of two sets obtained using enhanced global and local selection (Definitions 2 and 3~\cite{Stripinis2018a}).}} &   \\[9pt]
				\cline{0-0}\cline{4-7}
				\direct & & & \multirow{15}{*}{\rotatebox{90}{\parbox{10cm}{Uses \cref*{def:potOptRect}}}} & \multicolumn{2}{l|}{} &  \multirow{15}{*}{\rotatebox{90}{\parbox{10.5cm}{Balance parameter $\epsilon$}}} \\[1pt]
				\cline{0-0}\cline{5-6}
				\vspace{5pt}
				\texttt{DIRECT-} \texttt{restart} 	&&&& \multicolumn{2}{l|}{\parbox{7.5cm}{Uses two different $\varepsilon$ values during the selection: $0$ if there is an improvement in the solution, and $0.01$ otherwise.}} &   \\[6pt]
				\cline{0-0}\cline{5-6}
				\vspace{5pt}
				\directm 	&&&& \multicolumn{2}{l|}{\parbox{7.5cm}{Uses the median value $f_{\rm median}$ in \eqref{eqn:potOptRect3}.}} &   \\[6pt]
				\cline{0-0}\cline{5-6}
				\vspace{5pt}
				\directz 	&&&& \multicolumn{2}{l|}{\parbox{7.5cm}{Uses the infinity norm in \eqref{eq:distance} and selects at most one POH from each group having the same measure $(\delta^k_i)$}} &  \\[9pt]
				\cline{0-1}\cline{5-6}
				\vspace{5pt}
				\directrev 	& \vspace{0.05pt} \parbox{2.45cm}{Hyper-rectangular partitions-based on 1-dimensional trisection}    &&& \multicolumn{2}{l|}{\parbox{7.5cm}{Selects at the most one POH from each group of $(\delta^k_i)$. Additional minimization procedure \fmincon{} is employed.}} &  \\[30pt]
				\cline{0-1}\cline{5-6}
				\vspace{5pt}
				\directav 	&&&& \multicolumn{2}{l|}{\parbox{7.5cm}{Uses the average value $f_{\rm average}$ in \cref{eqn:potOptRect4}.}} &  \\[6pt]
				\cline{0-0}\cline{5-6}
				\vspace{5pt}
				\dirmin 	& \multirow{8}{*}{\rotatebox{90}{ \parbox{4.5cm}{Hyper-rectangular partitions-based on $n$-dimensional trisection}}} &&& \multicolumn{2}{l|}{\parbox{7.5cm}{\fmincon{} is performed from each selected POH.}} &  \\[6pt]
				\cline{0-0}\cline{5-6}
				\vspace{5pt}
				\plor 		&&&& \multicolumn{2}{l|}{\parbox{7.5cm}{The set of POH is reduced to just two: with the maximal $(\delta^k_{\rm max})$ and the minimal $(\delta^k_{\rm min})$ measures.}} &  \\[6pt]
				\cline{0-0}\cline{5-6}
				\glbsolve &  &  &  &  &  &  \\[2pt]
				\cline{0-0}\cline{6-6}
				\texttt{glbSolve-} \texttt{sym} &  &  & &  & \multicolumn{1}{l|}{\multirow{2}{*}{\parbox{6cm}{Discards unnecessary hyper-rectangles for symmetric functions.}}} &  \\[2pt]
				\cline{0-0}
				\texttt{glbSolve-} \texttt{sym2} & &  & & \multicolumn{1}{l|}{\multirow{11}{*}{\rotatebox{90}{\parbox{8cm}{Uses the function \textit{conhull} to return all points on the convex hull, even redundant ones. \textit{conhull} is based on the extended version of the {\scriptsize GRAHAMSHULL} \cite{Preparata1985:book}  algorithm.}}}} & & \\[2pt]
				\cline{0-0}\cline{6-6}
				\vspace{5pt}
				\directmr &  &  & &  & \multicolumn{1}{l|}{\multirow{2}{*}{\parbox{6cm}{Performs the selection of POH on three different sets (\textquotedblleft levels\textquotedblright). \directmr$_{\scriptsize 075}$  uses different $\varepsilon$ values at each level.}}} &  \\[2pt]
				\cline{0-0}
				\directmr$_{\scriptsize 075}$ & &  & &  & & \\[2pt]
				\cline{0-0}\cline{6-6}

				\vspace{4pt}
				\gbdirect & &  & &  & \multirow{3}{*}{\parbox{6cm}{Uses an adaptive scheme for balancing the local and global search.\\ Local minimization procedure \fmincon{} is embedded into the \birmin{} algorithm.}} & \\[2pt]
				\cline{0-2}
				\birectgb & \multirow{4}{*}{\rotatebox{90}{\parbox{3.3cm}{Hyper-rectangular partitions-based on 1-dimensional bisection}}} & \multirow{3}{\textwidth}{\parbox{2.2cm}{Samples hyper-rectangle at two points lying on diagonals}} & &  & & \\[4pt]
				\cline{0-0}
				\birmin & &  & &  & & \\[4pt]
				\cline{0-0}\cline{6-6}
				\birect & &  & &  & & \\[4pt]
				
				\cline{0-0}\cline{3-3}\cline{6-6}
				\adc & & \multirow{1}{\textwidth}{\parbox{2.2cm}{Sample hyper-rectangle at two vertices}} & &  & \parbox{6cm}{Uses an adaptive scheme for balancing the local and global search.} & \\[24pt]
				
				\cline{0-2}\cline{6-6}
				\disimplc & \multirow{4}{*}{\rotatebox{90}{\parbox{2.7cm}{Simplicial partitions-based on $n$-dimensional trisection}}} & \multirow{2}{\textwidth}{\parbox{2.2cm}{Samples midpoints of the simplices}} & &  & & \\[2pt]
				\cline{0-0}\cline{6-6}
				\disimplcgb & &  & &  & \multicolumn{1}{l|}{\multirow{2}{*}{\parbox{6cm}{Uses an adaptive scheme for balancing the local and global search.}}} & \\[4pt]
				\cline{0-0}\cline{3-3}
				\disimplvgb & & \multirow{2}{\textwidth}{\parbox{2.2cm}{Samples at vertices of the simplices}} & &  & & \\[4pt]
				\cline{0-0}\cline{6-6}
				\disimplv & & & &  & & \\[4pt]
				\bottomrule
		\end{tabular}}
	\end{table}

	\subsection{\direct-type algorithms for generally constrained global optimization}
	\label{sec:general}

	The original \direct{} algorithm~\cite{Jones1993} only solves optimization problems with the variables' bounds.
	In this subsection, we consider a generally constrained global optimization problem of the form:
	\begin{equation}
		\label{eq:opt-problem}
		\begin{aligned}
			\min_{\mathbf{x} \in D} & \; f(\mathbf{x})\\
			\st \; & \mathbf{g} (\mathbf{x}) \leq \mathbf{0},\\
			& \mathbf{h} (\mathbf{x}) = \mathbf{0},
		\end{aligned}
	\end{equation}
	where $ f: \mathbb{R}^n \rightarrow \mathbb{R}$, $\mathbf{g}: \mathbb{R}^n \rightarrow \mathbb{R}^m$,  $\mathbf{h}: \mathbb{R}^n \rightarrow \mathbb{R}^r $ are (possibly non-linear) continuous functions.
	The feasible region is a non-empty set, consisting of points that satisfy all constraints, i.e., $D^{\rm feas} = D \cap \Omega \neq \emptyset$, where $\Omega = \{ \mathbf{x} \in \mathbb{R}^n: \mathbf{g} (\mathbf{x}) \leq \mathbf{0}, \mathbf{h} (\mathbf{x}) = \mathbf{0} \} $.
	As for the box-constrained problems, it is also assumed that the objective and all constraint functions are Lipschitz-continuous (with unknown Lipschitz constants) but can be non-linear, non-differentiable, non-convex, and multi-modal.

	The first \direct-type algorithm for problems with general constraints was introduced in~\cite{Jones2001}.
	Finkel in~\cite{Finkel2005} investigated three different constraint handling schemes within the \direct{} framework.
	The comparison revealed various disadvantages of the initial proposals.
	Recently, various new promising extensions for general global optimization problems were introduced (see, e.g.,~\cite{Basudhar2012,Costa2017,Liu2017,pillo2016,pillo2010,Stripinis2018b} and the references given therein).
	Below we briefly review approaches implemented in the current release of the \toolbox{} toolbox~(see \cref{tab:direct_classification}).

	An exact L1 penalty approach \directll~\cite{Finkel2004} is transforming the original constrained problem~\eqref{eq:opt-problem} in the form:
	\begin{equation}
		\label{eq:opt-problem-l1}
		\begin{aligned}
			& \min_{\mathbf{x}\in D} f(\mathbf{x})+\sum_{i=1}^{m} \max\{\gamma_{i} g_{i}(\mathbf{x}),0\}+\sum_{i=1}^{r} \gamma_{i+m} |h_{i}(\mathbf{x})|,
		\end{aligned}
	\end{equation}
	where $\gamma_{i}$ are penalty parameters.
	Experiments in~\cite{Finkel2005} showed promising results of this approach.
	Nevertheless, the biggest drawback is the users' requirement to set penalty parameters for each constraint function manually.
  In practice, choosing penalty parameters is an essential task and can significantly impact the algorithm's performance~\cite{Finkel2005,Liu2017,Paulavicius2014:book,Paulavicius2016:ol,Stripinis2018b}.

	In~\cite{Stripinis2018b}, we have introduced a new \direct-type extension based on the \directgl~\cite{Stripinis2018a} algorithm.
	The new \directce{} algorithm uses an auxiliary function approach that combines objective and constraint functions and does not require penalty parameters.
	The \directce{} algorithm works in two phases, where during the first phase, the algorithm finds feasible points and in the second phase improves a feasible solution.
	A separate step for handling infeasible initial points is beneficial when the feasible region is small compared to the entire search space.
	In the first phase, \directce{} samples the search space and minimizes the sum of constraint violations, i.e.:
	\begin{equation}
		\label{eq:constr-violation}
		\min_{\mathbf{x}\in D}\varphi(\mathbf{x}),
	\end{equation}
	where
	\begin{equation}
		\label{eq:cons_viol}
		\varphi(\mathbf{x})=\sum_{i=1}^{m} \max\{g_{i}(\mathbf{x}),0\}+\sum_{i=1}^{r} |h_{i}(\mathbf{x})|.
	\end{equation}
	The algorithm works in this phase until at least one feasible point $(\mathbf{x}\in D^{\rm feas}_{\varepsilon_{\varphi}})$ is found, where
	\begin{equation}
		\label{eq:dfeas}
		D^{\rm feas}_{\varepsilon_{\varphi}} = \{\mathbf{x} : 0 \leq \varphi(\mathbf{x}) \leq \varepsilon_{\varphi}, \mathbf{x} \in D \}.
	\end{equation}
	The $\varepsilon_{\varphi}$ is a small user-specified tolerance for the sum of constraint functions \cref{eq:cons_viol}.
	When feasible points are located, the effort is switched to improve the feasible solutions.
	In the second phase, \directce{} uses the transformed problem \eqref{eq:opt-problem}:
	\begin{equation}
		\label{eq:directce}
		\begin{aligned}
			& \min_{\mathbf{x}\in D} f(\mathbf{x}) + \tilde{\xi}(\mathbf{x},f_{\rm min}^{\rm feas}),\\
			& \tilde{\xi}(\mathbf{x},f_{\rm min}^{\rm feas}) =
			\begin{cases}
				0, & \mathbf{x} \in D^{\rm feas}_{\varepsilon_{\varphi}} \\
				0, & \mathbf{x} \in D_{\varepsilon_{\rm cons}}^{\rm inf}\\
				\begin{array}{l}
					\varphi(\mathbf{x}) + \Delta,
				\end{array}
				& \text{otherwise,}\\
			\end{cases}
		\end{aligned}
	\end{equation}
	where
	\begin{equation}
		D_{\varepsilon_{\rm cons}}^{\rm inf} = \{\mathbf{x} : f(\mathbf{x}) \leq f_{\rm min}^{\rm feas}, \varepsilon_{\varphi} < \varphi(\mathbf{x}) \leq \varepsilon_{\rm cons}, \mathbf{x} \in D \},
	\end{equation}
	and $\varepsilon_{\rm cons}$ is a small tolerance for constraint function sum, which automatically varies during the optimization process.
	An auxiliary function $\xi(\mathbf{x},f_{\rm min}^{\rm feas})$ depends on the sum of the constraint functions and the parameter $\Delta = |f(\mathbf{x}) - f_{\rm min}^{\rm feas}|$, equal to the absolute difference between the best feasible function value found so far $(f_{\rm min}^{\rm feas})$ and the objective value at an infeasible center point.
	The purpose of the parameter $\Delta$ is to forbid the convergence to infeasible regions by penalizing the objective value at infeasible points.
	In such a way, the formulation \eqref{eq:directce} does not require any penalty parameters and determines the convergence of the algorithm to a feasible solution.
	The value of $\xi(\mathbf{x},f_{\rm min}^{\rm feas})$ is updated when a smaller value of $f_{\rm min}^{\rm feas}$ is found.
	This way, the new \directce{} algorithm divides more hyper-rectangles with center points lying close to the boundaries of the feasible region, i.e., the potential solution.

	The proposed \directce{} algorithm has two extensions:
	The first one is \directc~(see \cref{tab:direct_classification}), which is a simplified version of \directce{} and instead of  \eqref{eq:directce} minimizes the transformed problem:
	\begin{equation}
		\label{eq:directc}
		\begin{aligned}
			& \min_{\mathbf{x}\in D} f(\mathbf{x}) + \xi(\mathbf{x},f_{\rm min}^{\rm feas}),\\
			& \xi(\mathbf{x},f_{\rm min}^{\rm feas})=
			\begin{cases}
				0, & \mathbf{x} \in D^{\rm feas}_{\varepsilon_{\varphi}} \\
				\begin{array}{l}
					\varphi(\mathbf{x}) + \Delta,
				\end{array}
				& \text{otherwise,}\\
			\end{cases}
		\end{aligned}
	\end{equation}
	Experimental investigation in~\cite{Stripinis2018b} showed that the algorithm has the most wins in this comparison and can solve about $50\%$ of the problems with the highest efficiency.
	Unfortunately, the \directc{} algorithm's efficiency decreases, solving more challenging problems (with non-linear constraints and $n \ge 4$), where \directce{} is significantly better.
  Therefore, \directc{} should be used only for simpler optimization problems (with linear constraints and $n \leq 4$).
	The second extension of the \directce{} algorithm is \directcemin, where the algorithm is incorporated with \matlab{} optimization solver \fmincon{}.
	In~\cite{Stripinis2018b}, we observed that embedding a local minimization procedure into \directcemin~(see \cref{tab:direct_classification}) significantly reduces the total number of function evaluations compared to \directce{} and can significantly improve the quality of the final solution.

	\subsubsection{\direct-type algorithms for linearly constrained global optimization}
	\label{sec:linear}
	Let us note that all previously described algorithms for a generally constrained problem can be directly applied to solve linearly constrained problems.
	In this section, we consider optimization problems with only linear constraints.

	In~\cite{Paulavicius2016:ol}, we have extended the original simplicial partitioning-based \texttt{DISIMPL} algorithm~\cite{Paulavicius2013:jogo,Paulavicius2014:book} for such problems with linear constraints.
	Simplices may cover a search space defined by linear constraints.
	Therefore, a simplicial approach may tackle such linear constraints very subtly.
	In such a way, the new algorithms (\disimpllc{} and ~\disimpllv)~\cite{Paulavicius2016:ol} perform the search only in the feasible region, in contrast to other \direct-type approaches.
	Nevertheless, the authors in~\cite{Paulavicius2016:ol} showed that the feasible region's calculation requires solving $2n+m$ linear $n$-dimensional systems, and such operation is exponential in complexity.
	Therefore, the proposed algorithm can be effectively used for relatively small $n$ and $m$ values.

	\subsection{\direct-type algorithms for problems with hidden constraints}
	\label{sec:hidden}

	Optimization problems with hidden constraints often occur when the objective function is not defined everywhere~\cite{Chen2016}.
	Typical examples of such situations are the simulation crashes~\cite{Digabel2015} and failure of computations within the objective function~\cite{Carter2001,Choi2000,David1996,Stoneking1992}.
	As in~\cite{Chen2016,Digabel2015}, we call these internal to $f$ constraints ``hidden constraints'' and assume that $f$ fails to return a value when evaluated at $\mathbf{x} \notin D^{\rm feas}$.
	Some authors alternatively may use other terms like ``crash,'' ``unknown,'' ``unspecified,'' and ``forgotten'' constraints~\cite{Bachoc2020,Digabel2015}.

	In this subsection, we consider the solution to the constrained global optimization problem:
	\begin{equation}
		\label{eq:opt-problem3}
		\begin{aligned}
			& \min_{\mathbf{x}\in D^{\rm feas}} && f(\mathbf{x}), 
		\end{aligned}
	\end{equation}
	where $f:\mathbb{R}^n \rightarrow \mathbb{R} \cup \{ \infty \}$ denotes an extended real-valued, most likely ``black-box'' objective function.
	A priori an unknown feasible region $D^{\rm feas}$ is defined as a non-empty set
	$$
	D^{\rm feas} = D \setminus D^{\rm hidden} \neq \emptyset,
	$$
	and $D^{\rm hidden}$ are not given by explicit formulae hidden constraints.
	Such a problem formulation leads to a complex and analytically undefined feasible region.
	Hidden constraints are typically handled by returning NaN or $\infty$ evaluating the objective function at $\mathbf{x} \notin D^{\rm feas}$.
	Therefore, using NaN hyper-rectangles with an infeasible center point would not be selected as potential optimal at all \cite{Stripinis2021}.
	In the case of $\infty$ (or any other high value), they would be left unexplored as long as there are the same size hyper-rectangles with feasible centers \cite{Finkel2005, Stripinis2021}.
	Unfortunately, most \direct-type algorithms cannot be directly (without any modifications) applied for the problem's \eqref{eq:opt-problem3} solution.
	In the current version of the \toolbox{} toolbox, four such algorithms are available (see \cref{tab:direct_classification}).


	One of the first proposed modifications for such problems was the barrier method (\directbarrier)~\cite{Gablonsky2001:phd}.
	The \directbarrier{} is relatively straightforward and assigns a predefined high value to infeasible hyper-rectangles.
	However, such an approach produces other well-known problems discussed and reviewed by a few authors~\cite{Finkel2005,Paulavicius2014:book,Stripinis2018b}.
	The main issue is that the barrier approach makes exploration around the edges of feasibility very slow.
	Significant penalties used by the barrier method ensure that no infeasible hyper-rectangle can be potentially optimal as long as there is the same measure hyper-rectangle with the feasible center midpoint.
	For \directbarrier, the priority is the examination of regions where feasible points are found already.
	Another critical issue concluded in~\cite{Finkel2005} is that hyper-rectangles, even with the sizeable feasible region, will not be explored in a reasonable number of function evaluations.
	To sum up, the barrier approach is not the best fit for the problem \eqref{eq:opt-problem3}.

	The second \direct-type approach for hidden constraints is based on Neighbourhood Assignment Strategy (NAS)~\cite{Gablonsky2001:phd}.
	\directnas's main idea is to assign the value at infeasible point $\mathbf{x}^{\rm inf} \notin D^{\rm feas}$ relative to the objective values attained in the feasible points from the neighborhood of $\mathbf{x}^{\rm inf}$.
	\directnas{} iterates over all infeasible midpoints by creating surrounding hyper-rectangles around them by keeping the same center points in every iteration.
	These hyper-rectangles are increased by doubling the length of each dimension.
	If more than one feasible center point inside the enlarged region, \directnas{} assigns the smallest function value to the infeasible midpoint plus a small epsilon $f(\mathbf{x}^{\rm feas}) + \epsilon f(\mathbf{x}^{\rm feas})$, where $\epsilon = 10^{-6}$ was proposed to use.
	If inside the enlarged region has no feasible points, \directnas{} assigns the largest objective function value found so far $f_{\rm max} + \lambda$, where $\lambda = 1$ was proposed to use.
	This strategy does not allow the \directnas{} algorithm to move beyond the feasible region by penalizing infeasible midpoints with large values.
	However, the algorithm's principal concern is the slow convergence caused by many additional calculations.

	Another recent idea to handle hidden constraints within the \direct{} framework is to use a subdividing step for infeasible hyper-rectangles.
	The proposed \directsub~\cite{Na2017} incorporates the previously mentioned barrier approach techniques.
	Specifically, if the center point is identified as infeasible, then \directsub{} assigns a considerable penalty value to it.
	An extra subdividing step is performed only in specific iterations, during which all infeasible hyper-rectangles are identified as potentially optimal and subdivided together with others POHs.
	The sub-dividing step can decompose the boundaries of the hidden constraints quite efficiently. 
	Still, \directsub{} has several apparent drawbacks.
	The algorithm performance depends on when (how often) the subdividing step is performed.
	Therefore, new subdivisions can grow drastically, especially for higher dimensionality problems.

	The most recent version for hidden constraints \directh~\cite{Stripinis2021} is based on our previous \directgl~\cite{Stripinis2018a} algorithm.
	For hyper-rectangles with infeasible midpoints, \directh{} assigns a value depending on how far the center is from the current best minima $\mathbf{x}_{\rm min}$.
	Such a technique does not require any additional computation.
	Simultaneously, distances from the $\mathbf{x}_{\rm min}$  point are already known, as they are used to selecting potential optimal hyper-rectangle schemes adapted from \directgl~\cite{Stripinis2018a}.
	In such a way, \directh{} does not penalize infeasible hyper-rectangles with large values (as was suggested by previous proposals), which are close to the $\mathbf{x}_{\rm min}$ and assure a faster and more comprehensive examination of hidden regions.
	Moreover, this approach employs additional procedures to efficiently handle infeasible initial points (see~\cite{Stripinis2021} for experimental justification).

	\subsection{Implementation of \direct-type algorithms within \toolbox}
	\label{sec:parralell}

	\subsubsection{Sequential implementation of the algorithms}
	\label{sec:sequantial}

	The performance of \direct-type algorithms highly depends on computer implementation.
	Most publicly available \direct{} implementations (see, e.g., \directfinkel~\cite{Finkel2004} and \glbsolve~\cite{Bjorkman1999}) use static data memory management \cite{Bjorkman1999,Finkel2004,Gablonsky2001:phd}.
  In the ``Impelmentation'' column of \cref{tab:direct_classification}, we provide information on which data structures were used in our implementations and whether a particular algorithm was implemented in parallel.
  Below we look at the main advantages and disadvantages of each of them.

	With static data management, all information received after the domain partitioning is stored in the contiguous memory blocks.
	This includes objective and constraint function values, index numbers, center point coordinates, side lengths of hyper-rectangles, and so on.
  Such implementation can quickly access the elements for further selection, sampling, and subdivision steps.
  An apparent drawback of the static data structure is unpredictable memory demand due to different characteristics of the optimization problems.
  Thus many \direct-type algorithmic implementations use large static arrays to store the current state of the space partitioning.
  If any array is insufficient to store the required information, this can lead to code failure.

  Another disadvantage of static data structures used in \direct{} implementations is that they require unnecessary recalculations in each iteration.
  One of the essential tasks in the \direct-type algorithms is the selection step.
  This step requires sorting all existing hyper-rectangles by the same size of diameter.
  Such sorting becomes especially inefficient when the optimization process is longer and the amount of data gets large, e.g., for higher dimensionality problems or when a solution with high accuracy is required.

  In \cite{He2002}, the authors proposed using dynamic data structures. Information received after space partitioning is sorted by hyper-rectangle diameters and stored in columns.
  All rectangles of the same diameter are stored in the column in any order.
  In \cite{He2002}, the authors mentioned the idea of sorting columns by function values in descending order or inserting all new data in sorted sequences separately.
  However, any of these ideas have not been investigated further.
  With dynamic data structures, the selection step is much more efficient.
  It can be performed only in the set consisting of the best function values from each column.
  Such implementation saves lots of time compared with the static data structure-based implementation.

  In \cite{Stripinis2020}, we have compared two different implementations (static and dynamic) of the same \directce{} algorithm.
  The dynamic implementation of the code required, on average, $62 \%$ less total execution time than static-based.
  The difference was even more significant when the number of function evaluations was high.

One of the apparent drawbacks of the dynamic data structure is unpredictable columns size.
Fully processed POHs must be removed from the previous columns and added to a new/existing column.
During the algorithm's execution, there can be many hyper-rectangle diameters.
Depending on the dimension of the problem, usually, the initial array is allocated of reasonably large size.
If the array provides insufficient size, new blocks of columns will be reallocated as needed.
In practice, only a few of these columns need reallocation at any given time.

	\subsubsection{Parallel implementations of the algorithms}
	\label{sec:parral}

	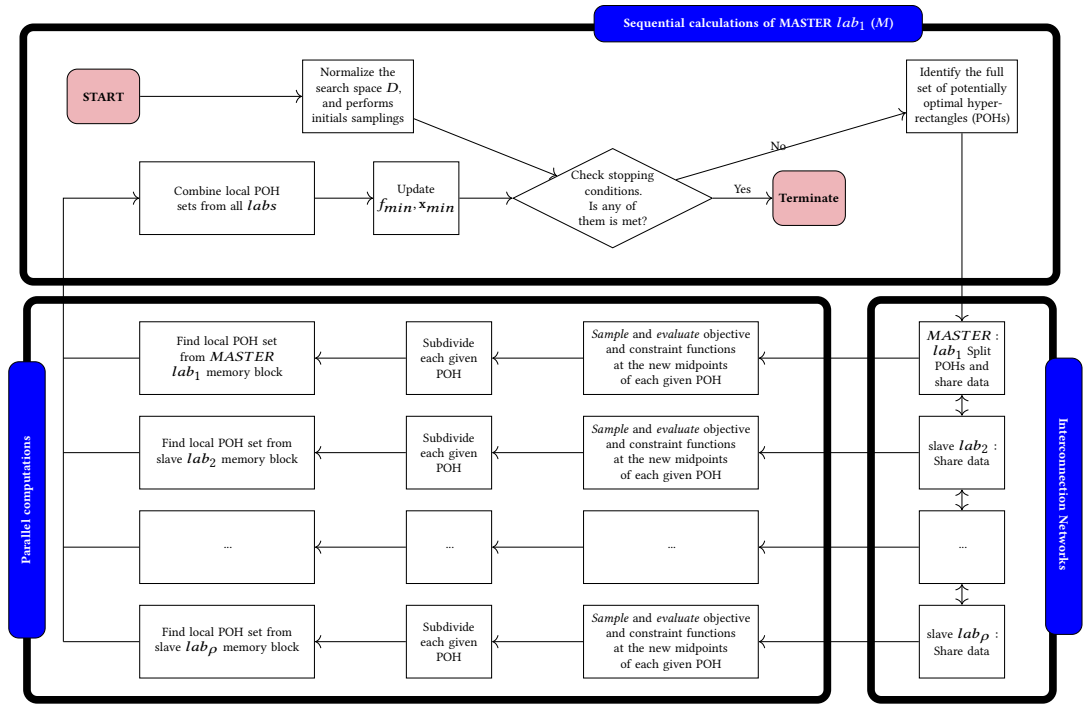
\begin{figure}[ht]
		\centering
		\resizebox{.95\textwidth}{!}{%
			\begin{tikzpicture}[node distance=1.5cm]
				\tiny
				\tikzstyle{envirom} = [rectangle, rounded corners, text width=5em, minimum width=2.5cm, minimum height=5.5cm,text centered, draw=black, mark size=1.5pt, line width=3pt]
				\tikzstyle{envipar} = [rectangle, rounded corners, text width=5em, minimum width=11cm, minimum height=5.5cm,text centered, draw=black, mark size=1.5pt, line width=3pt]
				\tikzstyle{envigpu} = [rectangle, rounded corners, text width=5em, minimum width=2cm, minimum height=6cm,text centered, draw=black, mark size=1.5pt, line width=3pt]
				\tikzstyle{proparalvertical}   = [rectangle, rounded corners, text width=2em, minimum width=0.5cm, minimum height=3.8cm,text centered, draw=black, fill=blue]
				\tikzstyle{infor} = [rectangle, rounded corners, text width=6em, minimum width=1cm, minimum height=0.5cm,text centered, draw=black, fill=blue]
				\tikzstyle{infos} = [rectangle, rounded corners, text width=25em, minimum width=1cm, minimum height=0.5cm,text centered, draw=black, fill=blue]
				\tikzstyle{sequan} = [rectangle, rounded corners, text width=5em, minimum width=14.2cm, minimum height=3.5cm,text centered, draw=black, mark size=1.5pt, line width=3pt]
				\tikzstyle{startstop} = [rectangle, rounded corners, text width=5em, minimum width=0.75cm, minimum height=0.75cm,text centered, draw=black, fill=venetianred!30]
				\tikzstyle{info} = [rectangle, rounded corners, text width=15em, minimum width=1cm, minimum height=0.5cm,text centered, draw=black, fill=blue]
				\tikzstyle{process}   = [rectangle, text width=6em, minimum width=1cm, minimum height=1cm, text centered, draw=black]
				\tikzstyle{proparal}   = [rectangle, text width=13em, minimum width=1cm, minimum height=1cm, text centered, draw=black]
				\tikzstyle{processas}   = [rectangle, text width=8em, minimum width=1cm, minimum height=1cm, text centered, draw=black]
				\tikzstyle{decision} = [diamond, aspect=2, inner sep=0pt, text width=7em, text centered, draw=black]

				\node (sequa)     [sequan] {};
				\node (sequand)   [infos, below of=sequa, yshift=3.3cm, xshift=3cm] {\textcolor{white}{\textbf{Sequential calculations of MASTER $lab_1$ ($M$)}}};
				\node (netwr)     [envirom, below of=sequa, xshift=5.8cm, yshift=-3.25cm] { };
				\node (netwre)    [envipar, left of=netwr, xshift=-5.85cm] { };
				\node (parpar)    [proparalvertical, above of=netwre, yshift=-1.5cm, xshift=-5.5cm] {\begin{turn}{-270}\textcolor{white}{\textbf{Parallel computations}}\end{turn}};
				\node (start)     [startstop, right of=sequa, xshift=-7.5cm, yshift=0.8cm] { \textbf{START} };
				\node (normalize) [processas, right of=start, xshift=2cm] { Normalize the search space $D$, and performs initials samplings};
				\node (seletio)   [processas, right of=normalize, xshift=6.8cm] { Identify the full set of potentially optimal hyper-rectangles (POHs)};

				\node (inte1)     [process, below of=seletio, yshift=-2.1cm] {$MASTER:$ $lab_1$ Split POHs and share data};
				\node (inte2)     [process, below of=inte1, yshift=0.2cm] {slave $lab_2:$ Share data};
				\node (inte3)     [process, below of=inte2, yshift=0.2cm] { ... };
				\node (inte4)     [process, below of=inte3, yshift=0.2cm] {slave $lab_\rho:$ Share data};
				\node (parpar)    [proparalvertical,right of=inte2, xshift=-0.1cm, yshift=-0.6cm] {\begin{turn}{-90}\textcolor{white}{\textbf{Interconnection Networks }}\end{turn}};
				\node (work1)     [proparal, left of=inte1, xshift=-2.5cm] {\textit{Sample} and \textit{evaluate} objective and constraint functions at the new midpoints of each given POH};
				\node (work2)     [proparal, left of=inte2, xshift=-2.5cm] {\textit{Sample} and \textit{evaluate} objective and constraint functions at the new midpoints of each given POH};
				\node (work3)     [proparal, left of=inte3, xshift=-2.5cm] {...};
				\node (work4)     [proparal, left of=inte4, xshift=-2.5cm] {\textit{Sample} and \textit{evaluate} objective and constraint functions at the new midpoints of each given POH};
				\node (hidde1)     [process, left of=work1, xshift=-1.55cm] {Subdivide each given POH};
				\node (hidde2)     [process, left of=work2, xshift=-1.55cm] {Subdivide each given POH};
				\node (hidde3)     [process, left of=work3, xshift=-1.55cm] {...};
				\node (hidde4)     [process, left of=work4, xshift=-1.55cm] {Subdivide each given POH};
				\node (poh1)      [proparal, left of=hidde1, xshift=-1.55cm] {Find local POH set from $MASTER$ $lab_1$ memory block};
				\node (poh2)      [proparal, left of=hidde2, xshift=-1.55cm] {Find local POH set from slave $lab_2$ memory block };
				\node (poh3)      [proparal, left of=hidde3, xshift=-1.55cm] {...};
				\node (poh4)      [proparal, left of=hidde4, xshift=-1.55cm] {Find local POH set from slave $lab_\rho$ memory block};
				\node (gather)    [proparal, above of=poh1, yshift=0.7cm] {Combine local POH sets from all $labs$};
				\node (update)    [process, right of=gather, xshift=1.1cm] { Update $f_{min}, \mathbf{x}_{min}$ };
				\node (stop)      [decision, right of=update, xshift=1.2cm] {Check stopping conditions. Is any of them is met?};
				\node (termin)    [startstop, right of=stop, xshift=1.2cm] { \textbf{Terminate} };
				\node (ii1)       [coordinate, left of=poh1, xshift=-0.75cm] {};
				\node (ii2)       [coordinate, left of=poh2, xshift=-0.75cm] {};
				\node (ii3)       [coordinate, left of=poh3, xshift=-0.75cm] {};
				\node (ii4)       [coordinate, left of=poh4, xshift=-0.75cm] {};

				\draw[->] (start)     -- (normalize);
				\draw[->] (normalize) -- (stop);
				\draw[->] (seletio)   -- (inte1);
				\draw[->]  (inte1)    -- (work1);
				\draw[->]  (inte2)    -- (work2);
				\draw[->]  (inte3)    -- (work3);
				\draw[->]  (inte4)    -- (work4);
				\draw[->]  (work1)    -- (hidde1);
				\draw[->]  (work2)    -- (hidde2);
				\draw[->]  (work3)    -- (hidde3);
				\draw[->]  (work4)    -- (hidde4);
				\draw[->]  (hidde1)   -- (poh1);
				\draw[->]  (hidde2)   -- (poh2);
				\draw[->]  (hidde3)   -- (poh3);
				\draw[->]  (hidde4)   -- (poh4);
				\draw[-]   (ii1)      -- (poh1);
				\draw[-]   (ii2)      -- (poh2);
				\draw[-]   (ii3)      -- (poh3);
				\draw[-]   (ii4)      -- (poh4);
				\draw[->]  (ii4)      |- (gather);
				\draw[->]  (gather)   -- (update);
				\draw[->]  (update)   -- (stop);
				\draw[<->]  (inte1)    -- (inte2);
				\draw[<->]  (inte2)    -- (inte3);
				\draw[<->]  (inte3)    -- (inte4);
				\draw[->] (stop)      -- node[anchor=east] {No} (seletio);
				\draw[->] (stop)      -- node[anchor=south] {Yes} (termin);
			\end{tikzpicture}
		}
		\caption{Flowchart diagram for the parallel implementations of selected \direct-type algorithms.}
		\label{fig:par}
	\end{figure}

	We use the  MathWorks official extension to the \matlab{} language~-- the Parallel Computing Toolbox~\cite{Matlab2020} for parallel implementations.
	The Parallel Computing Toolbox provides several parallel programming paradigms~\cite{Luszczek2009}, like threads, parallel \texttt{for}-loops, and \texttt{SPMD} (Single Program Multiple Data).
	In~\cite{Stripinis2020}, we concluded that the \texttt{SPMD}-based parallel implementation of the \directce{} is the most efficient and significantly outperforms the other two based on \texttt{parfor}-loops.



	Therefore, in the \toolbox{} toolbox, parallel implementations are based on the \texttt{SPMD} functionality within the Parallel Computing Toolbox, used to allocate the work across multiple labs in the \matlab{} software environment.
	Each lab stores information on its main memory block, and data is exchanged through the message passing over the interconnection network~\cite{Matlab2020}.
	The master-slave paradigm is used to implement dynamic load balancing.
	The flowchart of the parallel algorithmic framework is illustrated in \cref{fig:par}.
	One lab is the master, denoted by $lab_1$, and the other labs are slaves $lab_{i}, i = 2, \dots, \rho$. The master also acts as a slave.
	Each iteration must be done in a sequence to preserve the determinism.

  The master performs the following tasks:
\begin{itemize}
  \itemsep=0cm
  \item \textit{The initialization step}: normalizes the domain $(D)$ and evaluates the objective and constraint functions at the center point. Here, only the optimization problem and the information about the domain $D$ are shared with slaves $lab_{i}, i = 2, \dots, \rho$.
  \item \textit{At each iteration}:
  \begin{itemize}
    \item checks the stopping conditions and informs the slaves if any of them have been met.
    \item finds the {full} set of POHs by performing a selection step considering the combined set of {local} POHs.
    \item splits the {full} set of POHs  among all slaves and itself equally.
    \item gives instructions to the slaves having an excess of POHs (in their local memory) to share them with those who have a deficit, including itself.
    \item sends or receives POHs according to its instructions.
    \item performs the sampling, subdivision, and local selection steps as the slave.
    \item receives from slaves the information about their local POHs sets.
  \end{itemize}
\end{itemize}
\noindent The slaves perform the following tasks:
\begin{itemize}
  \item \textit{At each iteration}:
  \begin{itemize}
    \itemsep=0cm
    \item Send or receive POHs, according to the master's instructions.
    \item Perform the sampling and subdivision steps sequentially.
    \item Perform the selection using the information in their local memory and send local POHs to the master.
    \item Terminate when such an instruction from the master is received.
  \end{itemize}
\end{itemize}

	The master lab decides which hyper-rectangles will be sampled and subdivided and how these tasks will be distributed among all available slave labs.
	Additionally, the master lab is responsible for stopping the algorithm.
	The master lab also performs load balancing by distributing the selected hyper-rectangles to the rest of the slave labs.
	When the slave labs $(lab_i, i=1, \dots, \rho)$ receive tasks from the master lab, each sequentially performs the sampling and subdivision steps.
  Then finds a local set of POHs and sends local data back to the master lab for the further global selection step.
  After this, each slave becomes idle until further instructions are received.
	Suppose any of the termination conditions are satisfied.
	In that case, all slave labs receive the notification that the master lab has become inactive, and the slave labs will terminate themselves without further messaging.
	We refer to~\cite{Stripinis2020} for a more detailed description and analysis of parallel schemes.

  We should note that not all \direct-type implementations can use the latter scheme of parallelism.
  For example, implementations using the \textit{conhull} function, which returns all the points on a convex hull, cannot.
  To preserve the determinism, only the master should select POHs and have all data stored in its memory.
  The framework shown in \cref{fig:par} is inappropriate, and a new one should be developed.
  The \directnas{} algorithm has an additional expensive constraint handling step not addressed in the proposed parallel scheme.
  As summarized in \cref{tab:direct_classification}, currently, 17 \direct-type algorithms are implemented in parallel within the \toolbox{} toolbox.

	\section{\toolbox{} toolbox}
	\label{sec:tooblox}

	The sequential and parallel implementation of \direct-type algorithms presented in the previous sections forms the basis for our \toolbox{} toolbox.
	The toolbox consists of two main parts:
	\begin{itemize}
		\item \toolboxa{} - \matlab{} toolbox package containing implementations of \direct-type algorithms (from \cref{tab:direct_classification}), including an extensive \directlib{} library of the box and generally constrained test and practical engineering global optimization problems, often used for benchmarking \direct-type algorithms.
		\item \toolboxb{} - A single \matlab{} app installer containing everything necessary to install and run the \toolbox{} toolbox, including a graphical user interface (GUI).
	\end{itemize}

	\subsection{Graphical user interface}
\begin{figure}[t]
  \centering
  \includegraphics[width=0.95\textwidth]{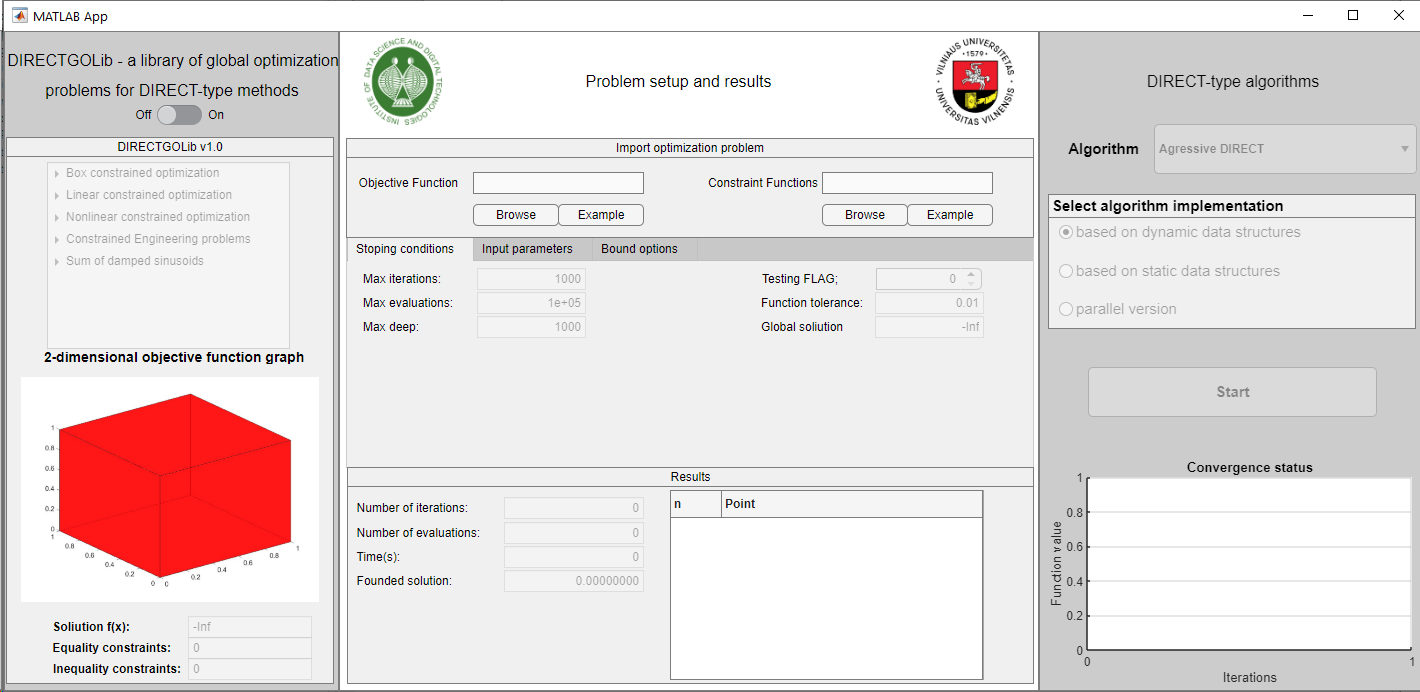}
  \caption{The graphical user interface (GUI) of \toolbox{} toolbox.}
  \label{fig:toolbox}
\end{figure}

	After installation (using \toolboxb{}), \toolbox{} can be launched from \matlab{} \textbf{APPS}, located in the toolbar.
	The graphical interface of the main \toolbox{} toolbox window is shown in \cref{fig:toolbox}.
	Application is divided into three main parts: i) selection of the problem's type from \directlib; ii) setting up an optimization problem and algorithmic options; iii) selection of \direct-type algorithm, his implementation, and convergence plot of obtained results.

	The first step is to specify the objective and constraint functions, loading them from the integrated \directlib{} library or selecting them from other sources.
	Examples of the structure needed are present.
	All test problems from the \directlib{} library have up to three key features: i) known globally optimal solutions, ii) a complete description of the problem, including objective and constraint functions (if any), and iii) problem visualization (only for two-dimensional problems).

	After selecting the optimization problem, the second step is to set up the bound constraints for each variable.
	First, the user needs to specify the algorithm and the type of implementation.
	Two implementations are based on different data structures (static and dynamic), and the third is a parallel version of the algorithm.
	For simplicity, some toolbox{} options are set to default values and not displayed in the GUI but can be changed in the toolbox settings.
	After the termination, the \textbf{Results} part displays the final solution and performance metrics.
	Additionally, the convergence process is shown in the \textbf{Convergence status} part.

\subsection{MATLAB toolbox}
\label{toolbox}

After installation of the \matlab{} toolbox (using \toolboxa), all implemented \direct-type algorithms and test problems can be freely accessed in the command window of \matlab.
Unlike using GUI, algorithms from the command line require more programming knowledge, and configurations must be done manually. All algorithms can be run using the same style and syntax:
\begin{minted}[bgcolor=bg]{matlab}
1. f_min = algorithm(P);
2. f_min = algorithm(P, OPTS);
3. f_min = algorithm(P, OPTS, D);
4. [f_min, x_min]  = algorithm(P, OPTS, D);
5. [f_min, x_min, history]  = algorithm(P, OPTS, D);
\end{minted}
The left side of the equations specifies the output parameters.
After the termination, the algorithm returns the best objective value (\texttt{f\_{min}}), solution point (\texttt{x\_{min}}), and history of the algorithmic performance during all iterations (\texttt{history}).
The information presented here is the iteration number, the total number of objective function evaluations, the current minimum value, and execution time.

The algorithm name (\texttt{algorithm}) and at least one input parameter are needed to specify on the right side.
The first one is the problem structure (\texttt{P}) consisting of an objective function:
\begin{minted}[bgcolor=bg]{matlab}
>> P.f = 'objfun';
\end{minted}
If the problem involves additional constraints, they also must be specified:
\begin{minted}[bgcolor=bg]{matlab}
>> P.constraint = 'confun';
\end{minted}
	The second parameter (\texttt{OPTS}) customizes the default algorithmic settings.
	The third parameter (\texttt{D}) is used to specify the bound constraints for each variable (see \cref{eq:opt-problem1}).


	\section{Experimental investigation of \direct-type algorithms in the \toolbox{} toolbox}
	\label{sec:experimentstests}

	This section presents the performance evaluation of \direct-type algorithms on test and engineering design problems from the \directlib{}~\cite{DIRECTGOLibv1.0,DIRECTGOLibZenodov1.0}.
	The \directlib{} library consists of the box and generally constrained test and practical engineering global optimization problems for various \direct-type algorithms benchmarking.
	Experimental results presented in this section are also available in digital form in the \texttt{Results/TOMS} directory of the Github repository \url{https://github.com/blockchain-group/DIRECTGO}~\cite{DGO2021}).
  The most recent version of \texttt{DIRECTGOLib v1.1}~\cite{DIRECTGOLibv1.1,DIRECTGOLibZenodov1.1} has a few additional test functions introduced in other parallel studies but has not been considered in this article.
  Despite this, our vision is to develop \texttt{DIRECTGOLib} further so that all algorithms can use the latest version of \texttt{DIRECTGOLib} without any additional preparations.


  We distinguish the following classes (types) of global optimization problems:
	\begin{itemize}
		\item \textbf{BC} -- \textbf{B}ox-\textbf{C}onstrained problems;
		\item \textbf{LC} -- \textbf{L}inearly-\textbf{C}onstrained problems;
		\item \textbf{GC} -- \textbf{G}enerally-\textbf{C}onstrained problems.
	\end{itemize}
	A summary of all optimization problems in \directlib{} and their properties is given in~\cref{sec:directlib}, \cref{tab:resultslow}.
	We note that some test problems have several variants, e.g., \textit{Bohachevsky}, \textit{Shekel}, and some of them, like \textit{Alpine}, \textit{Csendes}, and \textit{Griewank}, can be used by changing the problem's dimensionality.
  We used the following dimensions in our experimental setting: $n = 2, 5, 10, 15$.
  For some test problems, the second dimension $(n = 2)$ was skipped because the problem was then too easy to solve, and sometimes we skipped $n = 15$ because the resulting problem was too hard that none of the algorithms were able to solve it.
  The fourth column in \cref{tab:resultslow} indicates the exact dimensions used for all test problems.

	All computations were carried out on a 6-core computer with 8th Generation Intel R Core$^\textit{TM}$ i7-8750H @ 2.20GHz Processor, 16 GB of RAM, and \texttt{MATLAB R2020b}.
	Performance analysis was carried out using physical cores only and disabled hyper-threading.

	All global minima $ f^* $ are known for all test problems.
	Therefore, the investigated algorithms were stopped when it was generated such the point $ {\mathbf{x}} $ with whom the percent error
	\begin{equation}
		\label{eq:pe}
		\ pe = 100 \% \times
		\begin{cases}
			\frac{f({\mathbf{x}}) - f^*}{|f^*|},  & f^* \neq 0, \\
			f({\mathbf{x}}),  & f^* = 0,
		\end{cases}
	\end{equation}
	is smaller than the tolerance value $\varepsilon_{\rm pe}$, i.e., $pe \le \varepsilon_{\rm pe}$.
	Additionally, we stopped the tested algorithms when the number of function evaluations exceeded the prescribed maximal limit (equal to $2 \times 10^6$) or took more than $43,000.00$ seconds.
	In any of these situations, the final result is set to $2 \times 10^6$ to further process results.
	Two different values for $\varepsilon_{\rm pe}$ were considered: $10^{-2}$, $10^{-8}$.
	By default, algorithms were tested using the $\varepsilon_{\rm pe} = 10^{-2}$ value.
	Algorithms incorporating additional schemes to speed up the solution's refinement have been tested using the $\varepsilon_{\rm pe} = 10^{-8}$ value.

	Additionally, we analyze and compare the algorithms' performance by applying the data profiles~\cite{More2009} to the convergence test~\eqref{eq:pe}.
	The data profile is a popular and widely used tool for benchmarking and evaluating the performance of several algorithms (solvers) when run on a large problem set.
	Benchmark results are generated by running a certain algorithm $v$ (from a set of algorithms $\mathcal{V}$ under consideration) for each problem $u$ from a benchmark set $\mathcal{U}$ and recording the performance measure of interest.
  The performance measure could be, for example, the number of function evaluations, the computation time, the number of iterations, or the memory used.
	We used a number of function evaluations and the execution (computation) time criteria.

	The data profiles provide the percentage of problems that can be solved with a given budget of the desired performance measure.
	The data profile is defined
	\begin{equation}
		\lambda_{v}(\alpha) = \frac{1}{card(\mathcal{U})}\textrm{size} \left\{ u \in \mathcal{U} : t_{u,v} \le \alpha \right\},
		\label{eq:data-profile}
	\end{equation}
	where $t_{u,v} > 0$ is the number of performance measure required to solve problem $u$ by the algorithm $v$, and  $card(\mathcal{U})$ is the cardinality of $\mathcal{U}$.
	In our case, the $\lambda_{v}(\alpha)$ shows the percentage of problems that can be solved within $\alpha$ function evaluations, or seconds.

  	In the experimental studies, the developed \direct-type algorithms (within the \toolbox{} toolbox) were compared among themselves and with three \tomlab{} \direct-type solvers:
  	\begin{itemize}
  		\item \glbolve{} \cite{Holmstrom2010} -- implementation of the \direct{} algorithm \cite{Jones1993};
  		\item \glsolve{} \cite{Holmstrom2010} -- implementing an extended \direct{} version \cite{Jones1993,Jones2001} capable of handling linear and non-linear constrained problems;
  		\item \glccluster{} \cite{Holmstrom2010} -- implementation of the \direct{} algorithm \cite{Jones1993}, hybridized with local search subroutine and clustering techniques.
  	\end{itemize}

  We note that the \glccluster{} algorithm has a large number and different input parameters that may significantly impact the algorithm’s performance.
  Even a parameter such as the maximum allowed number of function evaluations can significantly impact the algorithm’s performance. Our aim was not to find the optimal parameters values, as this is a complex process, but to investigate how these algorithms compare using the default values (provided by \tomlab{} software developers).
  	When the algorithm reaches the default limit for the maximum number of function evaluations ($M_{\rm max}$) set in the default parameters, we restarted the algorithm using the final status from the previous run (``\textit{warm start}''~\cite{Holmstrom2010}) with doubled $2M_{\rm max}$.
  	In such a way, sometimes a default $M_{\rm max}$ value was doubled up to our maximal limit of evaluations $2 \times 10^6$ was reached.

  	The \texttt{Scripts/TOMS} directory of the Github repository (\url{https://github.com/blockchain-group/DIRECTGO}) provides four different \texttt{MATLAB} scripts for cycling through all different classes of test problems used in this paper.
  	The constructed scripts can be handy for reproducing the results presented here, as well as for comparison and evaluation of newly developed algorithms.


	\subsection{Comparison of \direct-type algorithms for box constrained optimization}
	\label{sec:box_constraints}

	\begin{table}

		\caption{The performance of \direct-type algorithms from \toolbox{} and \texttt{TOMLAB} based on the number of function evaluations ($f_{eval.}$),  the total execution time in seconds ($time$), and the total number of iterations ($iter.$) criteria on a set of box-constrained problems (from \directlib) using $\varepsilon_{\rm pe} = 10^{-2}$ in \eqref{eq:pe}.}
		\resizebox{\textwidth}{!}{
			\begin{tabular}[tb]{@{\extracolsep{\fill}}l|r|r|rrr|rrr|rrr|rrr}
				\toprule
				\multirow{2}{*}{Algorithm} & \textbf{Avg. \# local} & \multirow{2}{*}{\textbf{Failed}} & \multicolumn{3}{c|}{\textbf{Average results}}  & \multicolumn{3}{c|}{\textbf{Average results} $(n \le 4)$} & \multicolumn{3}{c|}{\textbf{Average results} $(n \ge 5)$} & \multicolumn{3}{c}{\textbf{Median results}}\\
				\cmidrule{4-15}
				& \textbf{searches} && $f_{eval.}$ & $time$ & $iter.$ & $f_{eval.}$ & $time$ & $iter.$ & $f_{eval.}$ & $time$ & $iter.$ & $f_{eval.}$ & $time$ & $iter.$ \\
				\midrule
				\direct & $-$ & $16/81$  & $452,244$	 & $263.29$	 & $1,645$ & $135,328$  & $220.05$   & $2,535$  & $638,666$   & $183.98$ & $1,121$   & $7,795$ 	 & $0.61$	  & $51$   \\
				\midrule
				\directrest		& $-$ & $23/81$ & $608,719$    & $5,518.88$ & $155$ & $268,291$ & $3,088.88$ & $84$ & $808,971$ & $6,948.29$ & $197$ & $12,861$ & $1.05$ & $48$	\\
				\midrule
				\directm & $-$ & $27/81$  & $745,865$	 & $340.83$	 & $2,505$ & $123,565$  & $384.15$   & $2,543$  & $1,111,924$ & $315.35$ & $2,482$ & $17,559$ 	 & $2.91$	  & $135$  \\
				\midrule
				\directz & $-$ & $24/81$ & $627,987$   & $4,486.48$ & $65,377$  & $134,552$ & $952.43$   & $22,558$ & $918,243$   & $6,565.33$ & $90,565$  & $7,185$     & $15.44$    & $437$    \\
				\midrule
				\directrev$^{*}$		& $3$ & $7/81$ & $223,483$    & $1,971.68$ & $27,093$ & $67,176$ & $437.96$ & $4,081$ & $315,429$ & $2,873.87$ & $40,629$ & $545$ & $0.13$ & $13$	\\
				\midrule
				\directav & $-$ & $41/81$ & $1,019,956$ & $830.72$   & $4,645$  & $203,905$   & $600.68$   & $3,699$  & $1,499,986$ & $966.03$   & $5,201$ & $2,000,000$ & $285.43$ & $176$  \\
				\midrule
				\dirmin$^{*}$		& $748$ & $5/81$ & $155,384$    & $18.45$ & $72$ & $67,178$ & $22.73$ & $91$ & $207,271$ & $15.94$ & $61$ & $361$ & $0.04$ & $1$ \\
				\midrule
				\plor & $-$ & $31/81$ & $775,748$   & $2,584.88$ & $57,574$  & $275,890$ & $2,162.11$ & $39,006$  & $1,069,782$   & $2,833.57$ & $68,496$ & $3,311$     & $1.12$     & $437$    \\
				\midrule
				\glbsolve & $-$ & $23/81$ & $624,411$ 	& $472.42$ 	& $2,326$ & $135,726$ & $335.40$  & $2,627$  & $911,873$ & $314.68$ & $2,149$ & $20,823$ 	  & $1.16$   & $54$ 	\\
				\midrule
				\symdirec & $-$ & $36/81$ & $923,592$   & $4,152.72$ & $10,128$  & $467,654$   & $1,288.24$ & $10,426$ & $1,191,790$ & $5,837.70$ & $9,953$  & $199,533$    & $446.61$   & $470$   \\
				\midrule
				\symdirect & $-$ & $35/81$ & $899,684$   & $5,283.04$ & $10,388$  & $603,737$   & $1,594.43$ & $10,931$ & $1,073,771$ & $7,452.82$  & $10,069$ & $124,961$   & $594.60$   & $355$ \\
				\midrule
				\directmr	& $-$ & $18/81$ & $502,032$    & $178.35$ & $2,689$ & $73,156$ & $9.05$ & $373$ & $754,313$ & $277.93$ & $4,051$ & $9,721$ & $0.52$ & $93$	\\
				\midrule
				\directmro$_{075}$		& $-$ & $16/81$ & $477,176$    & $241.40$ & $3,605$ & $69,755$ & $18.71$ & $520$ & $716,835$ & $372.40$ & $5,420$ & $8,547$ & $0.83$ & $103$	\\
				\midrule
				\birect & $-$ & $9/81$ & $255,671$   & $1,914.41$ & $6,829$  & $68,112$  & $914.36$   & $2,729$  & $366,000$   & $2,502.68$ & $9,241$  & $2,112$     & $1.28$   & $71$   \\
				\midrule
				\disimplcgb & $-$ & $46/81$ & $1,156,420$ & $3,903.78$ & $11,887$ & $224,138$   & $1,945.78$  & $16,156$ & $1,704,821$ & $5,055.54$ & $9,376$  & $2,000,000$ & $138.44$ & $24$ \\
				\midrule
				\disimplvgb & $-$ & $36/81$ & $898,336$   & $19,858.76$ & $2,109$ & $73,335$  & $1,623.49$  & $1,934$  & $1,383,631$ & $30,585.39$ & $2,211$ & $66,257$    & $4,298.97$  & $39$ \\
				\midrule
				\birectgb & $-$ & $13/81$ & $367,464$   & $2,451.11$ & $20,203$ & $70,721$    & $789.93$   & $7,937$   & $542,019$   & $3,428.28$ & $27,419$ & $5,782$     & $1.56$   & $153$   \\
				\midrule
				\birmin$^{*}$ & $1$ & $5/81$  & $125,541$ & $2,575.78$ & $13,106$ & $66,982$  & $1,433.43$ & $6,628$  & $159,987$ & $3,247.76$ & $16,917$ & $322$ & $0.07$ & $21$   \\
				\midrule
				\gbdirect & $-$ & $25/81$ & $671,030$   & $668.70$   & $8,092$  & $137,164$ & $1059.04$ & $15,501$ & $985,069$ & $439.08$ & $3,733$ & $22,541$ 	 & $2.40$ 	& $69$ \\
				\midrule
				\disimplc & $-$ & $46/81$ & $1,149,469$ & $4,257.54$ & $11,952$ & $215,702$   & $1,979.85$ & $13,649$ & $1,698,744$ & $5,593.43$ & $10,953$ & $2,000,000$ & $126.36$ & $22$ \\
				\midrule
				\disimplv & $-$ & $34/81$ & $844,074$   & $18,265.14$ & $709$   & $67,988$  & $1,436.35$  & $431$    & $1,300,595$ & $28,164.43$ & $872$   & $21,828$    & $667.34$    & $25$ \\
				\midrule
				\adc	 & $-$ 	& $30/81$ 	& $753,474$    & $16,537.93$  & $20,962$ & $74,665$ & $1,768.67$ & $8,999$ & $1,152,774$ & $25,225.73$ & $28,000$ & $8,868$ & $43.91$ & $603$	\\
				\midrule
				\directa & $-$ 	& $14/81$  	& $475,318$		& $20.71$	 & $95$    & $172,708$  & $11.37$    & $88$    & $653,324$   & $26.21$ 	  & $99$   & $65,253$ 	  & $2.34$	  & $44$   \\
				\midrule
				\directg & $-$ & $10/81$  & $310,535$	& $32.85$	 & $348$   & $84,905$   & $16.21$    & $309$   & $443,259$   & $42.64$    & $371$   & $9,835$ 	  & $0.45$	  & $51$   \\
				\midrule
				\directl & $-$ & $16/81$  & $430,885$	& $99.07$	 & $622$   & $68,918$   & $39.97$    & $362$   & $643,807$   & $133.84$   & $775$   & $9,601$ 	  & $0.52$	  & $46$   \\
				\midrule
				\directgl & $-$ & $3/81$  & $152,505$	& $11.78$	 & $99$   & $9,469$    & $0.68$     & $47$    & $236,645$   & $18.30$ 	  & $130$   & $7,737$ 	  & $0.33$	  & $36$   \\
				\midrule
				\glbolve & $-$ & $20/81$  & $534,930$	& $1,421.26$	 & $1,676$   & $201,103$    & $1,214.37$     & $2,576$    & $731,298$   & $1,542.98$ 	  & $1,147$   & $13,991$ 	  & $1.89$	  & $49$   \\
				\midrule
				\glccluster$^{*}$ & $21$ & $4/81$  & $116,281$	& $2,202.07$	 & $2$   & $68,612$    & $1,390.14$     & $2$    & $148,607$   & $2,760.65$ 	  & $2$   & $10,043$ 	  & $1.31$	  & $1$   \\

				\bottomrule
				\multicolumn{6}{l}{$^{*}$ -- a hybrid version of the algorithm, enriched with the local search subroutine}
			\end{tabular}
		}
		\label{tab:results}
	\end{table}

	\cref{tab:results} summarizes experimental results using $\varepsilon_{\rm pe} = 10^{-2}$.
	The smallest number of unsolved problems is achieved using \directgl{} $(3/81)$.
  At the same time, the second, third and fourth best algorithms are \glccluster{} $(4/81)$, \birmin{} $(5/81)$ and \dirmin{} $(5/81)$, hybrid versions enriched with the local search subroutines.
	In column `\textbf{Avg. \# local searches}' we report the average number of local searches performed by each hybridized algorithm.
	Hybridization of the \birmin{} algorithm allows solving more problems compared to, e.g., the globally biased version \birectgb{} $(12/81)$.
	Among traditional \direct-type algorithms, the second and third best algorithms are \birect{} and \directg{}.
	Both methods failed to solve $(9/81)$ and $(10/81)$ test problems accordingly.

	Furthermore, hybridization significantly reduces the total number of function evaluations (see \textbf{Average results} and \textbf{Median results} columns).
	The \birect{} and \plor{} were the most effective algorithms among the traditional algorithms based on the median number of function evaluations (see \textbf{Median results} column).
	However, for \plor, such performance needs to be interpreted correctly.
	As \plor{} restricts POH set to only two hyper-rectangles per iteration, a lower number of function evaluations are required to get closer to the solution for simpler (low-dimensional) problems.
	However, looking at the average number of function evaluations, even restricted to the simplest subset of problems (see \textbf{Average results $(n \le 4)$}), \plor{} performance is among the worst.
	\plor{} has failed on a larger number of simpler test problems than other approaches.
	In contrast, \directgl{} is only in eight place based on the median number of function evaluation criteria but is the only algorithm that solves all simpler $(n \le 4)$ problems and is the best performing algorithm, including hybridized versions.
	Moreover, on average, \directgl{} required approximately $35\%$ percent fewer evaluations of the objective function than the second-best, \birect{} algorithm, among all traditional \direct-type algorithms on a class of more challenging problems (see \textbf{Average results $(n \ge 5)$}).
	Not surprisingly, the hybridized \glccluster{} and \birmin{} algorithms within this class deliver the best average performance.
	Compared with the best traditional \direct-type algorithm, \directgl{}, \glccluster, and \birmin{} require approximately $37 \%$ and  $32 \%$ fewer evaluations.

	Based on the execution time ($time$), \directgl{} and \directa{} are the best among all traditional algorithms.
	\directa{} does not have a traditional POH selection procedure.
	Instead, the algorithm selects at least one candidate from each group of different measures.
	Therefore, the number of selected hyper-rectangles per iteration is larger, especially for higher dimensionality test problems.
	Consequently, the number of iterations ($iter.$) using \directa{} is among the smallest.
	Overall, \directgl{} showed the most promising performance among all tested traditional \direct-type algorithms.

	\begin{figure}
		\resizebox{\textwidth}{!}{
			\begin{tikzpicture}
				\begin{groupplot}[
					group style={
						group size=2 by 1,
						x descriptions at=edge bottom,
						y descriptions at=edge left,
						vertical sep=4pt,
						horizontal sep=7.5pt,
					},
					height=0.7\textwidth,width=0.5\textwidth
					]
					\nextgroupplot[
					xmode=log,
					title  = {Based on function evalutions},
					ymin=-0.01,ymax=1.01,
					ytick distance=0.1,
					xmin=9,xmax=2000000,
					xtick distance=10,
					xlabel = {Function evaluations},
					ylabel = {Proportion of problems solved},
					]

					\addplot[postaction={decoration={markings,mark=between positions 0 and 1 step 0.1 with {\node[diamond,draw=black,inner sep=1.25pt] {};}},decorate,},black,line width=0.5pt] table[x=T,y=DI] {figure_data/Profiler033.txt};
					\addplot[postaction={decoration={markings,mark=between positions 0 and 1 step 0.1 with {\node[star,star points=3,draw=armygreen,inner sep=1.25pt] {};}},decorate,},armygreen,line width=0.5pt] table[x=T,y=RS] {figure_data/Profiler033.txt};
					\addplot[postaction={decoration={markings,mark=between positions 0 and 1 step 0.1 with {\node[star,star points=5,draw=antiquebrass,inner sep=1.25pt] {};}},decorate,},antiquebrass,line width=0.5pt] table[x=T,y=MM] {figure_data/Profiler033.txt};
					\addplot[postaction={decoration={markings,mark=between positions 0 and 1 step 0.1 with {\node[star,star points=8,draw=blueryb,inner sep=1.25pt] {};}},decorate,},blueryb,line width=0.5pt] table[x=T,y=LL] {figure_data/Profiler033.txt};
					\addplot[postaction={decoration={markings,mark=between positions 0 and 1 step 0.1 with {\node[star,star points=11,draw=green,inner sep=1.25pt] {};}},decorate,},green,line width=0.5pt] table[x=T,y=RV] {figure_data/Profiler033.txt};
					\addplot[postaction={decoration={markings,mark=between positions 0 and 1 step 0.1 with {\node[regular polygon,regular polygon sides=3,draw=bluebell,inner sep=1.25pt] {};}},decorate,},bluebell,line width=0.5pt] table[x=T,y=AA] {figure_data/Profiler033.txt};
					\addplot[postaction={decoration={markings,mark=between positions 0 and 1 step 0.1 with {\node[regular polygon,regular polygon sides=4,draw=redwood,inner sep=1.25pt] {};}},decorate,},redwood,line width=0.5pt] table[x=T,y=DM] {figure_data/Profiler033.txt};
					\addplot[postaction={decoration={markings,mark=between positions 0 and 1 step 0.1 with {\node[regular polygon,regular polygon sides=5,draw=ao,inner sep=1.25pt] {};}},decorate,},ao,line width=0.5pt] table[x=T,y=PL] {figure_data/Profiler033.txt};
					\addplot[postaction={decoration={markings,mark=between positions 0 and 1 step 0.1 with {\node[regular polygon,regular polygon sides=6,draw=rosybrown,inner sep=1.25pt] {};}},decorate,},rosybrown,line width=0.5pt] table[x=T,y=GS] {figure_data/Profiler033.txt};
					\addplot[postaction={decoration={markings,mark=between positions 0 and 1 step 0.1 with {\node[regular polygon,regular polygon sides=7,draw=aureolin,inner sep=1.25pt] {};}},decorate,},aureolin,line width=0.5pt] table[x=T,y=SM] {figure_data/Profiler033.txt};
					\addplot[postaction={decoration={markings,mark=between positions 0 and 1 step 0.1 with {\node[rectangle,draw=DarkRed,inner sep=1.25pt] {};}},decorate,},DarkRed,line width=0.5pt] table[x=T,y=SN] {figure_data/Profiler033.txt};
					\addplot[postaction={decoration={markings,mark=between positions 0 and 1 step 0.1 with {\node[circle,draw=wheat,inner sep=1.25pt] {};}},decorate,},wheat,line width=0.5pt] table[x=T,y=MR] {figure_data/Profiler033.txt};
					\addplot[postaction={decoration={markings,mark=between positions 0 and 1 step 0.1 with {\node[diamond,draw=black,fill=black,inner sep=1.25pt] {};}},decorate,},black,line width=0.5pt] table[x=T,y=MO] {figure_data/Profiler033.txt};
					\addplot[postaction={decoration={markings,mark=between positions 0 and 1 step 0.1 with {\node[star,star points=3,draw=armygreen,inner sep=1.25pt] {};}},decorate,},armygreen,line width=0.5pt] table[x=T,y=BR] {figure_data/Profiler033.txt};
					\addplot[postaction={decoration={markings,mark=between positions 0 and 1 step 0.1 with {\node[star,star points=5,draw=antiquebrass,fill=antiquebrass,inner sep=1.25pt] {};}},decorate,},antiquebrass,line width=0.5pt] table[x=T,y=DG] {figure_data/Profiler033.txt};
					\addplot[postaction={decoration={markings,mark=between positions 0 and 1 step 0.1 with {\node[star,star points=8,draw=blueryb,fill=blueryb,inner sep=1.25pt] {};}},decorate,},blueryb,line width=0.5pt] table[x=T,y=VG] {figure_data/Profiler033.txt};
					\addplot[postaction={decoration={markings,mark=between positions 0 and 1 step 0.1 with {\node[star,star points=11,draw=green,fill=green,inner sep=1.25pt] {};}},decorate,},green,line width=0.5pt] table[x=T,y=GB] {figure_data/Profiler033.txt};
					\addplot[postaction={decoration={markings,mark=between positions 0 and 1 step 0.1 with {\node[regular polygon,regular polygon sides=3,draw=bluebell,fill=bluebell,inner sep=1.25pt] {};}},decorate,},bluebell,line width=0.5pt] table[x=T,y=BM] {figure_data/Profiler033.txt};
					\addplot[postaction={decoration={markings,mark=between positions 0 and 1 step 0.1 with {\node[regular polygon,regular polygon sides=4,draw=redwood,fill=redwood,inner sep=1.25pt,solid] {};}},decorate,},redwood,line width=0.5pt] table[x=T,y=GD] {figure_data/Profiler033.txt};
					\addplot[postaction={decoration={markings,mark=between positions 0 and 1 step 0.1 with {\node[regular polygon,regular polygon sides=5,draw=ao,fill=ao,inner sep=1.25pt,solid] {};}},decorate,},ao,line width=0.5pt] table[x=T,y=DS] {figure_data/Profiler033.txt};
					\addplot[postaction={decoration={markings,mark=between positions 0 and 1 step 0.1 with {\node[regular polygon,regular polygon sides=6,draw=rosybrown,fill=rosybrown,inner sep=1.25pt,solid] {};}},decorate,},rosybrown,line width=0.5pt] table[x=T,y=DV] {figure_data/Profiler033.txt};
					\addplot[postaction={decoration={markings,mark=between positions 0 and 1 step 0.1 with {\node[regular polygon,regular polygon sides=7,draw=aureolin,fill=aureolin,inner sep=1.25pt,solid] {};}},decorate,},aureolin,line width=0.5pt] table[x=T,y=AD] {figure_data/Profiler033.txt};
					\addplot[postaction={decoration={markings,mark=between positions 0 and 1 step 0.1 with {\node[rectangle,draw=DarkRed,fill=DarkRed,inner sep=1.25pt,solid] {};}},decorate,},DarkRed,line width=0.5pt] table[x=T,y=AG] {figure_data/Profiler033.txt};
					\addplot[postaction={decoration={markings,mark=between positions 0 and 1 step 0.1 with {\node[circle,draw=wheat,fill=wheat,inner sep=1.25pt] {};}},decorate,},wheat,line width=0.5pt] table[x=T,y=G] {figure_data/Profiler033.txt};
					\addplot[postaction={decoration={markings,mark=between positions 0 and 1 step 0.1 with {\node[cross out,draw=darkblue,inner sep=1.25pt] {};}},decorate,},darkblue,line width=0.5pt] table[x=T,y=L] {figure_data/Profiler033.txt};
					\addplot[black,line width=0.5pt] table[x=T,y=GL] {figure_data/Profiler033.txt};
					\addplot[black,dashed,line width=0.5pt] table[x=T,y=TS] {figure_data/Profiler033.txt};
					\addplot[black,densely dotted,line width=0.5pt] table[x=T,y=TC] {figure_data/Profiler033.txt};
					\nextgroupplot[
					legend pos=outer north east,
					xmode=log,
					legend style={font=\tiny},
					legend style={draw=none},
					legend style={row sep=0.5pt},
					title  = {Based on the execution time},
					ymin=-0.01,ymax=1.01,
					ytick distance=0.1,
					xmin=0.01,xmax=43000,
					xtick distance=10,
					xlabel = {Seconds},
					]

					\addlegendimage{legend image code/.code={\node [diamond,draw=black,inner sep=0.5pt] {};}}
					\addlegendentry{\direct}
					\addlegendimage{legend image code/.code={\node [star,star points=3,draw=armygreen,inner sep=0.5pt] {};}}
					\addlegendentry{\directrest}
					\addlegendimage{legend image code/.code={\node [star,star points=5,draw=antiquebrass,inner sep=0.5pt] {};}}
					\addlegendentry{\directm}
					\addlegendimage{legend image code/.code={\node [star,star points=8,draw=blueryb,inner sep=0.5pt] {};}}
					\addlegendentry{\directz}
					\addlegendimage{legend image code/.code={\node [star,star points=11,draw=green,inner sep=0.5pt] {};}}
					\addlegendentry{\directrev}
					\addlegendimage{legend image code/.code={\node [regular polygon,regular polygon sides=3,draw=bluebell,inner sep=0.5pt] {};}}
					\addlegendentry{\directav}
					\addlegendimage{legend image code/.code={\node [regular polygon,regular polygon sides=4,draw=redwood,inner sep=0.5pt] {};}}
					\addlegendentry{\dirmin}
					\addlegendimage{legend image code/.code={\node [regular polygon,regular polygon sides=5,draw=ao,inner sep=0.5pt] {};}}
					\addlegendentry{\plor}
					\addlegendimage{legend image code/.code={\node [regular polygon,regular polygon sides=6,draw=rosybrown,inner sep=0.5pt,solid] {};}}
					\addlegendentry{\glbsolve}
					\addlegendimage{legend image code/.code={\node [regular polygon,regular polygon sides=7,draw=aureolin,inner sep=0.5pt] {};}}
					\addlegendentry{\symdirec}
					\addlegendimage{legend image code/.code={\node [rectangle,draw=DarkRed,inner sep=0.5pt,solid] {};}}
					\addlegendentry{\symdirect}
					\addlegendimage{legend image code/.code={\node [circle,draw=wheat,inner sep=0.5pt] {};}}
					\addlegendentry{\directmr}
					\addlegendimage{legend image code/.code={\node [diamond,draw=black,fill=black,inner sep=0.5pt] {};}}
					\addlegendentry{\directmr$_{075}$}
					\addlegendimage{legend image code/.code={\node [star,star points=3,draw=armygreen,fill=armygreen,inner sep=0.5pt] {};}}
					\addlegendentry{\birect}
					\addlegendimage{legend image code/.code={\node [star,star points=5,draw=antiquebrass,fill=antiquebrass,inner sep=0.5pt] {};}}
					\addlegendentry{\disimplcgb}
					\addlegendimage{legend image code/.code={\node [star,star points=8,draw=blueryb,fill=blueryb,inner sep=0.5pt] {};}}
					\addlegendentry{\disimplvgb}
					\addlegendimage{legend image code/.code={\node [star,star points=11,draw=green,fill=green,inner sep=0.5pt] {};}}
					\addlegendentry{\birectgb}
					\addlegendimage{legend image code/.code={\node [regular polygon,regular polygon sides=3,draw=bluebell,fill=bluebell,inner sep=0.5pt] {};}}
					\addlegendentry{\birmin}
					\addlegendimage{legend image code/.code={\node [regular polygon,regular polygon sides=4,draw=redwood,fill=redwood,inner sep=0.5pt] {};}}
					\addlegendentry{\gbdirect}
					\addlegendimage{legend image code/.code={\node [regular polygon,regular polygon sides=5,draw=ao,fill=ao,inner sep=0.5pt] {};}}
					\addlegendentry{\disimplc}
					\addlegendimage{legend image code/.code={\node [regular polygon,regular polygon sides=6,draw=rosybrown,fill=rosybrown,inner sep=0.5pt,solid] {};}}
					\addlegendentry{\disimplv}
					\addlegendimage{legend image code/.code={\node [regular polygon,regular polygon sides=7,draw=aureolin,fill=aureolin,inner sep=0.5pt] {};}}
					\addlegendentry{\adc}
					\addlegendimage{legend image code/.code={\node [rectangle,draw=DarkRed,fill=DarkRed,inner sep=0.5pt,solid] {};}}
					\addlegendentry{\directa}
					\addlegendimage{legend image code/.code={\node [circle,draw=wheat,fill=wheat,inner sep=0.5pt] {};}}
					\addlegendentry{\directg}
					\addlegendimage{legend image code/.code={\node [cross out,draw=darkblue,inner sep=0.5pt] {};}}
					\addlegendentry{\directl}
					\addlegendimage{legend image code/.code={\node [draw=none,fill=none] {---};}}
					\addlegendentry{\directgl}
					\addlegendimage{legend image code/.code={\node [draw=none,fill=none] {- - -};}}
					\addlegendentry{\glbolve}
					\addlegendimage{legend image code/.code={\node [draw=none,fill=none] {.......};}}
					\addlegendentry{\glccluster}
					\addplot[postaction={decoration={markings,mark=between positions 0 and 1 step 0.1 with {\node[diamond,draw=black,inner sep=1.25pt] {};}},decorate,},black,line width=0.5pt] table[x=T,y=DI] {figure_data/Profiler022.txt};
					\addplot[postaction={decoration={markings,mark=between positions 0 and 1 step 0.1 with {\node[star,star points=3,draw=armygreen,inner sep=1.25pt] {};}},decorate,},armygreen,line width=0.5pt] table[x=T,y=RS] {figure_data/Profiler022.txt};
					\addplot[postaction={decoration={markings,mark=between positions 0 and 1 step 0.1 with {\node[star,star points=5,draw=antiquebrass,inner sep=1.25pt] {};}},decorate,},antiquebrass,line width=0.5pt] table[x=T,y=MM] {figure_data/Profiler022.txt};
					\addplot[postaction={decoration={markings,mark=between positions 0 and 1 step 0.1 with {\node[star,star points=8,draw=blueryb,inner sep=1.25pt] {};}},decorate,},blueryb,line width=0.5pt] table[x=T,y=LL] {figure_data/Profiler022.txt};
					\addplot[postaction={decoration={markings,mark=between positions 0 and 1 step 0.1 with {\node[star,star points=11,draw=green,inner sep=1.25pt] {};}},decorate,},green,line width=0.5pt] table[x=T,y=RV] {figure_data/Profiler022.txt};
					\addplot[postaction={decoration={markings,mark=between positions 0 and 1 step 0.1 with {\node[regular polygon,regular polygon sides=3,draw=bluebell,inner sep=1.25pt] {};}},decorate,},bluebell,line width=0.5pt] table[x=T,y=AA] {figure_data/Profiler022.txt};
					\addplot[postaction={decoration={markings,mark=between positions 0 and 1 step 0.1 with {\node[regular polygon,regular polygon sides=4,draw=redwood,inner sep=1.25pt] {};}},decorate,},redwood,line width=0.5pt] table[x=T,y=DM] {figure_data/Profiler022.txt};
					\addplot[postaction={decoration={markings,mark=between positions 0 and 1 step 0.1 with {\node[regular polygon,regular polygon sides=5,draw=ao,inner sep=1.25pt] {};}},decorate,},ao,line width=0.5pt] table[x=T,y=PL] {figure_data/Profiler022.txt};
					\addplot[postaction={decoration={markings,mark=between positions 0 and 1 step 0.1 with {\node[regular polygon,regular polygon sides=6,draw=rosybrown,inner sep=1.25pt] {};}},decorate,},rosybrown,line width=0.5pt] table[x=T,y=GS] {figure_data/Profiler022.txt};
					\addplot[postaction={decoration={markings,mark=between positions 0 and 1 step 0.1 with {\node[regular polygon,regular polygon sides=7,draw=aureolin,inner sep=1.25pt] {};}},decorate,},aureolin,line width=0.5pt] table[x=T,y=SM] {figure_data/Profiler022.txt};
					\addplot[postaction={decoration={markings,mark=between positions 0 and 1 step 0.1 with {\node[rectangle,draw=DarkRed,inner sep=1.25pt] {};}},decorate,},DarkRed,line width=0.5pt] table[x=T,y=SN] {figure_data/Profiler022.txt};
					\addplot[postaction={decoration={markings,mark=between positions 0 and 1 step 0.1 with {\node[circle,draw=wheat,inner sep=1.25pt] {};}},decorate,},wheat,line width=0.5pt] table[x=T,y=MR] {figure_data/Profiler022.txt};
					\addplot[postaction={decoration={markings,mark=between positions 0 and 1 step 0.1 with {\node[diamond,draw=black,fill=black,inner sep=1.25pt] {};}},decorate,},black,line width=0.5pt] table[x=T,y=MO] {figure_data/Profiler022.txt};
					\addplot[postaction={decoration={markings,mark=between positions 0 and 1 step 0.1 with {\node[star,star points=3,draw=armygreen,inner sep=1.25pt] {};}},decorate,},armygreen,line width=0.5pt] table[x=T,y=BR] {figure_data/Profiler022.txt};
					\addplot[postaction={decoration={markings,mark=between positions 0 and 1 step 0.1 with {\node[star,star points=5,draw=antiquebrass,fill=antiquebrass,inner sep=1.25pt] {};}},decorate,},antiquebrass,line width=0.5pt] table[x=T,y=DG] {figure_data/Profiler022.txt};
					\addplot[postaction={decoration={markings,mark=between positions 0 and 1 step 0.1 with {\node[star,star points=8,draw=blueryb,fill=blueryb,inner sep=1.25pt] {};}},decorate,},blueryb,line width=0.5pt] table[x=T,y=VG] {figure_data/Profiler022.txt};
					\addplot[postaction={decoration={markings,mark=between positions 0 and 1 step 0.1 with {\node[star,star points=11,draw=green,fill=green,inner sep=1.25pt] {};}},decorate,},green,line width=0.5pt] table[x=T,y=GB] {figure_data/Profiler022.txt};
					\addplot[postaction={decoration={markings,mark=between positions 0 and 1 step 0.1 with {\node[regular polygon,regular polygon sides=3,draw=bluebell,fill=bluebell,inner sep=1.25pt] {};}},decorate,},bluebell,line width=0.5pt] table[x=T,y=BM] {figure_data/Profiler022.txt};
					\addplot[postaction={decoration={markings,mark=between positions 0 and 1 step 0.1 with {\node[regular polygon,regular polygon sides=4,draw=redwood,fill=redwood,inner sep=1.25pt,solid] {};}},decorate,},redwood,line width=0.5pt] table[x=T,y=GD] {figure_data/Profiler022.txt};
					\addplot[postaction={decoration={markings,mark=between positions 0 and 1 step 0.1 with {\node[regular polygon,regular polygon sides=5,draw=ao,fill=ao,inner sep=1.25pt,solid] {};}},decorate,},ao,line width=0.5pt] table[x=T,y=DS] {figure_data/Profiler022.txt};
					\addplot[postaction={decoration={markings,mark=between positions 0 and 1 step 0.1 with {\node[regular polygon,regular polygon sides=6,draw=rosybrown,fill=rosybrown,inner sep=1.25pt,solid] {};}},decorate,},rosybrown,line width=0.5pt] table[x=T,y=DV] {figure_data/Profiler022.txt};
					\addplot[postaction={decoration={markings,mark=between positions 0 and 1 step 0.1 with {\node[regular polygon,regular polygon sides=7,draw=aureolin,fill=aureolin,inner sep=1.25pt,solid] {};}},decorate,},aureolin,line width=0.5pt] table[x=T,y=AD] {figure_data/Profiler022.txt};
					\addplot[postaction={decoration={markings,mark=between positions 0 and 1 step 0.1 with {\node[rectangle,draw=DarkRed,fill=DarkRed,inner sep=1.25pt,solid] {};}},decorate,},DarkRed,line width=0.5pt] table[x=T,y=AG] {figure_data/Profiler022.txt};
					\addplot[postaction={decoration={markings,mark=between positions 0 and 1 step 0.1 with {\node[circle,draw=wheat,fill=wheat,inner sep=1.25pt] {};}},decorate,},wheat,line width=0.5pt] table[x=T,y=G] {figure_data/Profiler022.txt};
					\addplot[postaction={decoration={markings,mark=between positions 0 and 1 step 0.1 with {\node[cross out,draw=darkblue,inner sep=1.25pt] {};}},decorate,},darkblue,line width=0.5pt] table[x=T,y=L] {figure_data/Profiler022.txt};
					\addplot[black,line width=0.5pt] table[x=T,y=GL] {figure_data/Profiler022.txt};
					\addplot[black,dashed,line width=0.5pt] table[x=T,y=TS] {figure_data/Profiler022.txt};
					\addplot[black,densely dotted,line width=0.5pt] table[x=T,y=TC] {figure_data/Profiler022.txt};

				\end{groupplot}		\end{tikzpicture}}
		\caption{Data profiles of \direct-type algorithms from \toolbox{} and \tomlab{} on the whole set of box-constrained optimization problems from \directlib using $\varepsilon_{\rm pe} = 10^{-2}$ in \eqref{eq:pe}.}
		\label{fig:box2}
	\end{figure}
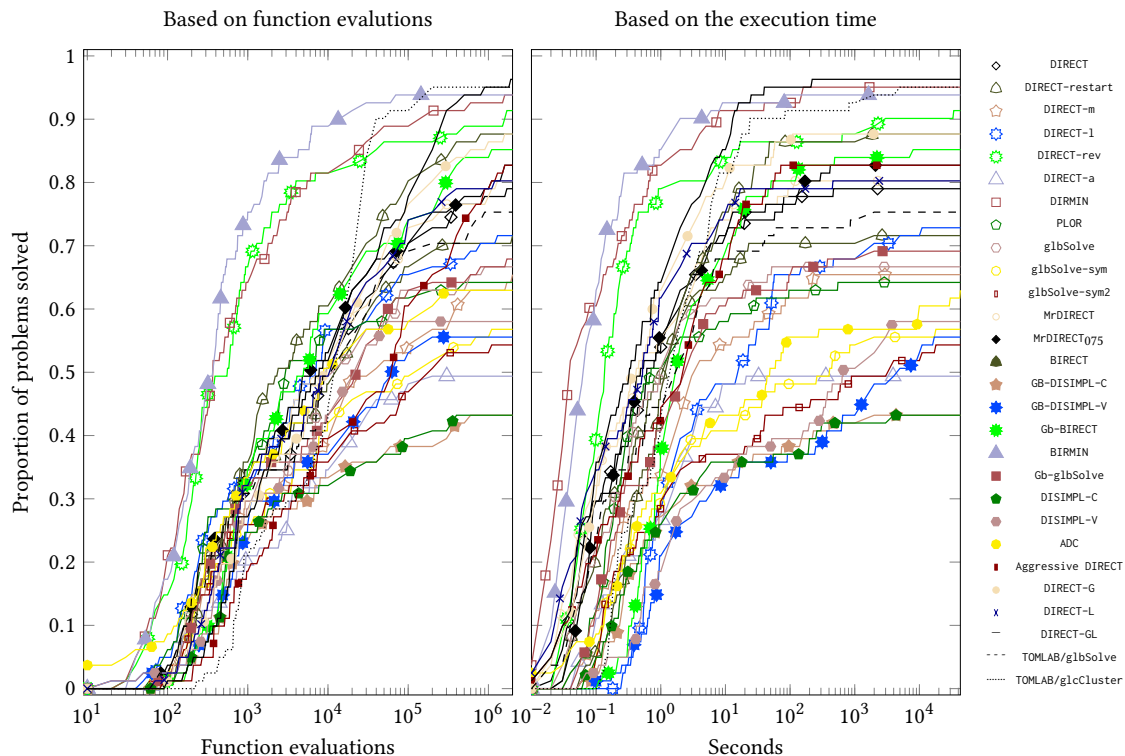

	The data profiles of all algorithms are shown in \cref{fig:box2}.
	Again, the data profiles confirm that the hybridized algorithms (enriched with \fmincon{} procedure) are more efficient than traditional \direct-type approaches.
	All three versions (\directrev, \dirmin, and \birmin) can solve about $60\%$ of problems in less than $1,000$ function evaluations (see the left panel of \cref{fig:box2}).
	However, by increasing the maximal budget of function evaluations, which is needed for more challenging problems, the \glccluster, and traditional \directgl{} algorithm outperformed all others, including hybridized ones.
	While \directgl{} delivers the best overall performance (based on function evaluations), data profiles in \cref{fig:box2} reveal that other \direct-type extensions (\plor, \directz, \birect, \birectgb) perform better when the maximal budget of function evaluations is $\leq 10^5$.
	Moreover, they can solve about $70\%$ of problems (mainly lower dimensionality) quicker.
	The two-step-based selection strategy in \directgl{} selects a more extensive set of POH.
	While for more straightforward problems, this is detrimental, it often helps to locate a global solution with higher accuracy faster.
	In terms of execution time (see the right panel in \cref{fig:box2}), \directgl{} is the fastest among the traditional \direct-type algorithms (excluding hybridized).

	Some algorithms have integrated schemes helping speed up the refinement of solutions.
  Therefore, we tested them with a much higher precision solution ($\varepsilon_{\rm pe} = 10^{-8}$).
	In \cref{tab:results2}, we summarize our experimental findings.
	\begin{table}
		\small
		\caption{The performance of selected \direct-type algorithms from \toolbox{} and \texttt{TOMLAB} based on the number of function evaluations ($f_{eval.}$),  the total execution time in seconds ($time$), and the total number of iterations ($iter.$) criteria on a set of box-constrained problems (from \directlib) using $\varepsilon_{\rm pe} = 10^{-8}$ in \eqref{eq:pe}}
		\resizebox{\textwidth}{!}{
			\begin{tabular}[tb]{@{\extracolsep{\fill}}l|r|r|rrr|rrr|rrr|rrr}
				\toprule
				\multirow{2}{*}{Algorithm} & \textbf{Avg. \# local} & \multirow{2}{*}{\textbf{Failed}}
				& \multicolumn{3}{c|}{\textbf{Average results}}  & \multicolumn{3}{c|}{\textbf{Average results} $(n \le 4)$} & \multicolumn{3}{c|}{\textbf{Average results} $(n \ge 5)$} & \multicolumn{3}{c}{\textbf{Median results}}\\
				\cmidrule{4-15}
				& \textbf{searches} && $f_{eval.}$ & $time$ & $iter.$ & $f_{eval.}$ & $time$     & $iter.$  & $f_{eval.}$ & $time$ & $iter.$ & $f_{eval.}$ & $time$ & $iter.$ \\
				\midrule
				\direct		& $-$   & $40/81$  & $1,066,171$ & $918.13$   & $7,823$  & $809,880$  & $1,483.54$ & $17,456$ & $1,216,930$ & $585.54$   & $2,157$  & $1,297,025$ & $268.16$ & $224$  	\\
				\midrule
				\directrest	& $-$ 	& $29/81$ & $778,771$    & $7,840.37$ & $200$ & $333,394$ & $4,377.27$ & $114$ & $1,035,618$ & $10,191.29$ & $252$ & $53,629$ & $6.00$ & $91$	\\
				\midrule
				\directrev$^{*}$  & $3$   & $15/81$ & $396,246$    & $1,095.13$ & $19,360$ & $200,403$  & $1,145.29$ & $9,595$ & $511,449$ 	& $1,065.63$   & $25,104$  & $844$	 	  & $0.16$   & $21$		\\
				\midrule
				\dirmin$^{*}$		& $1,001$ & $12/81$ & $345,763$    & $43.52$    & $164$ 	 & $135,323$  & $33.49$    & $106$    & $469,551$   & $49.41$    & $198$     & $501$ 	  & $0.09$ 	 & $2$		\\
				\midrule
				\glbsolve	& $-$   & $55/81$ & $1,423,737$  & $1,027.22$ & $9,609$  & $825,086$  & $1,720.33$ & $19,061$ & $1,775,885$ & $619.51$   & $4,050$  & $2,000,000$ & $340.86$ & $1,201$ 	\\
				\midrule
				\directmr	& $-$   & $38/81$ & $996,523$    & $2,021.44$ & $11,177$ & $706,737$  & $1,733.35$ & $21,678$ & $1,166,986$ & $2,190.90$   & $5,001$  & $542,125$   & $66.20$  & $406$	\\
				\midrule
				\directmro$_{075}$ & $-$ & $39/81$ & $1,036,760$ & $473.40$ & $6,438$  & $207,944$  & $104.02$    & $2,539$  & $1,524,299$ & $690.68$   & $8,732$  & $1,125,661$   & $141.77$ & $806$	\\
				\midrule
				\birmin$^{*}$		& $4$   & $12/81$ & $298,791$    & $4,827.15$ & $18,742$ & $136,309$  & $2,867.45$ & $11,448$ & $394,369$   & $5,979.91$ & $23,033$ & $379$       & $0.06$   & $21$   	\\
				\midrule
				\directl 	& $-$   & $21/81$  & $572,082$   & $151.41$	  & $952$    & $136,742$  & $86.78$    & $695$    & $828,164$   & $189.42$	 & $1,104$  & $23,069$	  & $1.45$   & $87$   	\\
				\midrule
				\directgl 	& $-$   & $6/81$  & $277,242$	 & $29.66$	   & $241$  & $85,123$   & $26.84$    & $252$  & $390,254$   & $31.32$	   & $235$   & $28,211$	  & $1.17$   & $60$   	\\
				\midrule
				\glccluster$^{*}$ & $76$  & $29/81$ & $753,378$	 & $16,307.54$ & $4$    & $338,957$  & $6,959.64$ & $3$    & $990,121$   & $22,553.22$ & $5$     & $66,391$	  & $54.86$  & $3$   	\\
				\bottomrule
				\multicolumn{6}{l}{$^{*}$ -- a hybrid version of the algorithm, enriched with the local search subroutine}
			\end{tabular}
		}
		\label{tab:results2}
	\end{table}
	First, the number of failed problems is much higher when higher accuracy is needed (see the \textbf{Failed} column in \cref{tab:results,tab:results2}).
	The smallest number of unsolved problems is achieved using \directgl{} $(6/81)$, where all failed test problems belong to the $(n \ge 5)$ class.
	The \directgl{} algorithm turns out to be more efficient even than hybrid methods.
	Overall, \directgl{} required approximately $7\%$ fewer function evaluations and took $32\%$ less time than the second and third best algorithms, \birmin{} and \dirmin, accordingly (see \textbf{Average results} column in \cref{tab:results2}).
	However, the \birmin{} algorithm has the best median value (see \textbf{Median results} column), solving at least half of the problems with the best performance.
	Finally, let us stress the inefficiency of the original \direct{} algorithm (\direct{} and \glbsolve{} implementations).
	As the median value is more than $2,000,000$, \glbsolve{} failed more than half of the test problems to solve.
	The performance of the \glccluster{} algorithm, which showed the best average results in the previous study, has also decreased significantly (see the \textbf{Average results} column in \cref{tab:results,tab:results2}).
	It turns out that the use of \fmincon{} with default parameters in hybridized methods (\directrev, \dirmin, and \birmin) is much more effective in finding a solution with higher accuracy than the \glccluster.
	Furthermore, \glccluster{} turned out to be the slowest algorithm among all involved in this study.

	\begin{figure}
		\resizebox{\textwidth}{!}{
			\begin{tikzpicture}
				\begin{groupplot}[
					group style={
						group size=2 by 1,
						x descriptions at=edge bottom,
						y descriptions at=edge left,
						vertical sep=4pt,
						horizontal sep=7.5pt,
					},
					height=0.4\textwidth,width=0.5\textwidth,
					xmode=log
					]
					\nextgroupplot[
					xmode=log,
					title  = {Based on function evalutions},
					ymin=-0.01,ymax=1.01,
					ytick distance=0.1,
					xmin=9,xmax=2000000,
					xtick distance=10,
					xlabel = {Function evaluations},
					ylabel = {Proportion of problems solved},
					]
					\addplot[postaction={decoration={markings,mark=between positions 0 and 1 step 0.1 with {\node[diamond,draw=black,inner sep=1.25pt] {};}},decorate,},black,line width=0.5pt] table[x=T,y=DI] {figure_data/Profiler044.txt};
					\addplot[postaction={decoration={markings,mark=between positions 0 and 1 step 0.1 with {\node[star,star points=3,draw=armygreen,inner sep=1.25pt] {};}},decorate,},armygreen,line width=0.5pt] table[x=T,y=RS] {figure_data/Profiler044.txt};
					\addplot[postaction={decoration={markings,mark=between positions 0 and 1 step 0.1 with {\node[star,star points=11,draw=green,inner sep=1.25pt] {};}},decorate,},green,line width=0.5pt] table[x=T,y=RV] {figure_data/Profiler044.txt};
					\addplot[postaction={decoration={markings,mark=between positions 0 and 1 step 0.1 with {\node[regular polygon,regular polygon sides=4,draw=redwood,inner sep=1.25pt] {};}},decorate,},redwood,line width=0.5pt] table[x=T,y=DM] {figure_data/Profiler044.txt};
					\addplot[postaction={decoration={markings,mark=between positions 0 and 1 step 0.1 with {\node[regular polygon,regular polygon sides=6,draw=rosybrown,inner sep=1.25pt] {};}},decorate,},rosybrown,line width=0.5pt] table[x=T,y=GS] {figure_data/Profiler044.txt};
					\addplot[postaction={decoration={markings,mark=between positions 0 and 1 step 0.1 with {\node[circle,draw=wheat,inner sep=1.25pt] {};}},decorate,},wheat,line width=0.5pt] table[x=T,y=MR] {figure_data/Profiler044.txt};
					\addplot[postaction={decoration={markings,mark=between positions 0 and 1 step 0.1 with {\node[diamond,draw=black,fill=black,inner sep=1.25pt] {};}},decorate,},black,line width=0.5pt] table[x=T,y=MO] {figure_data/Profiler044.txt};
					\addplot[postaction={decoration={markings,mark=between positions 0 and 1 step 0.1 with {\node[regular polygon,regular polygon sides=3,draw=bluebell,fill=bluebell,inner sep=1.25pt] {};}},decorate,},bluebell,line width=0.5pt] table[x=T,y=BM] {figure_data/Profiler044.txt};
					\addplot[postaction={decoration={markings,mark=between positions 0 and 1 step 0.1 with {\node[cross out,draw=darkblue,inner sep=1.25pt] {};}},decorate,},darkblue,line width=0.5pt] table[x=T,y=L] {figure_data/Profiler044.txt};
					\addplot[black,line width=0.5pt] table[x=T,y=GL] {figure_data/Profiler044.txt};
					\addplot[black,densely dotted,line width=0.5pt] table[x=T,y=TC] {figure_data/Profiler044.txt};
					\nextgroupplot[
					legend pos=outer north east,
					xmode=log,
					legend style={font=\tiny},
					legend style={draw=none},
					legend style={row sep=0.75pt},
					title  = {Based on the execution time},
					ymin=-0.01,ymax=1.01,
					ytick distance=0.1,
					xmin=0.01,xmax=43000,
					xtick distance=10,
					xlabel = {Seconds},
					]

					\addlegendimage{legend image code/.code={\node [diamond,draw=black,inner sep=0.5pt] {};}}
					\addlegendentry{\direct}
					\addlegendimage{legend image code/.code={\node [star,star points=3,draw=armygreen,inner sep=0.5pt] {};}}
					\addlegendentry{\directrest}
					\addlegendimage{legend image code/.code={\node [star,star points=11,draw=green,inner sep=0.5pt] {};}}
					\addlegendentry{\directrev}
					\addlegendimage{legend image code/.code={\node [regular polygon,regular polygon sides=4,draw=redwood,inner sep=0.5pt] {};}}
					\addlegendentry{\dirmin}
					\addlegendimage{legend image code/.code={\node [regular polygon,regular polygon sides=6,draw=rosybrown,inner sep=0.5pt,solid] {};}}
					\addlegendentry{\glbsolve}
					\addlegendimage{legend image code/.code={\node [circle,draw=wheat,inner sep=0.5pt] {};}}
					\addlegendentry{\directmr}
					\addlegendimage{legend image code/.code={\node [diamond,draw=black,fill=black,inner sep=0.5pt] {};}}
					\addlegendentry{\directmr$_{075}$}
					\addlegendimage{legend image code/.code={\node [regular polygon,regular polygon sides=3,draw=bluebell,fill=bluebell,inner sep=0.5pt] {};}}
					\addlegendentry{\birmin}
					\addlegendimage{legend image code/.code={\node [cross out,draw=darkblue,inner sep=0.5pt] {};}}
					\addlegendentry{\directl}
					\addlegendimage{legend image code/.code={\node [draw=none,fill=none] {---};}}
					\addlegendentry{\directgl}
					\addlegendimage{legend image code/.code={\node [draw=none,fill=none] {.......};}}
					\addlegendentry{\glccluster}
					\addplot[postaction={decoration={markings,mark=between positions 0 and 1 step 0.1 with {\node[diamond,draw=black,inner sep=1.25pt] {};}},decorate,},black,line width=0.5pt] table[x=T,y=DI] {figure_data/Profiler055.txt};
					\addplot[postaction={decoration={markings,mark=between positions 0 and 1 step 0.1 with {\node[star,star points=3,draw=armygreen,inner sep=1.25pt] {};}},decorate,},armygreen,line width=0.5pt] table[x=T,y=RS] {figure_data/Profiler055.txt};
					\addplot[postaction={decoration={markings,mark=between positions 0 and 1 step 0.1 with {\node[star,star points=11,draw=green,inner sep=1.25pt] {};}},decorate,},green,line width=0.5pt] table[x=T,y=RV] {figure_data/Profiler055.txt};
					\addplot[postaction={decoration={markings,mark=between positions 0 and 1 step 0.1 with {\node[regular polygon,regular polygon sides=4,draw=redwood,inner sep=1.25pt] {};}},decorate,},redwood,line width=0.5pt] table[x=T,y=DM] {figure_data/Profiler055.txt};
					\addplot[postaction={decoration={markings,mark=between positions 0 and 1 step 0.1 with {\node[regular polygon,regular polygon sides=6,draw=rosybrown,inner sep=1.25pt] {};}},decorate,},rosybrown,line width=0.5pt] table[x=T,y=GS] {figure_data/Profiler055.txt};
					\addplot[postaction={decoration={markings,mark=between positions 0 and 1 step 0.1 with {\node[circle,draw=wheat,inner sep=1.25pt] {};}},decorate,},wheat,line width=0.5pt] table[x=T,y=MR] {figure_data/Profiler055.txt};
					\addplot[postaction={decoration={markings,mark=between positions 0 and 1 step 0.1 with {\node[diamond,draw=black,fill=black,inner sep=1.25pt] {};}},decorate,},black,line width=0.5pt] table[x=T,y=MO] {figure_data/Profiler055.txt};
					\addplot[postaction={decoration={markings,mark=between positions 0 and 1 step 0.1 with {\node[regular polygon,regular polygon sides=3,draw=bluebell,fill=bluebell,inner sep=1.25pt] {};}},decorate,},bluebell,line width=0.5pt] table[x=T,y=BM] {figure_data/Profiler055.txt};
					\addplot[postaction={decoration={markings,mark=between positions 0 and 1 step 0.1 with {\node[cross out,draw=darkblue,inner sep=1.25pt] {};}},decorate,},darkblue,line width=0.5pt] table[x=T,y=L] {figure_data/Profiler055.txt};
					\addplot[black,line width=0.5pt] table[x=T,y=GL] {figure_data/Profiler055.txt};
					\addplot[black,densely dotted,line width=0.5pt] table[x=T,y=TC] {figure_data/Profiler055.txt};

				\end{groupplot}
		\end{tikzpicture}}
		\caption{Data profiles of selected \direct-type algorithms from \toolbox{} and \tomlab{} on the whole set of box-constrained optimization problems from \directlib using $\varepsilon_{\rm pe} = 10^{-8}$ in \eqref{eq:pe}.}
		\label{fig:box3}
	\end{figure}
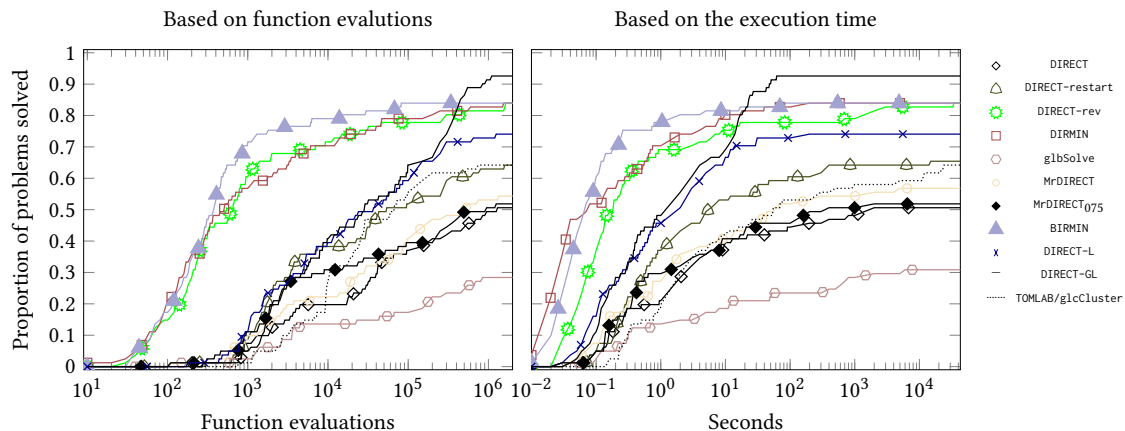
	Finally, the data profiles of the selected algorithms are shown in \cref{fig:box3}.
	Once again, the date profiles confirm that the hybridized algorithms are more efficient than traditional \direct-type approaches.
	However, by increasing the maximal budget of function evaluations, the traditional \directgl{} algorithm starts outperforming all algorithms, including hybridized ones.

	\subsection{Comparison of \direct-type algorithms for constrained global optimization}
	\label{sec:general_constraints}

	The comparison presented in this section was carried out using $80$ global optimization test problems with various constraints.
	In \directlib, $35$ test problems contain linear constraints, $39$ problems have non-linear constraints where $5$ include equality constraints.
	All necessary details about the test problems are given in \cref{sec:directlib},~\cref{tab:resultslow}.
	We used the same stopping condition in these experimental investigations as in the previous ones, and the value $\varepsilon_{\rm pe} = 10^{-2}$.

	Let us stress that $5$ of the test problems contain equality constraints, which we transform into inequality constraints as follows:
	\begin{equation}
		\label{eq:equolity-inequality}
		\mathbf{h} (\mathbf{x}) = 0 \rightarrow
		|\mathbf{h} (\mathbf{x}) - \varepsilon_{\rm h}| \leq 0,
	\end{equation}
	where $\varepsilon_{\rm h} > 0$ is a small tolerance for equality constraints.
	In our experiments, it was set to $10^{-8}$.

	\subsubsection{Test results on problems with hidden constraints}

	In the first part, we compared \direct-type versions devoted to problems with hidden constraints.
	We have used all constrained test problems but assumed that any information about the constraints is unavailable.
	In the experimental investigation, the hidden search area $(D^{\rm hidden})$ was defined as
	\begin{equation}
		\label{eq:dhidden}
		D^{\rm hidden} = \{\mathbf{x} \in D : \mathbf{g}(\mathbf{x}) \leq 0, \mathbf{h}(\mathbf{x}) = 0 \}
	\end{equation}
	Still, this information is unavailable for the tested algorithms and used only to determine whether a certain point is feasible or not.
	\begin{table}
		\small
		\caption{The performance of \direct-type algorithms from \toolbox{} and \texttt{TOMLAB} based on the number of function evaluations ($f_{eval.}$), the total execution time in seconds ($time$), and the total number of iterations ($iter.$) criteria on a set of constrained (hidden, general, and linear) optimization problems}
		\resizebox{\textwidth}{!}{
			\begin{tabular}[tb]{@{\extracolsep{\fill}}l|l|r|rrr|rrr|rrr|rrr}
				\toprule
				\multirow{2}{*}{Algorithm}  & \multirow{2}{*}{\textbf{Parameter}} & \multirow{2}{*}{\textbf{Failed}}
				& \multicolumn{3}{c}{\textbf{Average results}}  & \multicolumn{3}{c}{\textbf{Average results (Non-lin. constr.)}} & \multicolumn{3}{c}{\textbf{Average results (Lin. constr.)}} & \multicolumn{3}{c}{\textbf{Median results}}\\
				\cmidrule{4-15}
				&&& $f_{eval.}$ & $time$ & $iter.$ & $f_{eval.}$ & $time$ & $iter.$ & $f_{eval.}$ & $time$ & $iter.$ & $f_{eval.}$ & $time$ & $iter.$ \\
				\midrule
				\multicolumn{14}{c}{Performance of \direct-type algorithms for problems with hidden constraints}\\
				\midrule
				\directnas     & --  & $29/80$  & $711,720$   & $15,372.02$ & $20,049$ & $943,276$   & $20,453.13$ & $31,129$ & $414,006$   & $8,839.17$ & $5,802$  & $9,124$     & $18.68$	  & $144$  \\
				\directbarrier & --  & $46/80$  & $1,192,521$ & $2,908.88$  & $38,069$ & $1,260,510$ & $3,821.84$  & $51,635$ & $1,105,106$ & $1,735.09$ & $20,627$ & $2,000,000$ & $723.03$  & $2596$ \\
				\directsub     & $sub = 2$ & $53/80$  & $1,364,353$ & $860.50$ 	& $6,896$  & $1,316,302$ & $886.75$    & $6,042$  & $1,426,132$ & $826.76$   & $7,994$  & $2,000,000$ & $295.60$  & $65$   \\
				\directsub     & $sub = 3$ & $48/80$  & $1,240,771$ & $773.74$    & $3,623$  & $1,270,693$ & $773.81$    & $902$    & $1,202,301$ & $773.65$   & $7,121$  & $2,000,000$ & $286.48$  & $251$  \\
				\directsub 	   & $sub = 5$ & $47/80$  & $1,229,037$ & $1,359.66$  & $11,199$ & $1,367,032$ & $1,526.87$  & $10,660$ & $1,051,614$ & $1,144.67$ & $11,892$ & $2,000,000$ & $559.34$  & $1,294$\\
				\directh       & --  & $18/80$  & $470,807$   & $312.40$ 	& $104$    & $648,370$   & $491.81$    & $92$     & $242,513$   & $81.74$   & $119$    & $7,068$     & $0.57$	  & $32$   \\
				\midrule
				\multicolumn{14}{c}{Performance of \direct-type algorithms for generally constrained optimization problems}\\
				\midrule
				\directc       & --  & $11/80$  & $330,659$   & $67.46$     & $352$    & $480,785$   & $106.79$    & $542$    & $137,639$   & $16.90$    & $108$    & $3,759$ 	  & $0.34$	  & $37$   \\
				\directce      & --  & $7/80$   & $258,462$   & $40.02$     & $270$    & $380,598$   & $62.13$     & $368$    & $101,430$   & $11.59$    & $143$    & $9,768$     & $0.86$	  & $75$   \\
				\directcemin$^{*}$   & --  & $2/80$   & $62,233$ 	  & $10.72$ 		& $45$ 	   & $109,975$ 	 & $18.99$ 	   & $77$     & $852$     & $0.09$ 	 & $4$ 	& $124$ 	  & $0.04$	  & $1$    \\

				\directll      & $\gamma = 10$   & $43/80$  & $1,087,528$ & $228.17$    & $1,688$  & $1,206,763$ & $363.14$    & $2,617$  & $934,227$   & $54.63$    & $494$    & $2,000,000$ & $0.19$	  & $49$   \\
				\directll     & $\gamma = 10^2$ & $41/80$  & $1,051,478$ & $870.98$ & $3,369$ & $1,181,710$ & $1,425.34$ & $3,974$ & $884,038$ & $158.24$ & $2,591$ & $2,000,000$ & $2.49$	& $64$   \\
				\directll     & $\gamma = 10^3$ & $40/80$  & $1,042,671$ & $1,144.35$ & $8,203$ & $1,151,444$ & $1,249.36$ & $5,178$ & $902,820$ & $1,009.33$ & $12,093$ & $1,564,860$ & $62.27$	& $241$   \\
				\glsolve     & -- & $24/80$  & $607,397$ & $11,031.02$ & $13,568$ & $801,828$ & $14,550.11$ & $20,288$ & $357,415$ & $6,506.48$ & $4,927$ & $3,013$ & $2.75$ & $145$    \\
				\glccluster$^{*}$     & -- & $8/80$ & $207,226$ & $3,780.92$ & $2$ & $364,484$ & $6,718.99$ & $2$ & $5,038$ & $3.40$ & $1$ & $2,734$ & $1.40$ & $1$    \\
				\midrule
				\multicolumn{14}{c}{Performance of \direct-type algorithms devoted for problems with linear constraints only}\\
				\midrule
				\disimpllc     & -- & $5/35$  & N/A & N/A & N/A & N/A & N/A & N/A & $290,402$ & $6,424.68$ & $387$ & $443$ & $0.12$	& $27$   \\
				\disimpllv     & -- & $3/35$  & N/A & N/A & N/A & N/A & N/A & N/A & $171,738$ & $3,686.81$ & $24$  & $16$  & $0.01$	& $1$   \\
				\bottomrule
				\multicolumn{15}{l}{$^{*}$ -- a hybrid version of the algorithm, enriched with the local search subroutine} \\
				\multicolumn{15}{l}{N/A -- not available}
			\end{tabular}
		}
		\label{tab:results1}
	\end{table}
	Obtained experimental results are summarized in the upper part of~\cref{tab:results1} (see `Performance of \direct-type algorithms for problems with hidden constraints').
	First, let us note that $13$ out of the $80$ test problems contain complex constraints leading to a tiny feasible region.
	As the algorithms within this class do not use any information about constraint functions, none of the tested algorithms could find a single feasible point for these 13 test problems.

	The best among all \direct-type algorithms for problems with hidden constraints is \directh{} (failed to solve $(18/80)$), while the second-best is \directnas{} (failed to solve $(29/80)$).
	The median number of function evaluations (see $f_{eval}$ in \textbf{Median results}) is similar for both.
	Still, \directh{} is the best, primarily based on the number of iterations ($iter.$) and the execution $time$: \directh{} took around $4.5$ times fewer iterations and approximately $33$ times less execution time than the second-best \directnas, algorithm.
	The speed is the essential factor differentiating \directh{} from \directnas.

	For an extra subdividing step-based  \directsub, the user must define how often this step is activated.
	Unfortunately, the authors in \cite{Na2017} did not make any sensitivity analysis and guidance.
	In our experiments, we start the subdividing step at $sub^k, k = 1,2,\dots$ iterations.
	We tested three different values for the variable $sub$, i.e., $sub = 2, 3$, and $5$.
	Our experience showed that an extra subdividing step combined with a traditional barrier approach based on \directsub{} did not significantly improve performance (based on the number of \textbf{Failed} problems and the \textbf{Average results}) over the original \directbarrier.
	The most obvious difference is that, on average, \directsub{} subdivides much more POH per iteration because of an extra subdividing step, leading to a smaller number of iterations ($iter.$) and the execution $time$.
	\directsub{} algorithm suffers solving larger dimensionality and problems where $D^{\rm hidden}$ contains non-linear constraints (see \textbf{Average results (Non-lin. constr.)} column).
	Therefore, this limits \directsub{} applicability primarily to low-dimensional problems.

	Finally, in \cref{fig:Nonlinear}, the comparative performance of algorithms using the data profiles is demonstrated.
	They confirm that the \directh{} is the most effective optimizer in this class and has the highest efficiency based on the function evaluations and execution time.

	\subsubsection{Test results on problems with general constraints}

	Better performance of \direct-type algorithms can be expected when the constraint function information is known.
	Let us note that all \direct-type algorithms considered in the previous section are included in this analysis.
	The new experimental results are added in the middle part of~\cref{tab:results1} (see `Performance of \direct-type algorithms for generally constrained optimization problems).
	First, let us note that the recently proposed \directc{} and \directce{} algorithms can be much more successful (compared to the algorithms for problems with hidden) in solving problems containing complex and tiny feasible regions.
	The best average results among traditional \direct-type algorithms for constrained optimization were obtained using the \directce{} (failed to solve $7/80$).
	Interestingly, the performance based on the number of failed problems using \directc{} is worse, but based on the median results, it outperforms \directce{} quite clearly.
	It looks that \directc{} is effective on simpler problems, but the effectiveness drops in solving more complicated problems, e.g., higher dimensionality with non-linear constraints.
  We also note that the solution point is often located on the feasible region’s boundaries for optimization problems with general constraints.
  The common problem of some \direct-type algorithms in this class (\directbarrier, \directsub) is that hyper-rectangles with infeasible midpoints situated closely to the edges of feasibility are penalized with large values, resulting in a low probability of being elected as POH.
  In such situations, these algorithms converge very slowly.

	The best traditional \direct-type algorithm based on the median value was the \glsolve{} method (see \textbf{Median results} column in \cref{tab:results1}).
	However, this is the only category where this algorithm showed the best results.
	Overall, the \glsolve{} algorithm required $57 \%$ times more objective function evaluations than \directce.
	Furthermore, the \glsolve{} appears to be the slowest algorithm in this class.

	Hybridized \directcemin{} and \glccluster{} are the only two candidates among all approaches within this class.
	Again, incorporating the local minimization procedure into \directce{} improves the performance significantly, e.g., it reduces the overall number of function evaluations approximately seven times.
	Moreover, \directcemin{} fails to solve only two test problems.
  \glccluster{} fails to solve eight test instances and, on average, is around 353 times slower than the \directcemin{} method.

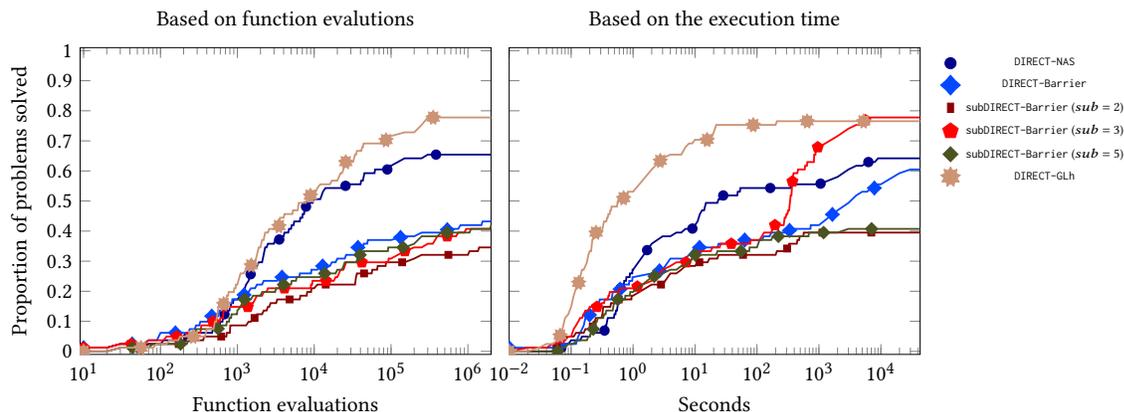
\begin{figure}
		\resizebox{\textwidth}{!}{
			\begin{tikzpicture}
				\begin{groupplot}[
					group style={
						group size=2 by 1,
						x descriptions at=edge bottom,
						y descriptions at=edge left,
						horizontal sep=7.5pt,
					},
					height=0.4\textwidth,width=0.5\textwidth,
					xmode=log
					]
					\nextgroupplot[
					xmode=log,
					title  = {Based on function evalutions},
					ymin=-0.01,ymax=1.01,
					ytick distance=0.1,
					xmin=9,xmax=2000000,
					xtick distance=10,
					xlabel = {Function evaluations},
					ylabel = {Proportion of problems solved},
					]
					\addplot[postaction={decoration={markings,mark=between positions 0 and 1 step 0.1 with {\node[circle,draw=darkblue,fill=darkblue,inner sep=1.25pt] {};}},decorate,},darkblue,line width=0.75pt] table[x=T,y=NAS] {figure_data/General.txt};
					\addplot[postaction={decoration={markings,mark=between positions 0 and 1 step 0.1 with {\node[star,star points=4,draw=blueryb,fill=blueryb,inner sep=1.25pt] {};}},decorate,},blueryb,line width=0.75pt] table[x=T,y=BAR] {figure_data/General.txt};
					\addplot[postaction={decoration={markings,mark=between positions 0 and 1 step 0.1 with {\node[draw=DarkRed,fill=DarkRed,inner sep=1.25pt] {};}},decorate,},DarkRed,line width=0.75pt] table[x=T,y=SU2] {figure_data/General.txt};
					\addplot[postaction={decoration={markings,mark=between positions 0 and 1 step 0.1 with {\node[regular polygon,regular polygon sides=5,draw=red,fill=red,inner sep=1.25pt] {};}},decorate,},red,line width=0.75pt] table[x=T,y=SU3] {figure_data/General.txt};
					\addplot[postaction={decoration={markings,mark=between positions 0 and 1 step 0.1 with {\node[diamond,draw=armygreen,fill=armygreen,inner sep=1.25pt] {};}},decorate,},armygreen,line width=0.75pt] table[x=T,y=SU5] {figure_data/General.txt};
					\addplot[postaction={decoration={markings,mark=between positions 0 and 1 step 0.1 with {\node[star,star points=9,fill=antiquebrass,draw=antiquebrass,inner sep=1.25pt] {};}},decorate,},antiquebrass,line width=0.75pt] table[x=T,y=H] {figure_data/General.txt};
					\nextgroupplot[
					legend pos=outer north east,
					xmode=log,
					legend style={font=\tiny},
					legend style={draw=none},
					legend style={row sep=0.5pt},
					title  = {Based on the execution time},
					ymin=-0.01,ymax=1.01,
					ytick distance=0.1,
					xmin=0.01,xmax=43000,
					xtick distance=10,
					xlabel = {Seconds},
					]
					\addlegendimage{legend image code/.code={\node [circle,draw=darkblue,fill=darkblue,inner sep=1.25pt] {};}}
					\addlegendentry{\directnas}
					\addlegendimage{legend image code/.code={\node [star,star points=4,draw=blueryb,fill=blueryb,inner sep=1.25pt] {};}}
					\addlegendentry{\directbarrier}
					\addlegendimage{legend image code/.code={\node [draw=DarkRed,fill=DarkRed,inner sep=1.25pt] {};}}
					\addlegendentry{\directsub{} ($sub=2$)}
					\addlegendimage{legend image code/.code={\node [regular polygon,regular polygon sides=5,draw=red,fill=red,inner sep=1.25pt] {};}}
					\addlegendentry{\directsub{} ($sub=3$)}
					\addlegendimage{legend image code/.code={\node [diamond,draw=armygreen,fill=armygreen,inner sep=1.25pt] {};}}
					\addlegendentry{\directsub{} ($sub=5$)}
					\addlegendimage{legend image code/.code={\node [star,star points=9,fill=antiquebrass,draw=antiquebrass,inner sep=1.25pt] {};}}
					\addlegendentry{\directh}
					\addplot[postaction={decoration={markings,mark=between positions 0 and 1 step 0.1 with {\node[circle,draw=darkblue,fill=darkblue,inner sep=1.25pt] {};}},decorate,},darkblue,line width=0.75pt] table[x=T,y=NAS] {figure_data/Hidden.txt};
					\addplot[postaction={decoration={markings,mark=between positions 0 and 1 step 0.1 with {\node[star,star points=4,draw=blueryb,fill=blueryb,inner sep=1.25pt] {};}},decorate,},blueryb,line width=0.75pt] table[x=T,y=BAR] {figure_data/Hidden.txt};
					\addplot[postaction={decoration={markings,mark=between positions 0 and 1 step 0.1 with {\node[draw=DarkRed,fill=DarkRed,inner sep=1.25pt] {};}},decorate,},DarkRed,line width=0.75pt] table[x=T,y=SU2] {figure_data/Hidden.txt};
					\addplot[postaction={decoration={markings,mark=between positions 0 and 1 step 0.1 with {\node[regular polygon,regular polygon sides=5,draw=red,fill=red,inner sep=1.25pt] {};}},decorate,},red,line width=0.75pt] table[x=T,y=SU3] {figure_data/Hidden.txt};
					\addplot[postaction={decoration={markings,mark=between positions 0 and 1 step 0.1 with {\node[diamond,draw=armygreen,fill=armygreen,inner sep=1.25pt] {};}},decorate,},armygreen,line width=0.75pt] table[x=T,y=SU5] {figure_data/Hidden.txt};
					\addplot[postaction={decoration={markings,mark=between positions 0 and 1 step 0.1 with {\node[star,star points=9,fill=antiquebrass,draw=antiquebrass,inner sep=1.25pt] {};}},decorate,},antiquebrass,line width=0.75pt] table[x=T,y=H] {figure_data/Hidden.txt};
				\end{groupplot}
		\end{tikzpicture}}
		\caption{Data profiles of \direct-type algorithms for problems with hidden constraints on the whole set of constrained optimization test problems. Explicit information about the constraints for all these algorithms was unknown.}
		\label{fig:Nonlinear}
	\end{figure}

\begin{figure}
		\resizebox{\textwidth}{!}{
			\begin{tikzpicture}
				\begin{groupplot}[
					group style={
						group size=2 by 1,
						x descriptions at=edge bottom,
						y descriptions at=edge left,
						horizontal sep=7.5pt,
					},
					height=0.4\textwidth,width=0.5\textwidth,
					xmode=log
					]
					\nextgroupplot[
					xmode=log,
					title  = {Based on function evalutions},
					ymin=-0.01,ymax=1.01,
					ytick distance=0.1,
					xmin=9,xmax=2000000,
					xtick distance=10,
					xlabel = {Function evaluations},
					ylabel = {Proportion of problems solved},
					]
					\addplot[postaction={decoration={markings,mark=between positions 0 and 1 step 0.1 with {\node[circle,draw=darkblue,fill=darkblue,inner sep=1.25pt] {};}},decorate,},darkblue,line width=0.75pt] table[x=T,y=NAS] {figure_data/General.txt};
					\addplot[postaction={decoration={markings,mark=between positions 0 and 1 step 0.1 with {\node[star,star points=4,draw=blueryb,fill=blueryb,inner sep=1.25pt] {};}},decorate,},blueryb,line width=0.75pt] table[x=T,y=BAR] {figure_data/General.txt};
					\addplot[postaction={decoration={markings,mark=between positions 0 and 1 step 0.1 with {\node[draw=DarkRed,fill=DarkRed,inner sep=1.25pt] {};}},decorate,},DarkRed,line width=0.75pt] table[x=T,y=SU2] {figure_data/General.txt};
					\addplot[postaction={decoration={markings,mark=between positions 0 and 1 step 0.1 with {\node[regular polygon,regular polygon sides=5,draw=red,fill=red,inner sep=1.25pt] {};}},decorate,},red,line width=0.75pt] table[x=T,y=SU3] {figure_data/General.txt};
					\addplot[postaction={decoration={markings,mark=between positions 0 and 1 step 0.1 with {\node[diamond,draw=armygreen,fill=armygreen,inner sep=1.25pt] {};}},decorate,},armygreen,line width=0.75pt] table[x=T,y=SU5] {figure_data/General.txt};
					\addplot[postaction={decoration={markings,mark=between positions 0 and 1 step 0.1 with {\node[star,star points=9,fill=antiquebrass,draw=antiquebrass,inner sep=1.25pt] {};}},decorate,},antiquebrass,line width=0.75pt] table[x=T,y=H] {figure_data/General.txt};
					\addplot[postaction={decoration={markings,mark=between positions 0 and 1 step 0.1 with {\node[circle,draw=saffron,inner sep=1.25pt] {};}},decorate,},saffron,line width=0.75pt] table[x=T,y=PL1] {figure_data/General.txt};
					\addplot[postaction={decoration={markings,mark=between positions 0 and 1 step 0.1 with {\node[diamond,draw=schoolbusyellow,inner sep=1.25pt] {};}},decorate,},schoolbusyellow,line width=0.75pt] table[x=T,y=PL2] {figure_data/General.txt};
					\addplot[postaction={decoration={markings,mark=between positions 0 and 1 step 0.1 with {\node[regular polygon,regular polygon sides=5,draw=skyblue,inner sep=1.25pt] {};}},decorate,},skyblue,line width=0.75pt] table[x=T,y=PL3] {figure_data/General.txt};
					\addplot[postaction={decoration={markings,mark=between positions 0 and 1 step 0.1 with {\node[draw=forestgreen,inner sep=1.25pt] {};}},decorate,},forestgreen,line width=0.75pt] table[x=T,y=C] {figure_data/General.txt};
					\addplot[postaction={decoration={markings,mark=between positions 0 and 1 step 0.1 with {\node[star,star points=7,draw=ufogreen,inner sep=1.25pt] {};}},decorate,},ufogreen,line width=0.75pt] table[x=T,y=CE] {figure_data/General.txt};
					\addplot[postaction={decoration={markings,mark=between positions 0 and 1 step 0.1 with {\node[star,star points=5,fill=ao,draw=ao,inner sep=1.25pt] {};}},decorate,},ao,line width=0.75pt] table[x=T,y=CEM] {figure_data/General.txt};
					\addplot[black,dashed,line width=0.5pt] table[x=T,y=TS] {figure_data/General.txt};
					\addplot[black,densely dotted,line width=0.5pt] table[x=T,y=TC] {figure_data/General.txt};
					\nextgroupplot[
					legend pos=outer north east,
					xmode=log,
					legend style={font=\tiny},
					legend style={draw=none},
					legend style={row sep=0.5pt},
					title  = {Based on the execution time},
					ymin=-0.01,ymax=1.01,
					ytick distance=0.1,
					xmin=0.01,xmax=43000,
					xtick distance=10,
					xlabel = {Seconds},
					]
					\addlegendimage{legend image code/.code={\node [circle,draw=darkblue,fill=darkblue,inner sep=1.25pt] {};}}
					\addlegendentry{\directnas}
					\addlegendimage{legend image code/.code={\node [star,star points=4,draw=blueryb,fill=blueryb,inner sep=1.25pt] {};}}
					\addlegendentry{\directbarrier}
					\addlegendimage{legend image code/.code={\node [draw=DarkRed,fill=DarkRed,inner sep=1.25pt] {};}}
					\addlegendentry{\directsub{} ($sub=2$)}
					\addlegendimage{legend image code/.code={\node [regular polygon,regular polygon sides=5,draw=red,fill=red,inner sep=1.25pt] {};}}
					\addlegendentry{\directsub{} ($sub=3$)}
					\addlegendimage{legend image code/.code={\node [diamond,draw=armygreen,fill=armygreen,inner sep=1.25pt] {};}}
					\addlegendentry{\directsub{} ($sub=5$)}
					\addlegendimage{legend image code/.code={\node [star,star points=9,fill=antiquebrass,draw=antiquebrass,inner sep=1.25pt] {};}}
					\addlegendentry{\directh}
					\addlegendimage{legend image code/.code={\node [circle,draw=saffron,inner sep=1.25pt] {};}}
					\addlegendentry{\directll{} ($\gamma=10$)}
					\addlegendimage{legend image code/.code={\node [diamond,draw=schoolbusyellow,inner sep=1.25pt] {};}}
					\addlegendentry{\directll{} ($\gamma=10^2$)}
					\addlegendimage{legend image code/.code={\node [regular polygon,regular polygon sides=5,draw=skyblue,inner sep=1.25pt] {};}}
					\addlegendentry{\directll{} ($\gamma=10^3$)}
					\addlegendimage{legend image code/.code={\node [draw=forestgreen,inner sep=1.25pt] {};}}
					\addlegendentry{\directc}
					\addlegendimage{legend image code/.code={\node [star,star points=7,draw=ufogreen,inner sep=1.25pt] {};}}
					\addlegendentry{\directce}
					\addlegendimage{legend image code/.code={\node [star,star points=5,fill=ao,draw=ao,inner sep=1.25pt] {};}}
					\addlegendentry{\directcemin}
					\addlegendimage{legend image code/.code={\node [draw=none,fill=none] {- - -};}}
					\addlegendentry{\glsolve}
					\addlegendimage{legend image code/.code={\node [draw=none,fill=none] {.......};}}
					\addlegendentry{\glccluster}
					\addplot[postaction={decoration={markings,mark=between positions 0 and 1 step 0.1 with {\node[circle,draw=darkblue,fill=darkblue,inner sep=1.25pt] {};}},decorate,},darkblue,line width=0.75pt] table[x=T,y=NAS] {figure_data/Hidden.txt};
					\addplot[postaction={decoration={markings,mark=between positions 0 and 1 step 0.1 with {\node[star,star points=4,draw=blueryb,fill=blueryb,inner sep=1.25pt] {};}},decorate,},blueryb,line width=0.75pt] table[x=T,y=BAR] {figure_data/Hidden.txt};
					\addplot[postaction={decoration={markings,mark=between positions 0 and 1 step 0.1 with {\node[draw=DarkRed,fill=DarkRed,inner sep=1.25pt] {};}},decorate,},DarkRed,line width=0.75pt] table[x=T,y=SU2] {figure_data/Hidden.txt};
					\addplot[postaction={decoration={markings,mark=between positions 0 and 1 step 0.1 with {\node[regular polygon,regular polygon sides=5,draw=red,fill=red,inner sep=1.25pt] {};}},decorate,},red,line width=0.75pt] table[x=T,y=SU3] {figure_data/Hidden.txt};
					\addplot[postaction={decoration={markings,mark=between positions 0 and 1 step 0.1 with {\node[diamond,draw=armygreen,fill=armygreen,inner sep=1.25pt] {};}},decorate,},armygreen,line width=0.75pt] table[x=T,y=SU5] {figure_data/Hidden.txt};
					\addplot[postaction={decoration={markings,mark=between positions 0 and 1 step 0.1 with {\node[star,star points=9,fill=antiquebrass,draw=antiquebrass,inner sep=1.25pt] {};}},decorate,},antiquebrass,line width=0.75pt] table[x=T,y=H] {figure_data/Hidden.txt};
					\addplot[postaction={decoration={markings,mark=between positions 0 and 1 step 0.1 with {\node[circle,draw=saffron,inner sep=1.25pt] {};}},decorate,},saffron,line width=0.75pt] table[x=T,y=PL1] {figure_data/Hidden.txt};
					\addplot[postaction={decoration={markings,mark=between positions 0 and 1 step 0.1 with {\node[diamond,draw=schoolbusyellow,inner sep=1.25pt] {};}},decorate,},schoolbusyellow,line width=0.75pt] table[x=T,y=PL2] {figure_data/Hidden.txt};
					\addplot[postaction={decoration={markings,mark=between positions 0 and 1 step 0.1 with {\node[regular polygon,regular polygon sides=5,draw=skyblue,inner sep=1.25pt] {};}},decorate,},skyblue,line width=0.75pt] table[x=T,y=PL3] {figure_data/Hidden.txt};
					\addplot[postaction={decoration={markings,mark=between positions 0 and 1 step 0.1 with {\node[draw=forestgreen,inner sep=1.25pt] {};}},decorate,},forestgreen,line width=0.75pt] table[x=T,y=C] {figure_data/Hidden.txt};
					\addplot[postaction={decoration={markings,mark=between positions 0 and 1 step 0.1 with {\node[star,star points=7,draw=ufogreen,inner sep=1.25pt] {};}},decorate,},ufogreen,line width=0.75pt] table[x=T,y=CE] {figure_data/Hidden.txt};
					\addplot[postaction={decoration={markings,mark=between positions 0 and 1 step 0.1 with {\node[star,star points=5,fill=ao,draw=ao,inner sep=1.25pt] {};}},decorate,},ao,line width=0.75pt] table[x=T,y=CEM] {figure_data/Hidden.txt};
					\addplot[black,dashed,line width=0.5pt] table[x=T,y=TS] {figure_data/Hidden.txt};
					\addplot[black,densely dotted,line width=0.5pt] table[x=T,y=TC] {figure_data/Hidden.txt};
				\end{groupplot}
		\end{tikzpicture}}
		\caption{Data profiles of \direct-type algorithms from \toolbox{} and \texttt{TOMLAB} on the whole set of constrained optimization test problems.}
		\label{fig:general}
	\end{figure}
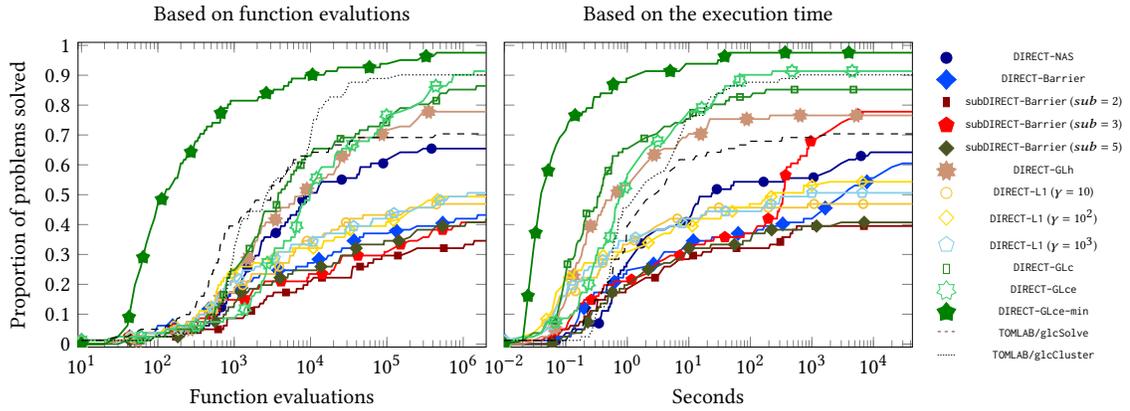

	Finally, in \cref{fig:general}, the comparative performance using the data profiles tool is demonstrated.
	Data profiles confirm the same trends, i.e., the hybridized versions are the best performing, and the overall advantage of methods that incorporate constrained information versus designed explicitly for problems with hidden constraints (see also \cref{fig:Nonlinear}).

	\subsubsection{Test results on problems with linear constraints}

\begin{figure}
		\resizebox{\textwidth}{!}{
			\begin{tikzpicture}
				\begin{groupplot}[
					group style={
						group size=2 by 1,
						x descriptions at=edge bottom,
						y descriptions at=edge left,
						horizontal sep=7.5pt,
					},
					height=0.4\textwidth,width=0.5\textwidth,
					xmode=log
					]
					\nextgroupplot[
					xmode=log,
					title  = {Based on function evalutions},
					ymin=-0.01,ymax=1.05,
					ytick distance=0.1,
					xmin=0.9,xmax=2000000,
					xtick distance=10,
					xlabel = {Function evaluations},
					ylabel = {Proportion of problems solved},
					]
					\addplot[postaction={decoration={markings,mark=between positions 0 and 1 step 0.1 with {\node[circle,draw=darkblue,fill=darkblue,inner sep=1.25pt] {};}},decorate,},darkblue,line width=0.75pt] table[x=T,y=NAS] {figure_data/Linear.txt};
					\addplot[postaction={decoration={markings,mark=between positions 0 and 1 step 0.1 with {\node[star,star points=4,draw=blueryb,fill=blueryb,inner sep=1.25pt] {};}},decorate,},blueryb,line width=0.75pt] table[x=T,y=BAR] {figure_data/Linear.txt};
					\addplot[postaction={decoration={markings,mark=between positions 0 and 1 step 0.1 with {\node[draw=DarkRed,fill=DarkRed,inner sep=1.25pt] {};}},decorate,},DarkRed,line width=0.75pt] table[x=T,y=SU2] {figure_data/Linear.txt};
					\addplot[postaction={decoration={markings,mark=between positions 0 and 1 step 0.1 with {\node[regular polygon,regular polygon sides=5,draw=red,fill=red,inner sep=1.25pt] {};}},decorate,},red,line width=0.75pt] table[x=T,y=SU3] {figure_data/Linear.txt};
					\addplot[postaction={decoration={markings,mark=between positions 0 and 1 step 0.1 with {\node[diamond,draw=armygreen,fill=armygreen,inner sep=1.25pt] {};}},decorate,},armygreen,line width=0.75pt] table[x=T,y=SU5] {figure_data/Linear.txt};
					\addplot[postaction={decoration={markings,mark=between positions 0 and 1 step 0.1 with {\node[star,star points=9,fill=antiquebrass,draw=antiquebrass,inner sep=1.25pt] {};}},decorate,},antiquebrass,line width=0.75pt] table[x=T,y=H] {figure_data/Linear.txt};
					\addplot[postaction={decoration={markings,mark=between positions 0 and 1 step 0.1 with {\node[circle,draw=saffron,inner sep=1.25pt] {};}},decorate,},saffron,line width=0.75pt] table[x=T,y=PL1] {figure_data/Linear.txt};
					\addplot[postaction={decoration={markings,mark=between positions 0 and 1 step 0.1 with {\node[diamond,draw=schoolbusyellow,inner sep=1.25pt] {};}},decorate,},schoolbusyellow,line width=0.75pt] table[x=T,y=PL2] {figure_data/Linear.txt};
					\addplot[postaction={decoration={markings,mark=between positions 0 and 1 step 0.1 with {\node[regular polygon,regular polygon sides=5,draw=skyblue,inner sep=1.25pt] {};}},decorate,},skyblue,line width=0.75pt] table[x=T,y=PL3] {figure_data/Linear.txt};
					\addplot[postaction={decoration={markings,mark=between positions 0 and 1 step 0.1 with {\node[draw=forestgreen,inner sep=1.25pt] {};}},decorate,},forestgreen,line width=0.75pt] table[x=T,y=C] {figure_data/Linear.txt};
					\addplot[postaction={decoration={markings,mark=between positions 0 and 1 step 0.1 with {\node[star,star points=7,draw=ufogreen,inner sep=1.25pt] {};}},decorate,},ufogreen,line width=0.75pt] table[x=T,y=CE] {figure_data/Linear.txt};
					\addplot[postaction={decoration={markings,mark=between positions 0 and 1 step 0.1 with {\node[star,star points=5,fill=ao,draw=ao,inner sep=1.25pt] {};}},decorate,},ao,line width=0.75pt] table[x=T,y=CEM] {figure_data/Linear.txt};
					\addplot[black,dashed,line width=0.5pt] table[x=T,y=TS] {figure_data/Linear.txt};
					\addplot[black,densely dotted,line width=0.5pt] table[x=T,y=TC] {figure_data/Linear.txt};
					\addplot[postaction={decoration={markings,mark=between positions 0 and 1 step 0.15 with {\node[regular polygon,regular polygon sides=9,fill=schoolbusyellow,draw=schoolbusyellow,inner sep=1.25pt] {};}},decorate,},schoolbusyellow,line width=0.75pt] table[x=T,y=DC] {figure_data/Linear.txt};
					\addplot[postaction={decoration={markings,mark=between positions 0 and 1 step 0.15 with {\node[star,star points=9,draw=selectiveyellow,inner sep=1.25pt] {};}},decorate,},selectiveyellow,line width=0.75pt] table[x=T,y=DV] {figure_data/Linear.txt};
					\nextgroupplot[
					legend pos=outer north east,
					xmode=log,
					legend style={font=\tiny},
					legend style={draw=none},
					legend style={row sep=0.5pt},
					title  = {Based on the execution time},
					ymin=-0.01,ymax=1.05,
					ytick distance=0.1,
					xmin=0.01,xmax=43000,
					xtick distance=10,
					xlabel = {Seconds},
					]
					\addlegendimage{legend image code/.code={\node [circle,draw=darkblue,fill=darkblue,inner sep=1.25pt] {};}}
					\addlegendentry{\directnas}
					\addlegendimage{legend image code/.code={\node [star,star points=4,draw=blueryb,fill=blueryb,inner sep=1.25pt] {};}}
					\addlegendentry{\directbarrier}
					\addlegendimage{legend image code/.code={\node [draw=DarkRed,fill=DarkRed,inner sep=1.25pt] {};}}
					\addlegendentry{\directsub{} ($sub=2$)}
					\addlegendimage{legend image code/.code={\node [regular polygon,regular polygon sides=5,draw=red,fill=red,inner sep=1.25pt] {};}}
					\addlegendentry{\directsub{} ($sub=3$)}
					\addlegendimage{legend image code/.code={\node [diamond,draw=armygreen,fill=armygreen,inner sep=1.25pt] {};}}
					\addlegendentry{\directsub{} ($sub=5$)}
					\addlegendimage{legend image code/.code={\node [star,star points=9,fill=antiquebrass,draw=antiquebrass,inner sep=1.25pt] {};}}
					\addlegendentry{\directh}
					\addlegendimage{legend image code/.code={\node [circle,draw=saffron,inner sep=1.25pt] {};}}
					\addlegendentry{\directll{} ($\gamma=10$)}
					\addlegendimage{legend image code/.code={\node [diamond,draw=schoolbusyellow,inner sep=1.25pt] {};}}
					\addlegendentry{\directll{} ($\gamma=10^2$)}
					\addlegendimage{legend image code/.code={\node [regular polygon,regular polygon sides=5,draw=skyblue,inner sep=1.25pt] {};}}
					\addlegendentry{\directll{} ($\gamma=10^3$)}
					\addlegendimage{legend image code/.code={\node [draw=forestgreen,inner sep=1.25pt] {};}}
					\addlegendentry{\directc}
					\addlegendimage{legend image code/.code={\node [star,star points=7,draw=ufogreen,inner sep=1.25pt] {};}}
					\addlegendentry{\directce}
					\addlegendimage{legend image code/.code={\node [star,star points=5,fill=ao,draw=ao,inner sep=1.25pt] {};}}
					\addlegendentry{\directcemin}
					\addlegendimage{legend image code/.code={\node [draw=none,fill=none] {- - -};}}
					\addlegendentry{\glsolve}
					\addlegendimage{legend image code/.code={\node [draw=none,fill=none] {.......};}}
					\addlegendentry{\glccluster}
					\addlegendimage{legend image code/.code={\node [regular polygon,regular polygon sides=9,fill=schoolbusyellow,draw=schoolbusyellow,inner sep=1.25pt] {};}}
					\addlegendentry{\disimpllc}
					\addlegendimage{legend image code/.code={\node [star,star points=9,draw=selectiveyellow,inner sep=1.25pt] {};}}
					\addlegendentry{\disimpllv}
					\addplot[postaction={decoration={markings,mark=between positions 0 and 1 step 0.1 with {\node[circle,draw=darkblue,fill=darkblue,inner sep=1.25pt] {};}},decorate,},darkblue,line width=0.75pt] table[x=T,y=NAS] {figure_data/Nonlinear.txt};
					\addplot[postaction={decoration={markings,mark=between positions 0 and 1 step 0.1 with {\node[star,star points=4,draw=blueryb,fill=blueryb,inner sep=1.25pt] {};}},decorate,},blueryb,line width=0.75pt] table[x=T,y=BAR] {figure_data/Nonlinear.txt};
					\addplot[postaction={decoration={markings,mark=between positions 0 and 1 step 0.1 with {\node[draw=DarkRed,fill=DarkRed,inner sep=1.25pt] {};}},decorate,},DarkRed,line width=0.75pt] table[x=T,y=SU2] {figure_data/Nonlinear.txt};
					\addplot[postaction={decoration={markings,mark=between positions 0 and 1 step 0.1 with {\node[regular polygon,regular polygon sides=5,draw=red,fill=red,inner sep=1.25pt] {};}},decorate,},red,line width=0.75pt] table[x=T,y=SU3] {figure_data/Nonlinear.txt};
					\addplot[postaction={decoration={markings,mark=between positions 0 and 1 step 0.1 with {\node[diamond,draw=armygreen,fill=armygreen,inner sep=1.25pt] {};}},decorate,},armygreen,line width=0.75pt] table[x=T,y=SU5] {figure_data/Nonlinear.txt};
					\addplot[postaction={decoration={markings,mark=between positions 0 and 1 step 0.1 with {\node[star,star points=9,fill=antiquebrass,draw=antiquebrass,inner sep=1.25pt] {};}},decorate,},antiquebrass,line width=0.75pt] table[x=T,y=H] {figure_data/Nonlinear.txt};
					\addplot[postaction={decoration={markings,mark=between positions 0 and 1 step 0.1 with {\node[circle,draw=saffron,inner sep=1.25pt] {};}},decorate,},saffron,line width=0.75pt] table[x=T,y=PL1] {figure_data/Nonlinear.txt};
					\addplot[postaction={decoration={markings,mark=between positions 0 and 1 step 0.1 with {\node[diamond,draw=schoolbusyellow,inner sep=1.25pt] {};}},decorate,},schoolbusyellow,line width=0.75pt] table[x=T,y=PL2] {figure_data/Nonlinear.txt};
					\addplot[postaction={decoration={markings,mark=between positions 0 and 1 step 0.1 with {\node[regular polygon,regular polygon sides=5,draw=skyblue,inner sep=1.25pt] {};}},decorate,},skyblue,line width=0.75pt] table[x=T,y=PL3] {figure_data/Nonlinear.txt};
					\addplot[postaction={decoration={markings,mark=between positions 0 and 1 step 0.1 with {\node[draw=forestgreen,inner sep=1.25pt] {};}},decorate,},forestgreen,line width=0.75pt] table[x=T,y=C] {figure_data/Nonlinear.txt};
					\addplot[postaction={decoration={markings,mark=between positions 0 and 1 step 0.1 with {\node[star,star points=7,draw=ufogreen,inner sep=1.25pt] {};}},decorate,},ufogreen,line width=0.75pt] table[x=T,y=CE] {figure_data/Nonlinear.txt};
					\addplot[postaction={decoration={markings,mark=between positions 0 and 1 step 0.1 with {\node[star,star points=5,fill=ao,draw=ao,inner sep=1.25pt] {};}},decorate,},ao,line width=0.75pt] table[x=T,y=CEM] {figure_data/Nonlinear.txt};
					\addplot[black,dashed,line width=0.5pt] table[x=T,y=TS] {figure_data/Nonlinear.txt};
					\addplot[black,densely dotted,line width=0.5pt] table[x=T,y=TC] {figure_data/Nonlinear.txt};
					\addplot[postaction={decoration={markings,mark=between positions 0 and 1 step 0.15 with {\node[regular polygon,regular polygon sides=9,fill=schoolbusyellow,draw=schoolbusyellow,inner sep=1.25pt] {};}},decorate,},schoolbusyellow,line width=0.75pt] table[x=T,y=DC] {figure_data/Nonlinear.txt};
					\addplot[postaction={decoration={markings,mark=between positions 0 and 1 step 0.15 with {\node[star,star points=9,draw=selectiveyellow,inner sep=1.25pt] {};}},decorate,},selectiveyellow,line width=0.75pt] table[x=T,y=DV] {figure_data/Nonlinear.txt};
				\end{groupplot}
		\end{tikzpicture}}
		\caption{Data profiles of \direct-type algorithms for problems with constraints on the subset of constrained optimization test problems containing linear constraints.}
		\label{fig:linear}
	\end{figure}
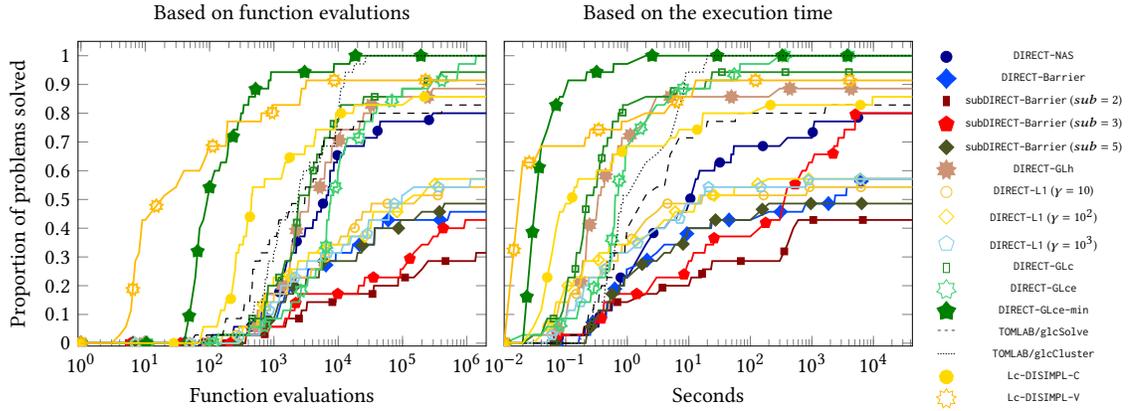

	In the final part, we test the performance of \direct-type algorithms on problems with linear constraints.
	We consider all previously tested algorithms and two specifically designed simplicial partitions-based \disimpllv{} and \disimpllc~\cite{Paulavicius2016:ol}.
	The main advantage of simplices is that they can cover a feasible region defined by linear constraints.
	Thus any infeasible areas are not involved in the search.
	Moreover, for most problems from \directlib, the solution is located at the intersection of linear constraints.
	Therefore \disimpllv{} finds it in the early or even in the first iteration.
	The data profiles (see \cref{fig:linear}) reveal the overall effectiveness of the \disimpllv{} algorithm for such problems with linear constraints.
	The latter algorithm even outperformed the hybridized \directcemin{} method, which is surprisingly enough.
	Nevertheless, the efficiency of simplicial partition-based algorithms suffers from the problem dimension.
	Three higher dimensionality linearly constrained problems were unsolved by the \disimpllv{} algorithm (see the bottom part of~\cref{tab:results1}, `Performance of \direct-type algorithms devoted to problems with linear constraints only').
	Moreover, simplicial partition-based implementations are pretty slow.
	For example, on average, \directce{} is approximately $307$ times faster than \disimpllv.
	Moreover, \directce, \directcemin, and \glccluster{} solved all test problems with linear constraints.

	\section{\toolbox{} performance on engineering problems}
	\label{sec:aplications}

	In this section, the algorithms from \toolbox{} and \texttt{TOMLAB} were tested on eleven engineering design problems: tension/compression spring, three-bar truss, NASA speed reducer, pressure vessel, welded beam, and six different versions of the general non-linear regression problem.
	The general non-linear regression problem is box-constrained, while the others involve different constraints.
	Only the most promising algorithms (based on~\cref{sec:experimentstests}) were considered.
	We used the same stopping rule as the global minimums are known for all these engineering problems.
  A detailed description of all engineering problems and mathematical formulations is given in \cref{sec:apendix}.

	\subsection{Tension/compression spring design problem}

	Here, we consider the tension-compression string design problem.
	This problem aims to minimize the string weight under the constraints on deflection, shear stress, surge frequency, and limits on the outside diameter. A detailed description of the practical problem can be found in \cite{Kazemi2011}, while in \cref{prob:Tension},
	we give a short description and mathematical formulation.

	A comparison of found solutions and performance metrics by the algorithms from \toolbox{} and \tomlab{} is shown in \cref{tab:tennison}.
	Four considered algorithms were able to solve this problem.
	Based on the number of function evaluation criteria, a slightly unexpectedly \directnas{} was significantly better than other algorithms.
	The surprise is that the algorithm does not use any information about the constraint functions.
	However, the \directnas{} algorithm was approximate four times slower than the second-best method \directce.
	Another surprise was that none of the hybridized algorithms performed well on this practical problem.
	Moreover, the \glccluster{} algorithm failed to find a solution within the given time limit.

	\begin{table}[ht]
		\caption{Performance of the \direct-type algorithms from \toolbox{} and \tomlab{} on a tension/compression design problem}
		\begin{tabular*}{\textwidth}{@{\extracolsep{\fill}}l|l|r|r|r|r}
			\toprule
			Algorithm (input)   & Parameter & Iterations & $f_{eval}$ & Time (s) & $f_{min}$	\\
			\midrule
			\directnas     			& -- & $297$    & $17,659$ 	 & $77.46$    & $0.012680$ \\
			\directbarrier 			& -- & $40,478$ & \textcolor{red}{$>2 \times 10^6$} & $2,465.81$ & \textcolor{red}{$0.012867$} \\
			\directsub{} & $sub = 2$ & $256$ & \textcolor{red}{$>2 \times 10^6$}   & $207.89$ & \textcolor{red}{$0.012708$} 	\\
			\directsub{} & $sub = 3$ & $2,187$ & \textcolor{red}{$>2 \times 10^6$}   & $198.94$ & \textcolor{red}{$0.012741$}		\\
			\directsub{} & $sub = 5$ & $15,626$ & \textcolor{red}{$>2 \times 10^6$} & $722.79$   & \textcolor{red}{$0.012867$} \\
			\directh       			& -- & $700$    & \textcolor{red}{$>2 \times 10^6$} & $214.81$   & \textcolor{red}{$0.012683$} \\
			\directc        		& -- & $675$    & $423,209$ 	& $47.53$    & $0.012680$ \\
			\directce       		& -- & $624$    & $178,115$   	& $20.45$    & $0.012680$ \\
			\directcemin$^{*}$    		& -- & $624$    & $178,115$   	& $20.70$    & $0.012680$ \\
			\directll{} & $\gamma = 10^1$ & $34,395$  & \textcolor{red}{$>2 \times 10^6$} & $3,265.38$ & \textcolor{red}{$0.012755$}    	\\
			\directll{} & $\gamma = 10^2$ & $33,897$  & \textcolor{red}{$>2 \times 10^6$} & $3,147.87$ & \textcolor{red}{$0.012867$}		\\
			\directll{} & $\gamma = 10^3$ & $33,571$  & \textcolor{red}{$>2 \times 10^6$} & $3,041.00$ & \textcolor{red}{$0.012867$} 	\\
			\glsolve        	& -- & $64$   & \textcolor{red}{$>2 \times 10^6$}      & $560.36$     	& \textcolor{red}{$0.014669$}    \\
			\glccluster$^{*}$        	& -- & $12$   & $1,528,205$      & \textcolor{red}{$>43,000.00$}     		& \textcolor{red}{$0.014669$}    \\
			\bottomrule
			\multicolumn{6}{l}{$^{*}$ -- a hybrid version of the algorithm, enriched with the local search subroutine}
		\end{tabular*}
		\label{tab:tennison}
	\end{table}

	\subsection{Three-bar truss design problem}

	Here, we consider the three-bar truss design problem.
	The goal is to minimize the volume subject to stress constraints.
	A detailed description of the problem is given in \cite{Ray2003}, while in \cref{prob:Three-bar}, we provide a brief description and mathematical formulation.

	A comparison of found solutions and performance metrics is shown in \cref{tab:three}.
	Here, hybridized \directcemin{} was the most efficient optimizer.
	However, none of the algorithms (with the proper input parameters) had any difficulty solving this problem, and they found the solution in less than one-second time.

	\begin{table}[ht]
		\caption{Performance of the \direct-type algorithms from \toolbox{} and \tomlab{} on a three-bar truss design problem}
		\begin{tabular*}{\textwidth}{@{\extracolsep{\fill}}l|l|r|r|r|r}
			\toprule
			Algorithm (input) & Parameter & Iterations & $f_{eval}$ & Time (s) & $f_{min}$	\\
			\midrule
			\directnas     		& -- & $29$ & $339$   & $0.05$ & $263.915790$ 	\\
			\directbarrier 		& -- & $13$ & $125$   & $0.02$ & $263.915790$ 	\\
			\directsub{} & $sub = 2$ & $512$ & \textcolor{red}{$>2 \times 10^6$}   & $141.63$ & \textcolor{red}{$283.223781$} 	\\
			\directsub{} & $sub = 3$ & $17$ & $333$   & $0.02$ & $263.915800$ 	\\
			\directsub{} & $sub = 5$ & $14$ & $161$   & $0.02$ & $263.915790$ 	\\
			\directh       		& -- & $12$ & $231$   & $0.03$ & $263.915790$    \\
			\directc        	& -- & $17$ & $727$   & $0.08$ & $263.911750$    \\
			\directce       	& -- & $33$ & $1,055$ & $0.13$ & $263.915790$    \\
			\directcemin$^{*}$    	& -- & $6$  & $93$    & $0.03$ & $263.895850$ 	\\
			\directll{} & $\gamma = 10^1$ & $1$  & $1$     & $0.01$ & \textcolor{red}{$199.705600^a$}    \\
			\directll{} & $\gamma = 10^2$ & $2$  & $11$    & $0.01$ & \textcolor{red}{$262.344700^a$}    \\
			\directll{} & $\gamma = 10^3$ & $17$ & $179$   & $0.03$ & $263.915790$    \\
			\glsolve        	& -- & $25$   & $647$      & $0.43$     & $263.910482$    \\
			\glccluster$^{*}$        	& -- & $1$   & $992$      & $0.55$     & $263.910482$    \\
			\bottomrule
			\multicolumn{6}{l}{$^{*}$ -- a hybrid version of the algorithm, enriched with the local search subroutine} \\
			\multicolumn{5}{l}{a -- result is outside the feasible region}
		\end{tabular*}
		\label{tab:three}
	\end{table}

	\subsection{NASA speed reducer design problem}

	Here we consider the NASA speed reducer design problem. The goal is to minimize the overall weight subject to constraints on the gear teeth' bending stress, surface stress, transverse deflection of the shaft, and stresses in the shafts.
	A detailed description of the problem can be found in \cite{Ray2003}, while in \cref{prob:nasa}, we provide a short description and mathematical formulation.

	A comparison of the found solutions and performance metrics is shown in \cref{tab:nasa}.
	Only three algorithms (\directce, \directcemin, and \glccluster) were able to tackle this problem.
	Again, the hybridized \directcemin{} algorithm showed the best performance.
	Note that the found solutions with \directll{} are better than the best-known value $f_{min}$.
	However, the reported solution point is outside the feasible region and violates some constraints.
	The \glsolve{} was very close to a solution, but could not find it with a required accuracy within the maximum time limit.

	\begin{table}[ht]
		\caption{Performance of the \direct-type algorithms from \toolbox{} and \tomlab{} on a NASA speed reducer design problem}
		\begin{tabular*}{\textwidth}{@{\extracolsep{\fill}}l|l|r|r|r|r}
			\toprule
			Algorithm (input) & Parameter & Iterations & $f_{eval}$ & Time (s) & $f_{min}$	\\
			\midrule
			\directnas     	& -- 		& $6,719$  	& $325,691$ 	 					& \textcolor{red}{$>43,000.00$} 	& \textcolor{red}{$3,006.874789$}  	\\
			\directbarrier 	& -- 		& $32,031$ 	& \textcolor{red}{$>2 \times 10^6$} & $2,625.88$  	& \textcolor{red}{$3,006.838136$} 	\\
			\directsub{} 	& $sub = 2$ & $64$ 		& \textcolor{red}{$>2 \times 10^6$} & $389.49$ 		& \textcolor{red}{$3,045.559256$}  	\\
			\directsub{} 	& $sub = 3$ & $243$ 	& \textcolor{red}{$>2 \times 10^6$} & $85.07$ 		& \textcolor{red}{$3,040.464809$} 	\\
			\directsub{} 	& $sub = 5$ & $11,112$ 	& \textcolor{red}{$>2 \times 10^6$} & $1,234.33$  	& \textcolor{red}{$3,006.838136$} 	\\
			\directh       	& -- 		& $528$    	& \textcolor{red}{$>2 \times 10^6$} & $169.74$    	& \textcolor{red}{$3,003.135167$} 	\\
			\directc       	& -- 		& $1,623$   & \textcolor{red}{$>2 \times 10^6$} & $373.40$    	& \textcolor{red}{$3,002.869474$}   	\\
			\directce      	& -- 		& $254$    	& $123,175$ 	 					& $9.70$  	   	& $2,996.572800$   \\
			\directcemin$^{*}$   	& -- 		& $55$     	& $10,229$ 	 						& $0.93$  	   	& $2,996.348212$   \\
			\directll{} & $\gamma = 10^1$ 	& $4$  		& $67$     							& $0.01$	 	& \textcolor{red}{$2,943.869936^a$}    \\
			\directll{} & $\gamma = 10^2$ 	& $4$  		& $67$    							& $0.02$ 		& \textcolor{red}{$2,982.462161^a$}    \\
			\directll{} & $\gamma = 10^3$ 	& $661$    	& $26,667$    						& $3.01$ 	   	& \textcolor{red}{$2,995.382568^a$}  \\
			\glsolve      	& -- 		& $66,305$   		& $1,504,889$      							& \textcolor{red}{$>43,000.00$}     		& \textcolor{red}{$2,996.659779$}    \\
			\glccluster$^{*}$    	& -- 		& $1$   	& $12,736$      			& $7.35$     		& $2,996.347954$    \\
			\bottomrule
			\multicolumn{6}{l}{$^{*}$ -- a hybrid version of the algorithm, enriched with the local search subroutine} \\
			\multicolumn{5}{l}{a -- result is outside the feasible region}
		\end{tabular*}
		\label{tab:nasa}
	\end{table}

	\subsection{Pressure vessel design problem}

	In this subsection, we consider a pressure vessel design problem, and the goal is to minimize the total cost of the material, form, and weld a cylindrical vessel.
	A detailed description of the problem can be found in \cite{Kazemi2011}, while in \cref{prob:Pressure}, we provide a short description and mathematical formulation.

	A comparison of the found solutions and performance metrics is shown in \cref{tab:pressure}.
	Five algorithms solved this problem: \directnas, \directh, \directce, \glccluster, and \directcemin{} was the most efficient optimizer again.
	\directnas{} is the best performing and outperformed the second-best by approximately $1.8$ times, among traditional \direct-type algorithms.
	However, the \directnas{} algorithm was about $26$ times slower than the second-best method (\directh).
	As in the previous case, the \directll{} returned a better than the best-know value $f_{min}$, but the solution points lay outside the feasible region.
	\begin{table}[ht]
		\caption{Performance of the \direct-type algorithms from \toolbox{} and \tomlab{} on a pressure vessel design problem}
		\begin{tabular*}{\textwidth}{@{\extracolsep{\fill}}l|l|r|r|r|r}
			\toprule
			Algorithm (input) & Parameter & Iterations & $f_{eval}$ & Time (s) & $f_{min}$	\\
			\midrule
			\directnas     	& -- 		& $273$    	& $31,081$ 	 						& $126.64$    	& $7,164.437307$ 	\\
			\directbarrier 	& -- 		& $26,379$ 	& \textcolor{red}{$>2 \times 10^6$} & $1,600.15$ 	& \textcolor{red}{$7,234.041903$}  \\
			\directsub{} 	& $sub = 2$ & $128$ 	& \textcolor{red}{$>2 \times 10^6$} & $190.11$ 		& \textcolor{red}{$7,234.402264$}  	\\
			\directsub{} 	& $sub = 3$ & $729$ 	& \textcolor{red}{$>2 \times 10^6$} & $99.05$ 		& \textcolor{red}{$7,234.222516$} 	\\
			\directsub{} 	& $sub = 5$	& $15,626$ 	& \textcolor{red}{$>2 \times 10^6$}	& $1,328.98$  	& \textcolor{red}{$7,234.041903$}  \\
			\directh       	& -- 		& $252$    	& $55,837$ 	 						& $4.80$  	   	& $7,164.437300$  \\
			\directc        & -- 		& $2,358$  	& \textcolor{red}{$>2 \times 10^6$}	& $433.87$    	& \textcolor{red}{$7,224.704257$}  \\
			\directce       & -- 		& $322$    	& $88,585$ 							& $8.52$  	   	& $7,164.437301$  \\
			\directcemin$^{*}$     & -- 		& $1$  	   	& $134$ 		 						& $0.24$  	   	& $7,163.739570$  \\
			\directll{} & $\gamma = 10^1$ 	& $86$  	& $2,117$     						& $0.34$	 	& \textcolor{red}{$7,025.940549^a$}    \\
			\directll{} & $\gamma = 10^2$ 	& $86$  	& $2,099$    						& $0.33$ 		& \textcolor{red}{$7,037.428049^a$}    \\
			\directll{} & $\gamma = 10^3$    & $87$     	& $2,295$ 	 						& $0.23$  	   	& \textcolor{red}{$7,152.303079^a$}  \\
			\glsolve      	& -- 		& $123$   		& \textcolor{red}{$>2 \times 10^6$}       	& $529.69$     		& \textcolor{red}{$8,260.982616$}    \\
			\glccluster$^{*}$     	& -- 		& $1$   		& $10,026$      							& $5.74$     		& $7163.739569$    \\
			\bottomrule
			\multicolumn{6}{l}{$^{*}$ -- a hybrid version of the algorithm, enriched with the local search subroutine} \\
			\multicolumn{5}{l}{a -- result is outside the feasible region}
		\end{tabular*}
		\label{tab:pressure}
	\end{table}

	\subsection{Welded beam design problem}
	The fifth engineering problem is the welded beam design.
	The goal is to minimize a welded beam for a minimum cost, subject to seven constraints.
	The detailed description is presented in \cite{Mirjalili2016,Mirjalili2014}, while in \cref{prob:Welded}, we provide a short description and mathematical formulation.

	A comparison of the algorithms is shown in \cref{tab:welder}.
	In total, six algorithms were able to solve the problem, and once again, the \directcemin{} was the most efficient one.
	Again, the \directnas{} algorithm showed the best performance (based on the total number of function evaluations) among traditional \direct-type algorithms but significantly suffered based on the execution time.

	\begin{table}[ht]
		\caption{Performance of the \direct-type algorithms from \toolbox{} and \tomlab{} on a welded beam design problem}
		\begin{tabular*}{\textwidth}{@{\extracolsep{\fill}}l|l|r|r|r|r}
			\toprule
			Algorithm (input) & Parameter & Iterations & $f_{eval}$ & Time (s) & $f_{min}$	\\
			\midrule
			\directnas     	& -- 		& $698$     & $86,863$ 							& $4,692.00$  	& $1.724970$  	\\
			\directbarrier 	& -- 		& $22,934$ 	& \textcolor{red}{$>2 \times 10^6$} & $1,363.24$ 	& \textcolor{red}{$1.728488$} 	\\
			\directsub{} 	& $sub = 2$ & $128$		& \textcolor{red}{$>2 \times 10^6$} & $152.58$ 		& \textcolor{red}{$1.728060$}  	\\
			\directsub{} 	& $sub = 3$ & $729$ 	& \textcolor{red}{$>2 \times 10^6$} & $128.37$ 		& \textcolor{red}{$1.728043$} 	\\
			\directsub{} 	& $sub = 5$ & $10,954$ 	& \textcolor{red}{$>2 \times 10^6$} & $860.06$   	& \textcolor{red}{$1.728037$} 	\\
			\directh       	& -- 		& $189$  	& $158,747$    						& $11.70$    	& $1.724970$  	\\
			\directc        & -- 		& $211$ 	& $108,683$   						& $9.19$     	& $1.724970$  	\\
			\directce       & -- 		& $366$  	& $104,191$   						& $9.80$     	& $1.724970$  	\\
			\directcemin$^{*}$    & -- 		& $3$  		& $163$       						& $0.06$     	& $1.724884$ 	\\
			\directll{} & $\gamma = 10^1$ 	& $21,143$  	& \textcolor{red}{$>2 \times 10^6$} & $1,995.01$	 	& \textcolor{red}{$1.728491$}    \\
			\directll{} & $\gamma = 10^2$ 	& $20,879$  	& \textcolor{red}{$>2 \times 10^6$} & $1,814.33$ 		& \textcolor{red}{$1.728491$}    \\
			\directll{} & $\gamma = 10^3$ 	& $20,767$  & \textcolor{red}{$>2 \times 10^6$} & $1,679.76$ 	& \textcolor{red}{$1.728488$} 	\\
			\glsolve      	& -- 		& $86$   		& \textcolor{red}{$>2 \times 10^6$}      							& $526.06$     		& \textcolor{red}{$2.473711$}     \\
			\glccluster$^{*}$    	& -- 		& $1$   	& $9,884$      						& $5.85$     	& $1.724852$    \\
			\bottomrule
			\multicolumn{6}{l}{$^{*}$ -- a hybrid version of the algorithm, enriched with the local search subroutine}
		\end{tabular*}
		\label{tab:welder}
	\end{table}

	\subsection{General non-linear regression problem}

	In the final part, a general non-linear regression design problem is considered in the form of fitting a sum of damped sinusoids to a series of observations.
	The detailed description of the problem can be found in \cite{Gillard2017,Paulavicius2019:eswa,Sergeyev2016a}, while in \cref{prob:general},
	we provide a short description and mathematical formulation.
	The problem is multi-modal and is considered challenging, especially with the increase in the number of samples ($T$).
	The higher number of sinusoids ($\varsigma$) leads to a more accurate but, at the same time, more challenging optimization problem.

	\begin{table}[htb]
		\caption{Performance of the \direct-type algorithms from \toolbox{} and \tomlab{} on a non-linear regression design problem}
		\resizebox{\textwidth}{!}{
			\begin{tabular}[tb]{@{\extracolsep{\fill}}l|rrrr|rrrr|rrrr|}
				\toprule
				Algorithm & Iter. & $f_{eval}$ & Time(s) & $f_{min}$ & Iter. & $f_{eval}$ & Time(s) & $f_{min}$ & Iter. & $f_{eval}$ & Time(s) & $f_{min}$ \\
				\midrule
				\multicolumn{1}{c}{} & \multicolumn{4}{c}{$\varsigma=1, T=10$}  & \multicolumn{4}{c}{$\varsigma=2, T=10$}  & \multicolumn{4}{c}{$\varsigma=3, T=10$} \\
				\midrule
				\direct   & $31$  & $501$    & $0.03$  & $0.000021$ & $9,328$ & \textcolor{red}{$>2 \times 10^6$} & $362.84$ & \textcolor{red}{$0.001228$} & $3,186$  & \textcolor{red}{$>2 \times 10^6$} & $215.58$   & \textcolor{red}{$0.007295$}  \\
				\directrest & $31$ & $501$ & $0.03$ & $0.000021$ & $327$ & \textcolor{red}{$>2 \times 10^6$} & $24,247.41$ & \textcolor{red}{$0.018320$} & $1,690$ & $1,945,435$ & \textcolor{red}{$>43,000.00$} & \textcolor{red}{$0.034164$}  \\

				\directm & $31$ & $501$ & $0.08$ & $0.000021$ & $10,174$ & \textcolor{red}{$>2 \times 10^6$} & $692.12$ & \textcolor{red}{$0.001229$} & $3,254$ & \textcolor{red}{$>2 \times 10^6$} & $233.60$ & \textcolor{red}{$0.007284$}  \\
				\directz & $121$ & $1,329$ & $3.38$ & $0.000021$ & $103,675$ & \textcolor{red}{$>2 \times 10^6$} & $2,761.66$ & \textcolor{red}{$0.004189$} & $122,325$ & \textcolor{red}{$>2 \times 10^6$} & $4,207.92$ & \textcolor{red}{$0.097183$}  \\
				\directrev$^{*}$ & $1$ & $103$ & $0.13$ & $9.29 \times 10^{-14}$ & $1$ & $343$ & $0.12$ & $1.68 \times 10^{-12}$ & $2$ & $2,218$ & $0.52$ & $8.37 \times 10^{-9}$  \\
				\directav & $31$ & $501$ & $0.06$ & $0.000021$ & $11,505$ & \textcolor{red}{$>2 \times 10^6$} & $799.78$ & \textcolor{red}{$0.000691$} & $3,585$ & \textcolor{red}{$>2 \times 10^6$} & $338.38$ & \textcolor{red}{$0.007331$}  \\
				\dirmin$^{*}$ & $3$ & $415$ & $0.18$ & $4.36 \times 10^{-13}$ & $1$ & $395$ & $0.07$ & $1.67 \times 10^{-12}$ & $2$ & $2,567$ & $0.33$ & $8.37 \times 10^{-12}$  \\
				\plor 	  & $35$  & $305$ 	 & $0.04$  & $0.000021$ & $203,119$ & \textcolor{red}{$>2 \times 10^6$} & $6,746.14$ & \textcolor{red}{$0.018320$} & $93,457$ & \textcolor{red}{$>2 \times 10^6$} & $2,287.23$ & \textcolor{red}{$0.009475$}  \\
				\glbsolve & $31$  & $517$    & $0.17$  & $0.000021$ & $2,805$ & $1,978,935$     & $190.77$   & $0.000097$ & $1,415$  & \textcolor{red}{$>2 \times 10^6$} & $134.80$   & \textcolor{red}{$0.021190$}  \\
				\symdirec & $27$ & $369$ & $0.22$ & $0.000021$ & $13,177$ & \textcolor{red}{$>2 \times 10^6$} & $859.96$ & \textcolor{red}{$0.283888$} & $9,085$ & \textcolor{red}{$>2 \times 10^6$} & $1,572.94$ & \textcolor{red}{$0.144776$}  \\
				\symdirect & $25$ & $331$ & $0.23$ & $0.000021$ & $13,180$ & \textcolor{red}{$>2 \times 10^6$} & $901.01$ & \textcolor{red}{$0.283888$} & $9,086$ & \textcolor{red}{$>2 \times 10^6$} & $1,661.46$ & \textcolor{red}{$0.144776$}  \\
				\directmr & $106$ & $1,615$ & $0.48$ & $0.000021$ & $15,563$ & \textcolor{red}{$>2 \times 10^6$} & $558.19$ & \textcolor{red}{$0.000289$} & $8,796$ & \textcolor{red}{$>2 \times 10^6$} & $391.51$ & \textcolor{red}{$0.021800$}  \\
				\directmr$_{075}$ & $106$ & $1,493$ & $0.57$ & $0.000021$ & $17,801$ & \textcolor{red}{$>2 \times 10^6$} & $649.95$ & \textcolor{red}{$0.000665$} & $15,189$ & \textcolor{red}{$>2 \times 10^6$} & $604.55$ & \textcolor{red}{$0.021733$}  \\
				\birect   & $58$  & $616$    & $0.43$  & $0.000094$ & $8,729$ & $471,930$ & $438.66$ & $0.000096$ & $36,981$ & \textcolor{red}{$>2 \times 10^6$} & $5,185.32$ & \textcolor{red}{$0.004731$}  \\
				\disimplcgb & $240$ & $4,504$ & $3.56$ & $0.000084$ & $11,407$ & $871,606$ & $11,908.40$ & $0.000099$ & $124,966$ & \textcolor{red}{$>2 \times 10^6$} & $18,480.38$ & \textcolor{red}{$0.022092$}  \\
				\disimplvgb & $223$ & $2,563$ & $2.61$ & $0.000093$ & $13,021$ & $204,286$ & \textcolor{red}{$>43,000.00$} & \textcolor{red}{$0.004809$} & $14$ & $5,676$ & $18,021.63$ & \textcolor{red}{$0.470434^{\beta}$}  \\

				\birectgb & $57$ & $616$ & $0.29$ & $0.000095$ & $12,546$ & $356,724$ & $293.07$ & $0.000097$ & $45,716$ & \textcolor{red}{$>2 \times 10^6$} & $9,055.27$ & \textcolor{red}{$0.003564$}  \\
				\birmin$^{*}$   & $5$   & $101$    & $0.03$  & $1.04 \times 10^{-13}$ & $15$ & $316$ & $0.03$ & $3.82 \times 10^{-7}$ & $85,458$ & $722,243$ & $5,438.98$ & $8.40 \times 10^{-9}$  \\
				\gbdirect & $31$ & $501$ & $0.25$ & $0.000021$ & $16,421$ & \textcolor{red}{$>2 \times 10^6$} & $1,384.93$ & \textcolor{red}{$0.000353$} & $9,656$ & \textcolor{red}{$>2 \times 10^6$} & $489.68$ & \textcolor{red}{$0.009475$}  \\
				\disimplc & $223$ & $4,856$ & $5.88$ & $0.000084$ & $7,417$ & $664,424$ & $36,619.83$ & $0.000099$ & $6,620$ & $1,296,672$ & \textcolor{red}{$>43,000.00$} & \textcolor{red}{$0.011085$}  \\
				\disimplv & $150$ & $2,344$ & $2.18$ & $0.000091$ & $7,506$ & $188,518$ & \textcolor{red}{$>43,000.00$} & \textcolor{red}{$0.004860$} & $6$ & $5,447$ & $13,480.53$ & \textcolor{red}{$0.474346^{\beta}$} \\
				\adc 	& $549$ & $2,459$ & $1.81$ & $0.000089$ & $126,297$ & $398,856$ & \textcolor{red}{$>43,000.00$} & \textcolor{red}{$0.004755$} & $97,369$ & $213,600$ & \textcolor{red}{$>43,000.00$} & \textcolor{red}{$0.096970$}  \\
				\directa  & $31$  & $5,061$  & $0.41$  & $0.000087$ & $311$   & \textcolor{red}{$>2 \times 10^6$} & $68.24$   & \textcolor{red}{$0.005848$} & $213$    & \textcolor{red}{$>2 \times 10^6$} & $64.41$   & \textcolor{red}{$0.099088$}  \\
				\directg & $36$ & $863$ & $0.14$ & $0.000088$ & $2,604$ & \textcolor{red}{$>2 \times 10^6$} & $124.51$ & \textcolor{red}{$0.000289$} & $1,974$ & \textcolor{red}{$>2 \times 10^6$} & $110.58$ & \textcolor{red}{$0.075381$}  \\
				\directl & $71$ & $3,087$ & $0.32$ & $0.000085$ & $1,576$ & $427,475$ & $43.53$ & $0.000098$ & $3,294$ & \textcolor{red}{$>2 \times 10^6$} & $333.59$ & \textcolor{red}{$0.000158$}  \\
				\directgl & $27$  & $1,617$  & $0.29$  & $0.000021$ & $184$   & $62,965$ 		& $3.32$ 	 & $0.000099$ & $652$ 	 & $543,207$ 	   & $36.50$ 	& $0.000086$  \\
				\glbolve & $32$ & $501$ & $0.05$ & $0.000021$ & $2,768$ & $1,993,523$ & $3,282.31$ & $0.000098$ & $1,209$ & \textcolor{red}{$>2 \times 10^6$} & $1,566.85$ & \textcolor{red}{$0.031797$} \\
				\glccluster$^{*}$ & $1$ & $5,688$ & $0.51$ & $9.77 \times 10^{-10}$ & $1$ & $12,065$ & $1.53$ & $2.11 \times 10^{-9}$ & $2$ & $47,903$ & $259.97$ & $0.000001$  \\
				\midrule
				Success rate $(\%)$ & \multicolumn{4}{c|}{$100.00$} & \multicolumn{4}{c|}{$42.86$} & \multicolumn{4}{c|}{$17.86$}	 \\
				\bottomrule
				\multicolumn{13}{l}{$^{\beta}$ -- algorithm crash, lack of memory}\\
				\multicolumn{13}{l}{$^{*}$ -- a hybrid version of the algorithm, enriched with the local search subroutine}\\
				\multicolumn{5}{l}{N/A -- not available}
			\end{tabular}
			}
	\label{tab:sinusoids}
	\end{table}

		\begin{table}[htb]
			\caption{Performance of the \direct-type algorithms from \toolbox{} and \tomlab{} on a non-linear regression design problem}
			\resizebox{\textwidth}{!}{
			\begin{tabular}[tb]{@{\extracolsep{\fill}}l|rrrr|rrrr|rrrr|}
				\toprule
				Algorithm & Iter. & $f_{eval}$ & Time(s) & $f_{min}$ & Iter. & $f_{eval}$ & Time(s) & $f_{min}$ & Iter. & $f_{eval}$ & Time(s) & $f_{min}$ \\
				\midrule
					\multicolumn{1}{c}{} & \multicolumn{4}{c}{$\varsigma=1, T=100$}  & \multicolumn{4}{c}{$\varsigma=2, T=100$}  & \multicolumn{4}{c}{$\varsigma =3, T=100$} \\
				\midrule
				\direct   	& $89$ & $1,923$ & $0.15$ & $0.000025$ & $9,191$ & \textcolor{red}{$>2 \times 10^6$} & $367.06$ & \textcolor{red}{$0.001696$} & $3,722$ & \textcolor{red}{$>2 \times 10^6$} & $309.91$ & \textcolor{red}{$0.034532$} \\
				\directrest & $88$ & $1,839$ & $0.15$ & $0.000026$ & $379$ & $916,449$ & \textcolor{red}{$>43,000.00$} &\textcolor{red}{ $0.018320$} & $1,153$ & \textcolor{red}{$>2 \times 10^6$} & $4,625.37$ & \textcolor{red}{$0.059241$} \\
				\directm 	& $88$ & $1,807$ & $0.28$ & $0.000026$ & $9,596$ & \textcolor{red}{$>2 \times 10^6$} & $627.22$ & \textcolor{red}{$0.001697$} & $3,736$ & \textcolor{red}{$>2 \times 10^6$} & $330.63$ & \textcolor{red}{$0.034538$} \\
				\directz    & $167$ & $2,037$ & $4.81$ & $0.000026$ & $106,353$ & \textcolor{red}{$>2 \times 10^6$} & $2,543.40$ & \textcolor{red}{$0.006783$} & $125,804$ & \textcolor{red}{$>2 \times 10^6$} & $4,253.70$ & \textcolor{red}{$0.114946$} \\
				\directrev$^*$  & $1$ & $95$ & $0.13$ & $9.30 \times 10^{-14}$ & $1$ & $318$ & $0.09$ & $1.13 \times 10^{-12}$ & $7$ & $3,677$ & $0.54$ & $5.65 \times 10^{-11}$ \\
				\directav  	& $88$ & $1,799$ & $0.22$ & $0.000026$ & $11,545$ & \textcolor{red}{$>2 \times 10^6$} & $821.49$ & \textcolor{red}{$0.001255$} & $3,809$ & \textcolor{red}{$>2 \times 10^6$} & $501.22$ & \textcolor{red}{$0.034546$} \\
				\dirmin$^{*}$		& $3$ & $416$ & $0.19$ & $1.15 \times 10^{-13}$ & $1$ & $415$ & $0.08$ & $5.83 \times 10^{-10}$ & $3$ & $3,653$ & $0.53$ & $6.12 \times 10^{-11}$ \\
				\plor 	  	& $6,811$ & $45,423$ & $10.90$ & $0.000025$ & $203,159$ & \textcolor{red}{$>2 \times 10^6$} & $8,147.64$ & \textcolor{red}{$0.020558$} & $91,574$ & \textcolor{red}{$>2 \times 10^6$} & $2,395.47$ & \textcolor{red}{$0.041690$} \\
				\glbsolve 	& $37$ & $669$ & $0.07$ & $0.000026$ & $2,602$ & $1,713,841$ & $195.49$ & $0.000096$ & $1,165$ & \textcolor{red}{$>2 \times 10^6$} & $169.23$ & \textcolor{red}{$0.041937$} \\
				\symdirec	& $62$ & $1,239$ & $0.29$ & $0.000026$ & $13,098$ & \textcolor{red}{$>2 \times 10^6$} & $873.50$ & \textcolor{red}{$0.298576$} & $9,235$ & \textcolor{red}{$>2 \times 10^6$} & $1,579.90$ & \textcolor{red}{$0.161607$} \\
				\symdirect	& $59$ & $1,149$ & $0.26$ & $0.000026$ & $13,098$ & \textcolor{red}{$>2 \times 10^6$} & $894.10$ & \textcolor{red}{$0.298576$} & $9,225$ & \textcolor{red}{$>2 \times 10^6$} & $1,633.59$ & \textcolor{red}{$0.161607$} \\
				\directmr	& $124$ & $1,991$ & $0.60$ & $0.000026$ & $15,324$ & \textcolor{red}{$>2 \times 10^6$} & $547.68$ & \textcolor{red}{$0.001146$} & $9,328$ & \textcolor{red}{$>2 \times 10^6$} & $436.52$ & \textcolor{red}{$0.043663$} \\
				\directmr$_{075}$ & $124$ & $1,857$ & $0.56$ & $0.000026$ & $17,627$ & \textcolor{red}{$>2 \times 10^6$} & $680.61$ & \textcolor{red}{$0.001247$} & $15,660$ & \textcolor{red}{$>2 \times 10^6$} & $670.92$ & \textcolor{red}{$0.034531$} \\
				\birect   	& $86$ & $1,042$ & $0.22$ & $0.000096$ & $13,028$ & $675,332$ & $238.64$ & $0.000097$ & $36,473$ & \textcolor{red}{$>2 \times 10^6$} & $6,874.16$ & \textcolor{red}{$0.019857$} \\
				\disimplcgb & $226$ & $4,098$ & $2.90$ & $0.000098$ & $27,508$ & \textcolor{red}{$>2 \times 10^6$} & $5,240.49$ & \textcolor{red}{$0.003371$} & $121,510$ & \textcolor{red}{$>2 \times 10^6$} & $21,508.76$ & \textcolor{red}{$0.044617$} \\
				\disimplvgb & $227$ & $2,353$ & $2.04$ & $0.000004$ & $11,838$ & $168,651$ & \textcolor{red}{$>43,000.00$} & \textcolor{red}{$0.008066$} & $14$ & $5,676$ & $19,133.02$ & \textcolor{red}{$0.488079^{\beta}$} \\
				\birectgb	& $85$ & $1,042$ & $0.26$ & $0.000097$ & $15,029$ & $414,268$ & $370.36$ & $0.000091$ & $48,312$ & \textcolor{red}{$>2 \times 10^6$} & $10,437.74$ & \textcolor{red}{$0.014167$} \\
				\birmin$^{*}$   	& $5$ & $105$ & $0.02$ & $1.15 \times 10^{-13}$ & $15$ & $292$ & $0.06$ & $2.23 \times 10^{-7}$ & $89,360$ & $730,887$ & $6,871.27$ & $2.30 \times 10^{-9}$ \\
				\gbdirect	& $86$ & $1,691$ & $0.44$ & $0.000026$ & $16,218$ & \textcolor{red}{$>2 \times 10^6$} & $1,316.70$ & \textcolor{red}{$0.000712$} & $13,339$ & \textcolor{red}{$>2 \times 10^6$} & $599.68$ & \textcolor{red}{$0.041748$} \\
				\disimplc	& $214$ & $4,556$ & $5.31$ & $0.000098$ & $7,518$ & $650,936$ & $32,115.01$ & $0.000088$ & $7,148$ & $1,309,272$ & \textcolor{red}{$>43,000.00$} & \textcolor{red}{$0.014357$} \\
				\disimplv 	& $179$ & $3,024$ & $2.76$ & $0.000087$ & $10,124$ & $244,502$ & \textcolor{red}{$>43,000.00$} & \textcolor{red}{$0.008066$} & $6$ & $5,447$ & $11,406.93$ & \textcolor{red}{$0.488079^{\beta}$} \\
				\adc		& $644$ & $2,848$ & $2.46$ & $0.000088$ & $127,872$ & $355,695$ & \textcolor{red}{$>43,000.00$} & \textcolor{red}{$0.007043$} & $73,097$ & $190,302$ & \textcolor{red}{$>43,000.00$} & \textcolor{red}{$0.115773$} \\
				\directa  	& $36$ & $6,921$ & $0.52$ & $0.000091$ & $311$ & \textcolor{red}{$>2 \times 10^6$} & $91.84$ & \textcolor{red}{$0.006979$} & $213$ & \textcolor{red}{$>2 \times 10^6$} & $98.15$ & \textcolor{red}{$0.117012$} \\
				\directg	& $40$ & $1,045$ & $0.11$ & $0.000045$ & $2,361$ & \textcolor{red}{$>2 \times 10^6$} & $151.77$ & \textcolor{red}{$0.000786$} & $2,031$ & \textcolor{red}{$>2 \times 10^6$} & $147.34$ & \textcolor{red}{$0.094152$} \\
				\directl	& $67$ & $2,833$ & $0.30$ & $0.000093$ & $689$ & $294,415$ & $23.61$ & $0.000099$ & $3,617$ & \textcolor{red}{$>2 \times 10^6$} & $382.54$ & \textcolor{red}{$0.000248$} \\
				\directgl 	& $29$ & $1,829$ & $0.18$ & $0.000025$ & $188$ & $65,645$ & $3.93$ & $0.000097$ & $251$ & $158,989$ & $9.99$ & $0.000099$ \\
				\glbolve 	& $89$ & $1,997$ & $0.21$ & $0.000026$ & $2,495$ & $1,666,955$ & $2,160.47$ & $0.000097$ & $1,162$ & \textcolor{red}{$>2 \times 10^6$} & $1,557.18$ & \textcolor{red}{$0.041937$}  \\
				\glccluster$^{*}$ & $1$ & $6,092$ & $0.57$ & $1.28 \times 10^{-9}$ & $1$ & $12,034$ & $1.18$ & $4.47 \times 10^{-10}$ & $1$ & $18,008$ & $4.91$ & $0.000008$  \\
				\midrule
				Success rate $(\%)$ & \multicolumn{4}{c|}{$100.00$} & \multicolumn{4}{c|}{$39.28$} & \multicolumn{4}{c|}{$17.86$}	 \\
				\bottomrule
				\multicolumn{13}{l}{$^{\beta}$ -- algorithm crash, lack of memory}\\
				\multicolumn{13}{l}{$^{*}$ -- a hybrid version of the algorithm, enriched with the local search subroutine}\\
				\multicolumn{5}{l}{N/A -- not available}
			\end{tabular}
		}
		\label{tab:sinusoid}
	\end{table}
	Our experiments have used three different values for $\varsigma = 1, 2,$ and $3$ (correspond to $3, 6,$ and $9$-dimensional problems) and two different values, $T = 10$ and $T = 100$ for each dimension $n$ as was done in \cite{Paulavicius2019:eswa}.

	The obtained results are summarized in \cref{tab:sinusoids,tab:sinusoid}.
	Solving the lowest dimension $(n=3)$ cases (corresponding to $\varsigma=1, T=10$ and $\varsigma=1, T=100$) all algorithms located solution correctly (the success rate is $100 \%$).
	Among the hybridized methods, the \birmin{} and \directrev{} proved the most effective.
	Among the traditional algorithms, the \plor{} algorithm is the most efficient in solving the first problem (with $\varsigma=1, T=10$ parameters).
	However, by increasing the number of samples $T$, \plor{} performed worst among all \direct-type algorithms.
	The most efficient algorithm for the case with $\varsigma=1$ and $T=100$ was \glbsolve.

	For the higher dimensionality case $(n=6)$, more than half of the algorithms failed to find the correct solution.
	The success rates for these two cases are  $42.86 \%$ and $39.28 \%$, respectively.
	The \directgl{} algorithm has shown a significant advantage among the traditional \direct-type algorithms.
	The best two performing hybridized methods were the \birmin{} and \directrev.

	Finally, most \direct-type methods have faced significant challenges in solving two variants of the highest dimensional $(n=9)$ case.
	The success rates on these two cases are only $17.86 \%$ and $17.86 \%$, respectively.
  While hybridized methods had no significant difficulties, among the traditional, only the \directgl{} algorithm solved both variants.

	\section{Conclusion}
	\label{sec:conclusions}

	This paper has introduced a new open-source \direct-type \matlab{} toolbox (\toolbox) for derivative-free global optimization.
	The new toolbox combines various state-of-the-art \direct-type algorithms for the global solution of box-constrained, generally-constrained, and optimization problems with hidden constraints.
	All algorithms were implemented using two different data structures: static and dynamic.
	Additionally, several parallel schemes were adopted to promising algorithms.
	Furthermore, an online test library \directlib, containing $119$ global optimization test and engineering problems, has been presented.

	The performance of various algorithms within \toolbox{} has been investigated via a detailed numerical study using the test problems from \directlib.
	A further 11 examples of using the \toolbox{} for engineering design optimization have been investigated.
	The results demonstrate the promising performance of \toolbox{} in tackling these challenging problems.
	We also gave guidance on which algorithms to use for specific optimization problems.

	Motivated by the promising performance, we plan to extend this work to facilitate the broader adoption of \toolbox{}.
	We plan to include newly appearing promising \direct-type algorithms within this toolbox continuously.
	Another direction is extending the developed algorithms using a hybrid CPU-GPU scheme.
	Finally, we will consider advanced data structures for better organization and reduced communication overhead.

	\section*{Source code statement}
	All implemented \direct-type algorithms (\toolbox{} toolbox) are available at the GitHub repository: \url{https://github.com/blockchain-group/DIRECTGO} and can be used under the MIT license.
  We welcome contributions and corrections to it.

	\section*{Data statement}
	\texttt{DIRECTGOLib} - \direct{} \textbf{G}lobal \textbf{O}ptimization test problems \textbf{Lib}rary is designed as a \textit{continuously-growing} open-source GitHub repository (\url{https://github.com/blockchain-group/DIRECTGOLib}) to which anyone can easily contribute.
   Therefore, the most recent version is slightly different from the one used in these studies.
   The exact data underlying this article (\directlib) can be accessed either on GitHub or at Zenodo (connected with GitHub):
   \begin{itemize}
     \item at GitHub: \url{https://github.com/blockchain-group/DIRECTGOLib/tree/v1.0},
     \item at Zenodo: \url{https://doi.org/10.5281/zenodo.6491863},
   \end{itemize}
    and used under the MIT license.

	\bibliographystyle{ACM-Reference-Format}
	\bibliography{main}

	\appendix

	\section{\directlib{} library}
	\label{sec:directlib}

	A summary of all optimization problems in \directlib~\cite{DIRECTGOLibv1.0,DIRECTGOLibZenodov1.0} and their properties are given in~\cref{tab:resultslow}.
	The first column denotes the problem type, and the second contains the problem name.
	The third column contains the source of the problem.
	The fourth column specifies the dimension ($n$) of the used problems.
  When problems are of various dimensionality, all considered dimensions are listed.
  The fifth through eighth columns ($\mathbf{g}$,  $\mathbf{h}$, $a$, $D$) specify the number of inequality ($\mathbf{g}$) and equality ($\mathbf{h}$) constraints, the number of active constraints ($a$), and two versions of optimization domains, default ($D$) and perturbed ($\tilde{D}$), respectively.
  The default domains are taken from the literature.
	However, whenever the global minimum point lay for at least one algorithm from \toolbox{} at the initial sampling point, the default feasible region ($D$) was perturbed and used in the experimental analysis.
  Here the sign ``-'' means that $\tilde{D}$ is the same as $D$.
	Finally, the last column contains the best-known optimal solution value ($f^*$).

	{\scriptsize
		\begin{longtable}[tb]{@{\extracolsep{\fill}}llrp{0.02\textwidth}rrrp{0.3\textwidth}rr}
			\label{tab:resultslow}\\
			\caption{Key characteristics of the \directlib~\cite{DIRECTGOLibv1.0,DIRECTGOLibZenodov1.0} test problems for global optimization}\\
			\toprule
			\multicolumn{2}{c}{Problem} & \multirow{2}{*}{Source} &	\multicolumn{6}{c}{Problem properties}	&  \multirow{2}{*}{$ f^* $} 	\\
			\cmidrule{0-1} \cmidrule{4-9}
			Type & Name &  & $n$ & $ \mathbf{g} $ & $ \mathbf{h} $ & $ a $ & $ D $  & $ \tilde{D} $ & \\
			\midrule
			\endfirsthead
			\multicolumn{10}{c}{{\bfseries Table \thetable\ Continued from previous page}} \\
			\toprule
			\multicolumn{2}{c}{Problem} & \multirow{2}{*}{Source} &	\multicolumn{6}{c}{Problem properties}	&  \multirow{2}{*}{$ f^* $} 	\\
			\cmidrule{0-1} \cmidrule{4-9}
			Type & Name &  & $n$ & $ \mathbf{g} $ & $ \mathbf{h} $ & $ a $ & $D$ (Default) & $D$ (Perturbed) & \\
			\midrule
			\endhead
			\midrule
			\multicolumn{10}{c}{{Continued on next page}}\\
			\midrule
			\endfoot
			\multicolumn{10}{c}{{Concluded}}\\
			\bottomrule
			\endlastfoot

			\textbf{BC} & \textit{Ackley} & \cite{Hedar2005,Derek2013} & $2,5,10$ & $0$ & $0$ & $0$ & $[-15, 30]^n$ & $[-15, 35]^n$ & $0.0000$ \\
			& \textit{Alpine} & \cite{Clerc1999} & $5,10,15$ & $0$ & $0$ & $0$ & $[0, 10]^n$ & $-$ & $-2.8081^n$ \\
			& \textit{Beale} & \cite{Hedar2005,Derek2013} & $2$ & $0$ & $0$ & $0$ & $[-4.5, 4.5]^n$ & $-$ & $0.0000$ \\
			& \textit{Bohachevsky$1$} & \cite{Hedar2005,Derek2013} & $2$ & $0$ & $0$ & $0$ & $[-100, 100]^n$ & $[-100, 110]^n$ & $0.0000$ \\
			& \textit{Bohachevsky$2$} & \cite{Hedar2005,Derek2013} & $2$ & $0$ & $0$ & $0$ & $[-100, 100]^n$ & $[-100, 110]^n$ & $0.0000$ \\
			& \textit{Bohachevsky$3$} & \cite{Hedar2005,Derek2013} & $2$ & $0$ & $0$ & $0$ & $[-100, 100]^n$ & $[-100, 110]^n$ & $0.0000$ \\
			& \textit{Booth} & \cite{Hedar2005,Derek2013} & $2$ & $0$ & $0$ & $0$ & $[-10, 10]^n$ & $-$ & $0.0000$\\
			& \textit{Branin} & \cite{Hedar2005,Dixon1978,Derek2013} & $2$ & $0$ & $0$ & $0$ & $[-5, 10] \times [0, 15]$ & $-$ & $0.3978$\\
			& \textit{Bukin6} & \cite{Derek2013} & $2$ & $0$ & $0$ & $0$ & $[-15, 5] \times [-3,3]$ & $-$ & $0.0000$ \\
			& \textit{Colville} & \cite{Hedar2005,Derek2013} & $4$ & $0$ & $0$ & $0$ & $[-10, 10]^n$ & $-$ & $0.0000$ \\
			& \textit{Cross\_in\_Tray} & \cite{Derek2013} & $2$ & $0$ & $0$ & $0$ & $[-10, 10]^n$ & $-$ & $-2.0626$ \\
			& \textit{Csendes} & \cite{Clerc1999} & $5,10,15$ & $0$ & $0$ & $0$ & $[-10, 10]^n$ & $[-10, 21]^n$ & $0.0000$ \\
			& \textit{Dixon\_and\_Price} & \cite{Hedar2005,Derek2013} & $2,5,10$ & $0$ & $0$ & $0$ & $[-10, -10]^n$ & $-$ & $0.0000$ \\
			& \textit{Drop\_wave} & \cite{Derek2013} & $2$ & $0$ & $0$ & $0$ & $[-5.12, -5.12]^n$ & $[-5.12, -6.12]^n$ & $-1.0000$ \\
			& \textit{Easom} & \cite{Hedar2005,Derek2013} & $2$ & $0$ & $0$ & $0$ & $[-100, 100]^n$ & $-$ & $-1.0000$\\
			& \textit{Eggholder} & \cite{Derek2013} & $2$ & $0$ & $0$ & $0$ & $[-512, 512]^n$ & $-$ & $-959.6406$ \\
			& \textit{Goldstein\_and\_Price} & \cite{Hedar2005,Dixon1978,Derek2013} & $2$ & $0$ & $0$ & $0$ & $[-2, 2]^n$ & $-$ & $3.0000$ \\
			& \textit{Griewank} & \cite{Hedar2005,Derek2013} & $5,10,15$ & $0$ & $0$ & $0$ & $[-600, 600]^n$ & $[-600, 700]^n$ & $0.0000$ \\
			& \textit{Hartman$3$} & \cite{Hedar2005,Derek2013} & $3$ & $0$ & $0$ & $0$ & $[0, 1]^n$ & $-$ & $-3.8627$\\
			& \textit{Hartman$6$} & \cite{Hedar2005,Derek2013} & $6$ & $0$ & $0$ & $0$ & $[0, 1]^n$ & $-$ & $-3.3223$\\
			& \textit{Holder\_Table} & \cite{Derek2013} & $2$ & $0$ & $0$ & $0$ & $[-10, 10]^n$ & $-$ & $-19.2085$ \\
			& \textit{Hump} & \cite{Hedar2005,Derek2013} & $2$ & $0$ & $0$ & $0$ & $[-5, 5]^n$ & $-$ & $-1.0316$\\
			& \textit{Langermann} & \cite{Derek2013} & $2$ & $0$ & $0$ & $0$ & $[0, 10]^n$ & $-$ & $-4.1558$\\
			& \textit{Levy} & \cite{Hedar2005,Derek2013} & $5,10,15$ & $0$ & $0$ & $0$ & $[-5, 5]^n$ & $-$ & $0.0000$ \\
			& \textit{Matyas} & \cite{Hedar2005,Derek2013} & $2$ & $0$ & $0$ & $0$ & $[-10, 10]^n$ & $[-10, 15]^n$ & $0.0000$ \\
			& \textit{McCormick} & \cite{Derek2013} & $2$ & $0$ & $0$ & $0$ & $[-1.5, 4] \times [-3,4]$ & $-$ & $-1.9132$\\
			& \textit{Michalewicz} & \cite{Hedar2005,Derek2013} & $2$ & $0$ & $0$ & $0$ & $[0, \pi]^n$ & $-$ & $-1.8013$ \\
			& \textit{Michalewicz} & \cite{Hedar2005,Derek2013} & $5$ & $0$ & $0$ & $0$ & $[0, \pi]^n$ & $-$ & $-4.6876$\\
			& \textit{Michalewicz} & \cite{Hedar2005,Derek2013} & $10$ & $0$ & $0$ & $0$ & $[0, \pi]^n$ & $-$ & $-9.6601$\\
			& \textit{Perm} & \cite{Hedar2005,Derek2013} & $8$ & $0$ & $0$ & $0$ & $[-i, i]^n$ & $-$ & $0.0000$ \\
			& \textit{Permdb} & \cite{Hedar2005,Derek2013} & $5$ & $0$ & $0$ & $0$ & $[-i, i]^n$ & $-$ & $0.0000$ \\
			& \textit{Powell} & \cite{Hedar2005,Derek2013} & $4$ & $0$ & $0$ & $0$ & $[-4, 4]^n$ & $[-4, 5]^n$ & $0.0000$ \\
			& \textit{Power\_Sum} & \cite{Hedar2005,Derek2013} & $4$ & $0$ & $0$ & $0$ & $[0, 4]^n$ & $-$ & $0.0000$ \\
			& \textit{Qing} & \cite{Clerc1999} & $5,10,15$ & $0$ & $0$ & $0$ & $[-500, 500]^n$ & $-$ & $0.0000$ \\
			& \textit{Rastrigin} & \cite{Hedar2005,Derek2013} & $2,5,10$ & $0$ & $0$ & $0$ & $[-5.12, 5.12]^n$ & $[-6.12, 5.12]^n$ & $0.0000$ \\
			& \textit{Rosenbrock} & \cite{Hedar2005,Dixon1978,Derek2013} & $5,10,15$ & $0$ & $0$ & $0$ & $[-5, 10]^n$ & $-$ & $0.0000$ \\
			& \textit{Rotated\_H\_Ellip} & \cite{Derek2013} & $5,10,15$ & $0$ & $0$ & $0$ & $[-65.536, 65.536]^n$ & $-65.536, 66.536]^n$ & $0.0000$ \\
			& \textit{Schwefel} & \cite{Hedar2005,Derek2013} & $2,5,10$ & $0$ & $0$ & $0$ & $[-500, 500]^n$ & $-$ & $0.0000$ \\
			& \textit{Shekel$5$} & \cite{Hedar2005,Derek2013} & $4$ & $0$ & $0$ & $0$ & $[0, 10]^n$ & $-$ & $-10.1531$ \\
			& \textit{Shekel$7$} & \cite{Hedar2005,Derek2013} & $4$ & $0$ & $0$ & $0$ & $[0, 10]^n$ & $-$ & $-10.4029$ \\
			& \textit{Shekel$10$} & \cite{Hedar2005,Derek2013} & $4$ & $0$ & $0$ & $0$ & $[0, 10]^n$ & $-$ & $-10.5364$ \\
			& \textit{Shubert} & \cite{Hedar2005,Derek2013} & $2$ & $0$ & $0$ & $0$ & $[-10, 10]^n$ & $-$ & $-186.7309$ \\
			& \textit{Sphere} & \cite{Hedar2005,Derek2013} & $5,10,15$ & $0$ & $0$ & $0$ & $[-5, 5]^n$ & $[-5.12, 6.12]^n$ & $0.0000$ \\
			& \textit{Styblinski\_Tang} & \cite{Clerc1999} & $5,10,15$ & $0$ & $0$ & $0$ & $[-5, 5]^n$ & $-$ & $-39.1661n$ \\
			& \textit{Sum\_of\_Powers} & \cite{Derek2013} & $5,10,15$ & $0$ & $0$ & $0$ & $[-1, 1]^n$ & $[-1, 2.5]^n$ & $0.0000$ \\
			& \textit{Sum\_Square} & \cite{Derek2013} & $5,10,15$ & $0$ & $0$ & $0$ & $[-10, 10]^n$ & $[-10, 15]^n$ & $0.0000$ \\
			& \textit{Trid$6$} & \cite{Hedar2005,Derek2013} & $6$ & $0$ & $0$ & $0$ & $[-36, 36]^n$ & $-$ & $-50.0000$ \\
			& \textit{Trid$10$} & \cite{Hedar2005,Derek2013} & $10$ & $0$ & $0$ & $0$ & $[-100, 100]^n$ & $-$ & $-210.0000$ \\
			& \textit{Zakharov} & \cite{Hedar2005,Derek2013} & $2,5,10$ & $0$ & $0$ & $0$ & $[-5, 10]^n$ & $[-5, 11]^n$ & $0.0000$ \\
			\midrule
			\textbf{LC} & \textit{Bunnag$1$} & \cite{Vaz2009} & $4$ & $1$ & $0$ & $1$ & $[0, 3]^n$ & $-$ & $0.1111$\\
			& \textit{Bunnag$2$} & \cite{Vaz2009} & $4$ & $2$ & $0$ & $2$ & $[0, 4]^n$ & $-$ & $-6.4052$ \\
			& \textit{Bunnag$3$} & \cite{Vaz2009} & $5$ & $3$ & $0$ & $1$ & $[0, 3] \times [0, 2] \times [0, 4] \times [0, 4] \times [0, 2]$ & $-$ & $-16.3692$ \\
			& \textit{Bunnag$4$} & \cite{Vaz2009} & $6$ & $2$ & $0$ & $1$ & $[0, 1]^5 \times [0, 20]$ & $-$ & $-213.0470$ \\
			& \textit{Bunnag$5$} & \cite{Vaz2009} & $6$ & $5$ & $0$ & $1$ & $[0, 2] \times [0, 8] \times [0, 2] \times [0, 1]^2 \times [0, 2]$ & $-$ & $-11.0000$ \\
			& \textit{Bunnag$6$} & \cite{Vaz2009} & $10$ & $11$ & $0$ & $3$ & $[0, 1]^n$ & $-$ & $-268.0146$ \\
			& \textit{Bunnag$7$} & \cite{Vaz2009} & $10$ & $5$ & $0$ & $0$ & $[0, 1]^n$ & $-$ & $-39.0000$ \\
			& \textit{G$01$} & \cite{Suganthan2005} & $13$ & $9$ & $0$ & $6$ & $[0, 10]^9 \times [0, 100]^3 \times [0, 10]$ & $-$ & $-15.0000$ \\
			& \textit{Genocop$9$} & \cite{Vaz2009} & $3$ & $5$ & $0$ & $2$ & $[0, 10]^n$ & $-$ & $-2.4714$ \\
			& \textit{Genocop$10$} & \cite{Vaz2009} & $4$ & $5$ & $0$ & $0$ & $[0, 3] \times [0, 10]^2 \times [0, 1]$ & $-$ & $-4.5280$ \\
			& \textit{Genocop$11$} & \cite{Vaz2009} & $6$ & $5$ & $0$ & $0$ & $[0, 5] \times [0, 8] \times [0, 5] \times [0, 1]^2 \times [0, 2]$ & $-$ & $-11.0000$ \\
			& \textit{Horst$1$} & \cite{Horst1995:book} & $2$ & $3$ & $0$ & $1$ & $[0, 3] \times [0, 2]$ & $-$ & $-1.0625$ \\
			& \textit{Horst$2$} & \cite{Horst1995:book} & $2$ & $3$ & $0$ & $2$ & $[0, 2.5] \times [0, 2]$ & $-$ & $-6.8995$ \\
			& \textit{Horst$3$} & \cite{Horst1995:book} & $2$ & $3$ & $0$ & $0$ & $[0, 1] \times [0, 1.5]$ & $-$ & $-0.4444$ \\
			& \textit{Horst$4$} & \cite{Horst1995:book} & $3$ & $4$ & $0$ & $2$ & $[0.5, 2] \times [0, 3] \times [0, 2.8]$ & $-$ & $-6.0858$ \\
			& \textit{Horst$5$} & \cite{Horst1995:book} & $3$ & $4$ & $0$ & $0$ & $[0, 1.2]^2 \times[0, 1.7]$ & $-$ & $-3.7220$ \\
			& \textit{Horst$6$} & \cite{Horst1995:book} & $3$ & $7$ & $0$ & $2$ & $[0, 6] \times [0, 5.0279] \times [0, 2.6]$ & $-$ & $-32.5784$ \\
			& \textit{Horst$7$} & \cite{Horst1995:book} & $3$ & $4$ & $0$ & $2$ & $[0, 6] \times [0, 3]^2$ & $-$ & $-52.8769$ \\
			& \textit{hs$021$} & \cite{Vaz2009} & $2$ & $1$ & $0$ & $1$ & $[2, 50] \times [-50, 10]$ & $-$ & $-99.9599$ \\
			& \textit{h$s021$mod} & \cite{Vaz2009} & $7$ & $3$ & $0$ & $1$ & $[2, 50] \times [-50, 50] \times [0, 50] \times [2, 10] \times [-10, 10] \times [-10, 0] \times [0, 10]$ & $-$ & $4.0400$ \\
			& \textit{hs$024$} & \cite{Vaz2009} & $2$ & $3$ & $0$ & $2$ & $[0, 5]^n$ & $-$ & $-1.0000$ \\
			& \textit{hs$035$} & \cite{Vaz2009} & $3$ & $1$ & $0$ & $1$ & $[0, 3]^n$ & $-$ & $0.1111$\\
			& \textit{hs$036$} & \cite{Vaz2009} & $3$ & $1$ & $0$ & $1$ & $[0, 20] \times [0, 11] \times [0, 15]$ & $-$ & $-3300.0000$ \\
			& \textit{hs$037$} & \cite{Vaz2009} & $3$ & $2$ & $0$ & $1$ & $[0, 42]^n$ & $-$ & $-3456.0000$ \\
			& \textit{hs$038$} & \cite{Vaz2009} & $4$ & $2$ & $0$ & $0$ & $[-10, 10]^n$ & $-$ & $0.0000$\\
			& \textit{hs$044$} & \cite{Vaz2009} & $4$ & $6$ & $0$ & $2$ & $[0, 5]^n$ & $-$ & $-15.0000$ \\
			& \textit{hs$076$} & \cite{Vaz2009} & $4$ & $3$ & $0$ & $1$ & $[0, 1] \times [0, 3] \times [0, 1]^2$ & $-$ & $-4.6818$ \\
			& \textit{P$9$} & \cite{Floudas1999book1} & $3$ & $9$ & $0$ & $2$ & $[10^{-5}, 3] \times [10^{-5}, 4]^2$ & $-$ & $-13.4020$\\
			& \textit{P$14$} & \cite{Floudas1999book1} & $3$ & $4$ & $0$ & $2$ & $[10^{-5}, 3] \times [10^{-5}, 4] \times [0, 1]$ & $-$ & $-4.51420$\\
			& \textit{s$224$} & \cite{Vaz2009} & $2$ & $4$ & $0$ & $1$ & $[0, 6] \times [0, 11]$ & $-$ & $-304.0000$ \\
			& \textit{s$231$} & \cite{Vaz2009} & $2$ & $2$ & $0$ & $0$ & $[-10, 10]^n$ & $-$ & $0.0000$\\
			& \textit{s$232$} & \cite{Vaz2009} & $2$ & $3$ & $0$ & $2$ & $[0, 100]^n$ & $-$ & $-1.0000$\\
			& \textit{s$250$} & \cite{Vaz2009} & $3$ & $2$ & $0$ & $1$ & $[0, 20] \times [0, 11] \times [0, 42]$ & $-$ & $-3300.0000$ \\
			& \textit{s$251$} & \cite{Vaz2009} & $3$ & $1$ & $0$ & $1$ & $[0, 42]^n$ & $-$ & $-3456.0000$ \\
			& \textit{zecevic$2$} & \cite{Vaz2009} & $3$ & $2$ & $0$ & $1$ & $[0, 10]^n$ & $-$ & $-4.1249$\\
			\midrule
			\textbf{GC} & \textit{circle} & \cite{Vaz2009} & $3$ & $10$ & $0$ & $3$ & $[0, 10]^n$ & $-$ & $4.5742$\\
			& \textit{G$02$} & \cite{Suganthan2005} & $20$ & $2$ & $0$ & $1$ & $[0, 10]^n$ & $-$ & $-0.8036$\\
			& \textit{G$04$} & \cite{Suganthan2005} & $5$ & $6$ & $0$ & $2$ & $[78, 102] \times [33, 45]\times [27, 45]^3$ & $-$ & $-30665.5386$ \\
			& \textit{G$06$} & \cite{Suganthan2005} & $2$ & $2$ & $0$ & $2$ & $[13, 100] \times [0, 100]$ & $-$ & $-6961.8138$ \\
			& \textit{G$07$} & \cite{Suganthan2005} & $10$ & $8$ & $0$ & $6$ & $[-10, 10]^n$ & $-$ & $24.3062$ \\
			& \textit{G$08$} & \cite{Suganthan2005} & $2$ & $2$ & $0$ & $0$ & $[0, 10]^n$ & $-$ & $-0.0958$\\
			& \textit{G$09$} & \cite{Suganthan2005} & $7$ & $4$ & $0$ & $2$ & $[-10, 10]^n$ & $-$ & $680.6300$ \\
			& \textit{G$10$} & \cite{Suganthan2005} & $8$ & $6$ & $0$ & $6$ & $[100, 10,000] \times [1,000, 10,000]^2 \times [10, 1,000]^5$ & $-$ & $7049.2480$ \\
			& \textit{G$12$} & \cite{Suganthan2005} & $3$ & $1$ & $0$ & $0$ & $[0.2, 10]^n$ & $-$ & $-1.0000$\\
			& \textit{G$16$} & \cite{Suganthan2005} & $5$ & $38$ & $0$ & $4$ & $[704.4148, 906.3855] \times [68.6, 288.88] \times [0, 134.75] \times [193, 287.0966] \times [25, 84.1988]$ & $-$ & $-1.9051$ \\
			& \textit{G$18$} & \cite{Suganthan2005} & $9$ & $13$ & $0$ & $6$ & $[0, 10]^n$ & $-$ & $-0.8660$\\
			& \textit{G$19$} & \cite{Suganthan2005} & $15$ & $5$ & $0$ & $0$ & $[0, 10]^n$ & $-$ & $32.6555$\\
			& \textit{G$24$} & \cite{Suganthan2005} & $2$ & $2$ & $0$ & $2$ & $[0, 3] \times [0, 4]$ & $-$ & $-5.5080$\\
			& \textit{Goldstein\_and\_PriceC} & \cite{Na2017} & $2$ & $2$ & $0$ & $1$ & $[-2, 2]^n$ & $-$ & $3.5389$\\
			& \textit{Gomez} & \cite{Birgin2010} & $2$ & $1$ & $0$ & $1$ & $[-1, 1]^n$ & $-$ & $-0.9711$\\
			& \textit{Himmelblau} & \cite{Cagnina2008} & $5$ & $5$ & $0$ & $2$ & $[78, 102] \times [33, 45] \times [27, 45]^3$ & $-$ & $-31025.5602$ \\
			& \textit{P$1$} & \cite{Floudas1999book1} & $5$ & $0$ & $3$ & $3$ & $[-5, 5]^n$ & $-$ & $0.0293$\\
			& \textit{P$2$a} & \cite{Floudas1999book1} & $5$ & $10$ & $0$ & $5$ & $[0, 500]^5$ & $-$ & $-400.0000$\\
			& \textit{P$2$b} & \cite{Floudas1999book1} & $5$ & $10$ & $0$ & $5$ & $[0, 500]^5$ & $-$ & $-600.0000$\\
			& \textit{P$2$c} & \cite{Floudas1999book1} & $5$ & $10$ & $0$ & $4$ & $[0, 500]^5$ & $-$ & $-750.0000$\\
			& \textit{P$2$d} & \cite{Floudas1999book1} & $5$ & $12$ & $0$ & $5$ & $[0, 100] \times [0, 200] \times [0, 100] \times [0, 200] \times [1, 3]$ & $-$ & $-400.0000$\\
			& \textit{P$3$a} & \cite{Floudas1999book1} & $6$ & $1$ & $4$ & $5$ & $[0, 1]^4 \times [10^{-5}, 16]^2$ & $-$ & $0.3888$\\
			& \textit{P$3$b} & \cite{Floudas1999book1} & $2$ & $1$ & $0$ & $1$ & $[10^{-5}, 16]^n$ & $-$ & $0.3888$\\
			& \textit{P$4$} & \cite{Floudas1999book1} & $2$ & $1$ & $0$ & $1$ & $[0, 6] \times [0, 4]$ & $-$ & $-6.6666$\\
			& \textit{P$5$} & \cite{Floudas1999book1} & $3$ & $2$ & $0$ & $2$ & $[0, 9.422] \times [0, 5.903] \times [0, 267.42]$ & $-$ & $201.1600$\\
			& \textit{P$6$} & \cite{Floudas1999book1} & $2$ & $1$ & $0$ & $1$ & $[0, 115.8] \times [10^{-5}, 30]$ & $-$ & $376.2900$\\
			& \textit{P$7$} & \cite{Floudas1999book1} & $2$ & $4$ & $0$ & $1$ & $[-2, 2]^n$ & $-$ & $-2.8284$\\
			& \textit{P$8$} & \cite{Floudas1999book1} & $2$ & $2$ & $0$ & $1$ & $[-8, 10] \times [0, 10]$ & $-$ & $-118.7000$\\
			& \textit{P$10$} & \cite{Floudas1999book1} & $2$ & $2$ & $0$ & $2$ & $[0, 1]^n$ & $-$ & $0.7417$\\
			& \textit{P$11$} & \cite{Floudas1999book1} & $2$ & $1$ & $0$ & $1$ & $[0, 1]^n$ & $-$ & $-0.5000$\\
			& \textit{P$12$} & \cite{Floudas1999book1} & $1$ & $2$ & $0$ & $0$ & $[0, 2]$ & $-$ & $-16.7390$\\
			& \textit{P$13$} & \cite{Floudas1999book1} & $3$ & $0$ & $2$ & $2$ & $[10^{-5}, 34] \times [10^{-5}, 17] \times [100, 300]$ & $-$ & $189.3500$\\
			& \textit{P$15$} & \cite{Floudas1999book1} & $3$ & $0$ & $3$ & $3$ & $[10^{-5}, 12.5] \times [10^{-5}, 37.5] \times [0, 50]$ & $-$ & $0.0000$\\
			& \textit{P$16$} & \cite{Floudas1999book1} & $2$ & $6$ & $0$ & $0$ & $[1, 3] \times [1, 4]$ & $-$ & $0.7049$\\
			& \textit{s$365$mod} & \cite{Vaz2009} & $7$ & $9$ & $0$ & $5$ & $[0, 19]^n$ & $-$ & $52.1399$\\
			& \textit{Tproblem} & \cite{Finkel2005} & $2,3,4,5,$ $6,7,8$ & $1$ & $0$ & $1$ & $[-4, 4]^n$ & $-$ & $-n$\\
			& \textit{zy$2$} & \cite{Vaz2009} & $2$ & $3$ & $0$ & $1$ & $[0, 10]^n$ & $-$ & $2.0000$\\
			& \textit{zecevic$3$} & \cite{Vaz2009} & $2$ & $2$ & $0$ & $1$ & $[0, 10]^n$ & $-$ & $97.3094$\\
			& \textit{zecevic$4$} & \cite{Vaz2009} & $4$ & $2$ & $0$ & $1$ & $[0, 10]^n$ & $-$ & $7.5575$\\

			\bottomrule

		\end{longtable}
	}

	\section{The mathematical formulation of engineering problems}
	\label{sec:apendix}

	\subsection{Tension/compression spring design problem}
	\label{prob:Tension}

	The design variables of the tension/compression spring design problem~\cite{Kazemi2011} are the number of the wire diameter $x_1$, the winding diameter $x_2$, and active coils of the spring $x_3$.
	The objective function and the mechanical constraints are given by:

	\begin{equation*}
		\label{eq:tension}
		\begin{aligned}
			\min 	& \; f(\mathbf{x})=x_1^2x_2(x_3+2)\\
			\st \;  & g_1(\mathbf{x})= 1- \frac{x_2^3x_3}{71875x_1^4} \leq 0,~g_2(\mathbf{x})= \frac{x_2(4x_2-x_1)}{12566x_1^3(x_2-x_1)}+ \frac{2.46}{12566x_1^2}-1 \leq 0,~ g_3(\mathbf{x})= 1- \frac{140.54x_1}{x_3x_2^2} \leq 0, \\
			& g_4(\mathbf{x})= \frac{x_1+x_2}{1.5}-1 \leq 0
		\end{aligned}
	\end{equation*}
	where $0.05 \leq x_1 \leq 0.2,~0.25 \leq x_2 \leq 1.3,~2 \leq x_3 \leq 15$.
	The best known solution $\mathbf{x}^* =$ ($0.05169591, 0.35688327, 11.29337893$), where $f(\mathbf{x}^*) = 0.01267867$.
	Two of the constraint functions are active ($g_1$ and $g_2$).

	\subsection{Three-bar truss design problem}
	\label{prob:Three-bar}
	The three-bar truss design problem~\cite{Ray2003} has two design variables and three constraints.
	The optimization problem is formulated as follows:
	\begin{equation*}
		\label{eq:three-bar}
		\begin{aligned}
			\min 	& \; f(\mathbf{x})=100(2 \sqrt{2} x_1+x_2)\\
			\st \;  & g_1(\mathbf{x})= \frac{\sqrt{2} x_1+x_2}{\sqrt{2} x_1^2+2x_1x_2}2-2 \leq 0,~  g_2(\mathbf{x})= \frac{x_2}{\sqrt{2} x_1^2+2x_1x_2}2-2 \leq 0,~g_3(\mathbf{x})= \frac{1}{x_1+\sqrt{2} x_2}2-2 \leq 0
		\end{aligned}
	\end{equation*}
	where $0 \leq x_1 \leq 1,~0 \leq x_2 \leq 1$.
	The best known solution $\mathbf{x}^* =$ ($0.78867531, 0.40824778$), where $f(\mathbf{x}^*) = 263.89584337$.
	One of the constraint functions is active ($g_1$).

	\subsection{NASA speed reducer design problem}
	\label{prob:nasa}

	The design variables of the NASA speed reducer design problem~\cite{Ray2003} are the face width $x_1$, the module of teeth $x_2$, the number of teeth on the pinion $x_3$, the length of the first shaft between the bearings $x_4$, the distance of the second shaft between the bearings $x_5$, the diameter of the first shaft $x_6$, and, finally, the width of the second shaft $x_7$.
	The optimization problem is formulated as follows:

	\begin{equation*}
		\label{eq:nasa}
		\begin{aligned}
			\min 	& \; f(\mathbf{x}) = 0.7854x_1x_2^2(3.3333x_3^2+14.9334x_3-43.0934) -1.508x_1(x_6^2+x_7^2)+7.4777(x_6^3+x_7^3) \\
			& +0.7854(x_4x_6^2+x_5x_7^2) \\
			\st \;  & g_1(\mathbf{x}) = \frac{27}{x_1x_2^2x_3}-1 \leq 0,~
			g_2(\mathbf{x})=\frac{397.5}{x_1x_2^2x_3^2} - 1 \leq 0,~ g_3(\mathbf{x})=\frac{1.93x_4^3}{x_2x_3x_6^4} - 1 \leq 0,~
			g_4(\mathbf{x})=\frac{1.93x_5^3}{x_2x_3x_7^4} - 1 \leq 0,\\
			& g_5(\mathbf{x})=  \frac{((\frac{745x_4}{x_2x_3})^2+16.9 \times 10^6)^{0.5}}{110x_6^3}  - 1 \leq 0,~g_6(\mathbf{x})=\frac{((\frac{745x_5}{x_2x_3})^2+157.5 \times 10^6)^{0.5}}{85x_7^3}  - 1 \leq 0, \\
			& g_7(\mathbf{x})=\frac{x_2x_3}{40} - 1 \leq 0,~
			g_8(\mathbf{x})=  \frac{5x_2}{x_1} - 1 \leq 0,~g_9(\mathbf{x})=\frac{x_1}{12x_2} - 1 \leq 0,~
			g_{10}(\mathbf{x})=\frac{1.5x_6+1.9}{x_4} - 1 \leq 0, \\
			& g_{11}(\mathbf{x})=\frac{1.1x_7+1.9}{x_5} - 1 \leq 0
		\end{aligned}
	\end{equation*}
	where $2.6 \leq x_1 \leq 3.6,~0.7 \leq x_2 \leq 0.8,~17 \leq x_3 \leq 28,~7.3 \leq x_4 \leq 8.3,~7.8 \leq x_5 \leq 8.3,~2.9 \leq x_6 \leq 3.9,~5 \leq x_7 \leq 5.5 $.
	The best known solution $\mathbf{x}^* =$ ($3.5, 0.7, 17, 7.3, 7.8, 3.35021467, 5.28668323$), where $f(\mathbf{x}^*) = 2996.34816924$.
	Three constraints are active ($g_5, g_6$ and $g_8$).

	\subsection{Pressure vessel design problem}
	\label{prob:Pressure}

	There are four design variables in the pressure vessel design problem~\cite{Kazemi2011}(in inches): the thickness of the pressure vessel $x_1$, the thickness of the head $x_2$, the inner radius of the vessel $x_3$, and the length of the cylindrical component $x_4$.
	The optimization problem is formulated as follows:
	\begin{equation*}
		\label{eq:perssure}
		\begin{aligned}
			\min 	& \; f(\mathbf{x})=0.6224x_1x_3x_4+1.7781x_2x_3^2+3.1661x_1^2x_4 +19.84x_1^2x_3\\
			\st \;  & g_1(\mathbf{x})= -x_1+0.0193x_3  \leq 0,~g_2(\mathbf{x})= -x_2+0.00954x_3 \leq 0,~g_3(\mathbf{x})=- \pi x_3^2x_4- \frac{4}{3} \pi x_3^3+1296000 \leq 0, \\
			& g_4(\mathbf{x})= x_4-240 \leq 0,~
			g_5(\mathbf{x})= 1.1-x_1 \leq 0,~g_6(\mathbf{x})= 0.6-x_2  \leq 0
		\end{aligned}
	\end{equation*}
	where $1 \leq x_1 \leq 1.375,~0.625 \leq x_2 \leq 1,~25 \leq x_3 \leq 150,~25 \leq x_4 \leq 240$.
	The best known solution $\mathbf{x}^* =$ ($1.1, 0.625, 56.99481865, 51.00125173$), where $f(\mathbf{x}^*) = 7163.73956887$.
	Three constraints are active ($g_1, g_3$ and $g_5$).

	\subsection{Welded beam design problem}
	\label{prob:Welded}

	The welded beam design problem~\cite{Mirjalili2016,Mirjalili2014} is to design a welded beam at minimum cost, subject to some constraints \cite{Mirjalili2016,Mirjalili2014}.
	The objective is to find a minimum fabrication cost.
	Considering the four design variables and constraints of shear stress $\tau$, bending stress in the beam $\sigma$, buckling load on the bar $P_c$, and end deflection on the beam $\delta$.
	The optimization model is summarized in the following equation:

	\begin{equation*}
		\label{eq:welder}
		\begin{aligned}
			\min 	& \; f(\mathbf{x}) = 1.10471 x_1^2 x_2+0.04811 x_3 x_4 (14+x_2)\\
			\st \;  & g_1(\mathbf{x}) = \tau (\mathbf{x}) - 13600 \leq 0,~g_2(\mathbf{x}) = \sigma (\mathbf{x}) - 3\times 10^4 \leq 0,~g_3(\mathbf{x}) = x_1 - x_4 \leq 0,~g_4(\mathbf{x}) = P - P_c (\mathbf{x}) \leq 0 \\
			& g_5(\mathbf{x}) = 0.10471 x_1^2 + 0.04811x_3 x_4 (14 + x_2) - 5 \leq 0,~g_6(\mathbf{x}) = \delta (\mathbf{x}) - 0.25 \leq 0,~g_7(\mathbf{x}) = 0.125 - x_1 \leq 0, \\
		\end{aligned}
	\end{equation*}
	with:
	\begin{equation*}
		\begin{aligned}
			& \tau (\mathbf{x}) = \sqrt{(\tau^1)^2+(\tau^1) (\tau^2) x_2/R + (\tau^2)^2},~P_c = \frac{4.013 E \sqrt{x_3^2 x_4"6 / 36} }{L^2}(1 - \frac{x_3}{2L} \sqrt{\frac{E}{4G}}),~R = \sqrt{\frac{x_2^2)}{4} + (\frac{x_1+x_3}{2})^2}, \\
			& \tau^1 = \frac{P}{\sqrt{2} x_1 x_"},~\tau^2 = \frac{M R}{J},~\sigma (\mathbf{x}) = \frac{6PL}{x_4 x_3^2},~J = 2(\sqrt{2}x_1 x_2 (\frac{x_2^2}{12} + \frac{1}{4}(x_1 x_3)^2)),~ \delta (\mathbf{x}) = \frac{4PL^3}{Ex_4 x_3^3},~M = P \lgroup L+ \frac{x_2}{2} \rgroup, \\
			& P = 6000, L = 14,~E = 3 \times 10^7,~G = 12 \times 10^6,
		\end{aligned}
	\end{equation*}
	where $0.1 \leq x_1 \leq 2, ~0.1 \leq x_2 \leq 10, ~0.1 \leq x_3 \leq 10, ~0.1 \leq x_3 \leq 2$.
	The best known solution $\mathbf{x}^* =$ ($0.20572963, 3.47048893, 9.03662399, 0.20572964$), where $f(\mathbf{x}^*) = 1.72485237$.
	One of the constraint functions is active ($g_3$).

	\subsection{General non-linear regression problem}
	\label{prob:general}

	Parameter estimation in the general non-linear regression model~\cite{Gillard2017,Paulavicius2019:eswa,Sergeyev2016a} can be reduced to solving the minimization problem:

	\begin{equation*}
		\label{eq:sinusoid}
		\begin{aligned}
			\min 	& \; f(\mathbf{x})= \sum_{t=1}^{T} (\kappa(t) - \phi(\mathbf{x},t))^2 \\
		\end{aligned}
	\end{equation*}
	with:
	\begin{equation*}
		\begin{aligned}
			&  \sum_{q=1}^{\varsigma} \kappa(t) = e^{t d_q} \sin(2 \pi t \omega_q + \theta_q),~\sum_{q=1}^{\varsigma} \phi(\mathbf{x},t) = e^{(x_{3(q - 1) + 1} t)} \sin(2 \pi t x_{3(q - 1) + 2} + x_{3(q - 1) + 3}) \\
		\end{aligned}
	\end{equation*}

	where $-1 \leq x_{3(q - 1) + 1} \leq 0,~0 \leq x_{3(q - 1) + 2, 3(q - 1) + 2} \leq 1, q = 1...\varsigma$.
	$\mathbf{d}$ is non-positive damping coefficients, \text{\boldmath$\omega$} is frequencies, and \text{\boldmath$\theta$} is phases of the sinusoids ($\varsigma$) (hereafter, $d_q \in [-1, 0],~\omega_q \in [0, 1],~\theta_q \in [0, 1],~ q = 1...\varsigma$).
	For signal approximation, the parameter $\mathbf{x}$ of the problem is determined to fit best the real-valued signal values observed in the uniformly distributed time moments $t = 1, 2,..., T$.
	The general non-linear regression problem is multi-modal, especially with the increase in the number of samples $T$.
	The increase of the sinusoid number $\varsigma$ leads to a more accurate but at the same time more challenging optimization problem.
  In our experimental study, six versions of the problem were considered.
  The sinusoid number was fixed to $\varsigma = 1, 2,$ and $3$ (corresponding to $3, 6,$ and $9$-dimensional cases), while the value of $T$ to $10$ and $100$.
	The best known solutions: i) $f(\mathbf{x}^*) = 0$ and $\mathbf{x}^* =$ ($-0.2$, $0.4$, $0.3$) for $n = 3$; ii) $f(\mathbf{x}^*) = 0$ and $\mathbf{x}^* =$ ($-0.3, 0.3, 0.1, -0.2, 0.4, 0.3$) for $n = 6$; iii) $f(\mathbf{x}^*) = 0$ and $\mathbf{x}^* =$ ($-0.4, 0.6, 0.2, -0.3, 0.3, 0.1, -0.2, 0.4, 0.3$) for $n = 9$.

	\section{Quick \toolbox{} User Guide}
	\label{sec:userguide}

	This section provides a brief user guide on how to use \toolbox{} software.
	The following subsections provide examples of using algorithms (their implementations) to solve box constrained and problems with various constraints, including the parallel usage of the algorithms.

	\subsection{Example of box constrained global optimization algorithm usage}
	\label{sec:boxconstrained}

	Any \direct-type algorithmic implementation from the \toolbox{} \matlab{} toolbox for box constrained global optimization can be called using the same style and syntax introduced in \cref{toolbox}.
  In this example we use the \plor{} algorithm and solve \textit{Bukin6} test problem (see \cref{tab:resultslow}).

  The \textit{Bukin6} test problem is defined in \matlab{} in the following way:

\begin{minted}[bgcolor=bg]{matlab}
function y = Bukin6(x)
    if nargin == 0                     % Extract info from the function
        y.nx = 2;                      % Dimension of the problem
        xl = [-15; -3];
        y.xl = @(i) xl(i);             % Lower bounds for each variable
        xu = [5; 3];
        y.xu = @(i) xu(i);             % Upper bounds for each variable
        y.fmin = @(nx) get_fmin(nx);   % Known solution value
        y.xmin = @(nx) get_xmin(nx);   % Known solution point
        return
    end
    term1 = 100*sqrt(abs(x(2) - 0.01*x(1)^2));
    term2 = 0.01*abs(x(1) + 10);
    y = term1 + term2;                 % Return function value at x
end

function fmin = get_fmin(~)
    fmin = 0;
end

function xmin = get_xmin(~)
    xmin = [-10; 1];
end
\end{minted}
	Each test problem in the \directlib{} stores the information about the problem structure together with the objective function.
	In this case, in the \texttt{Bukin6.m} file, the following information is stored: i) the dimensionality of the problem; ii) the lower and upper bounds for each variable; iii) the objective function value of the known solution; iv) the solution point.
	For some problems, the optimum might depend on the number of variables, therefore the solution values and points are returned as a functions for all test problems in \directlib.

  The optimization problem is passed to the algorithm as part of a \texttt{P} structure.
	For a box-constrained problem, only one field of the \texttt{P} structure is needed:
\begin{minted}[bgcolor=bg]{matlab}
>> P.f = 'Bukin6';
\end{minted}
  When a user wants to change the default algorithmic settings, the \texttt{OPTS} structure should be used:
\begin{minted}[bgcolor=bg]{matlab}
>> opts.maxevals = 50; % Maximal number of function evaluations
>> opts.maxits = 100;  % Maximal number of iterations
>> opts.testflag = 1;  % 1 if global minima is known, 0 otherwise
>> opts.tol = 0.01;    % Tolerance for termination if testflag = 1
\end{minted}
	Now we are ready to call the dynamic data structure based \plor{} implementation (\texttt{dPlor.m}) to solve this problem:
\begin{minted}[bgcolor=bg]{matlab}
>> [f_min, x_min, history] = dPlor(P, OPTS);
\end{minted}
	The iterative output stored in the \texttt{history} parameter contains the following information:
\begin{minted}[bgcolor=bg]{matlab}
>> history

history =

    1.0000     5.0000    16.7833    0.0023
    2.0000     7.0000    16.7833    0.0030
    3.0000    13.0000     5.6500    0.0039
    4.0000    19.0000     5.6500    0.0046
    5.0000    27.0000     1.9537    0.0053
    6.0000    33.0000     1.9537    0.0060
    7.0000    41.0000     0.7167    0.0070
    8.0000    47.0000     0.7167    0.0077
    9.0000    55.0000     0.3060    0.0086
\end{minted}
  Here, the first column shows the iteration number, while the second is the total number of function evaluations.
	The third column shows how the best objective function value improves at each iteration, while the last column shows the execution time in seconds.
	In this example, the \plor{} algorithm was terminated when the maximum number of function evaluations (\texttt{opts.maxevals = 50}) exceeded.

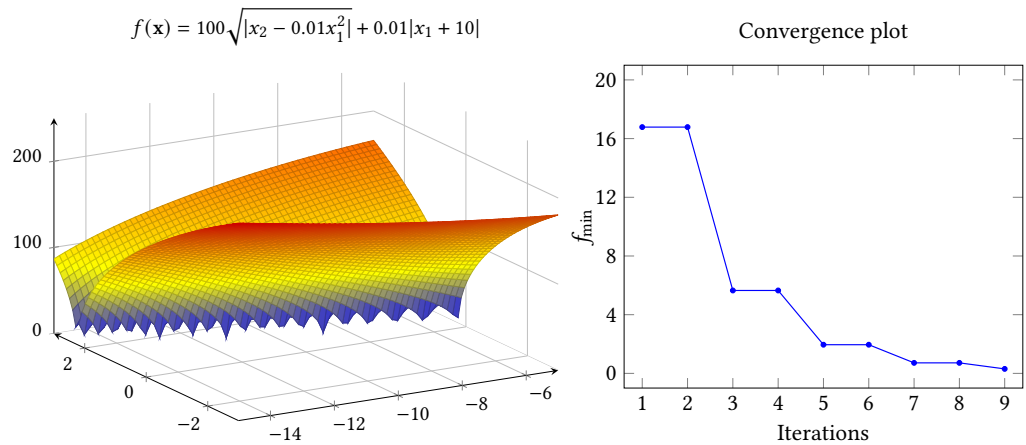
\begin{figure}[ht]
		\resizebox{.9\textwidth}{!}{
			\begin{tikzpicture}
				\begin{axis}[
					width=0.5\textwidth,
					height=0.35\textwidth,
					title = {$f(\mathbf{x}) = 100 \sqrt{| x_2 - 0.01x_1^2|} +0.01|x_1+10|$},
					scale only axis,
					ymin=-3, ymax=3,
					xmajorgrids,
					xmin=-15, xmax=-5,
					ymajorgrids,
					zmin=0, zmax=250,
					zmajorgrids,
					axis lines=left,
					grid=major,
					view={330}{25},
					]
					\addplot3[surf,samples=60,	domain=-15:-5, y domain=-3:3]
					{100*sqrt(abs(y-0.01*x^2))+0.01*abs(x+10)};
				\end{axis}
			\end{tikzpicture}
			\begin{tikzpicture}
				\begin{axis}[
					width=0.5\textwidth,height=0.425\textwidth,
					title = {Convergence plot},
					every axis/.append style={font=\large},
					ylabel style={yshift=-0.2cm},
					xlabel = {Iterations},
					ymin=0,ymax=20,
					xmin=1,xmax=9,
					xtick distance=1,
					ytick distance=4,
					ylabel = {$f_{\rm min}$},
					enlargelimits=0.05,
					]
					\addplot[blue,mark=*,mark size=1pt,line width=0.5pt] table[x=X,y=Y] {figure_data/plor.txt};
				\end{axis}
		\end{tikzpicture}}
		\caption{The \textit{Bukin6} test function (on the left side) and the convergence plot (on the right side) of the \plor{} algorithm in the first $9$ iterations.}
		\label{fig:bukin}
	\end{figure}

  The convergence plot is shown on the right panel of \cref{fig:bukin}, while the left panel illustrates the \textit{Bukin6} test function over its domain.

	\subsection{Example of constrained global optimization algorithm usage}
	\label{sec:generalconstrained}

  Any \direct-type algorithmic implementation for constrained global optimization problems can be used the same way as box-constrained problems.
  However, for constrained problems, implemented algorithms extract additional information from the function's definition, such as the number of inequality constraints, the number of equality constraints, and the constraint functions.
	Let us take the \textit{G06} problem (see \cref{tab:resultslow}) as an example, defined in the following way:
\begin{minted}[bgcolor=bg]{matlab}
function y = G06(x)
    if nargin == 0                            % Extract info from the function
        y.nx = 2;                             % Dimension of the problem
        y.ng = 2;                             % Number of g(x) constraints
        y.nh = 0;                             % Number of h(x) constraints
        xl = [13, 0];
        y.xl = @(i) xl(i);                    % Lower bounds for each variable
        y.xu = @(i) 100;
        y.fmin = @(nx) get_fmin(nx);          % Known solution value
        y.xmin = @(nx) get_xmin(nx);          % Known solution point
        y.confun = @(i) G06c(i);              % Constraint functions
    return
    end
    y = (x(1) - 10)^3 + (x(2) - 20)^3;        % Return function value at x
end

function [c, ceq] = G06c(x)
    c(1) = -(x(1) - 5)^2 - (x(2) - 5)^2 + 100;
    c(2) = (x(1) - 6)^2 + (x(2) - 5)^2 - 82.81;
    ceq = [];
end

function fmin = get_fmin(~)
    fmin = -6961.8138751273809248;
end

function xmin = get_xmin(~)
    xmin = [14.0950000002011322; 0.8429607896175201];
end
\end{minted}
	Same as in \cref{sec:boxconstrained}, the constrained optimization problem is passed to the algorithm as part of a \texttt{P} structure:
\begin{minted}[bgcolor=bg]{matlab}
>> P.f = 'G06';
\end{minted}
	Next, assume that a user wants to stop the search as soon as the known solution is within a $0.01 \%$ error.
	The \texttt{OPTS} structure should be specified as follows:
\begin{minted}[bgcolor=bg]{matlab}
>> opts.maxevals = 10000;           % Maximal number of function evaluations
>> opts.maxits = 1000;              % Maximal number of iterations
>> opts.testflag = 1;               % 1 if global minima is known, 0 otherwise
>> opts.tol = 0.01;                 % Tolerance for termination if testflag = 1
>> opts.showits = 1;                % Print iteration status
\end{minted}
	The desired solver (in this case implementation of the \directc{} algorithm) is run using:
\begin{minted}[bgcolor=bg]{matlab}
>> [f_min, x_min, history] = dDirect_GLc(P, OPTS);
\end{minted}
	Since \texttt{opts.showits = 1}, the optimization result after the each iteration is printed in the \matlab{} command window:

\begin{minted}[bgcolor=bg]{shell}
Phase_II - searching feasible point:
con viol: 2404.4400000000 fn evals: 5
con viol: 515.5511111111  fn evals: 7
...
con viol: 0.1374240038    fn evals: 123
con viol: 0.0000000000    fn evals: 159  f_min: -5612.1483164940
Phase_I - Improve feasible solution:
Iter: 1   f_min: -5886.5625227848   time(s): 0.05935   fn evals: 197
Iter: 2   f_min: -5931.8554991123   time(s): 0.06473   fn evals: 241
...
Iter: 13  f_min: -6873.0583159376   time(s): 0.13197   fn evals: 947
Iter: 14  f_min: -6901.5099081387   time(s): 0.13869   fn evals: 1027
Minima was found with Tolerance:    1
\end{minted}
	We see that the solution $f_{\rm min} = -6901.5099081387$ (within a $0.01 \%$ error) was found after $14$ iterations.

  	\subsection{Parallel algorithm usage}
  	\label{sec:parallelexample}
  	This section briefly explains how to use parallel versions of the algorithms.
    We can see which algorithms are implemented in parallel in \cref{tab:direct_classification}.
    Assume a user wishes to use parallel code for the \plor{} algorithm.
    First, a parallel implementation of the \plor{} algorithm (\texttt{parallel\_dPlor.m}) should be chosen.
  	Next, a user should specify the number of workers (computational threads).
    For parallel \plor, it is reasonable to select 2, as only two POH are selected per iteration.
    In this case, \matlab{} parallel pool size should be specified using the \texttt{parpool} command, after which the parallel algorithm should be executed:

\begin{minted}[bgcolor=bg]{matlab}
>> parpool(2);
>> [f_min, x_min, history] = parallel_dPlor(P, OPTS);
\end{minted}

  	By default, the \textit{parpool} command starts the \matlab{} pool on the local machine with one worker per physical CPU core.
    Using \texttt{parpool(2)}, we limit the number of workers to 2.
  	After this, the parallel code is executed using both workers.
  	However, it should be taken into account that creating parallel \textit{parpool} takes some time.
    Therefore, using the parallel \plor{} algorithm is inefficient in solving simple problems.
  	The use of parallel codes should address higher-dimensionality, more expensive optimization problems~\cite{Stripinis2020}.
  	When all necessary calculations in parallel mode are finished, the following command:

\begin{minted}[bgcolor=bg]{matlab}
>> delete(gcp);
\end{minted}
    shuts down the parallel pool.

\end{document}